\documentclass[12pt,a4paper]{article}

\usepackage{amsgen,amsmath,amstext,amsbsy,amsopn,amsfonts,amssymb,latexsym}

\usepackage{graphicx}

\renewcommand{\tilde}{\widetilde}
\newcommand{\C}{{\mathbb C}}
\newcommand{\HH}{{\mathbb H}}

\newcommand{\N}{{\mathbb N}}

\newcommand{\R}{{\mathbb R}}
\newcommand{\bbS}{{\mathbb S}}

\newcommand{\Z}{{\mathbb Z}}

\newcommand{\g}{{\mathfrak g}}
\newcommand{\h}{{\mathfrak h}}

\newcommand{\twomatrix}[4]{\left(\begin{array}{cc} {#1} & {#2} \\ {#3} & {#4} \end{array}\right)}

\newcommand{\twovec}[2]{\left(\begin{array}{c} {#1}\\{#2}\end{array}\right)}

\newcommand{\setof}[1]{ \left\{ {#1} \right\} }

\newcommand{\nonvoid}{\neq\emptyset}
\newcommand{\1}{{\bf1}}

\newcommand{\Aut}{\mathop{\rm Aut}\nolimits}

\newcommand{\spec}{\mathop{\rm spec}\nolimits}

\newcommand{\GL}{\mathop{\rm GL}\nolimits}

\newcommand{\sL}{\mathop{\mathfrak{sl}}\nolimits}

\newcommand{\so}{\mathop{\mathfrak{so}}\nolimits}

\newcommand{\su}{\mathop{\mathfrak{su}}\nolimits}

\newcommand{\SL}{\mathop{\rm SL}\nolimits}

\newcommand{\SLt}{\mathop{\tilde{\rm SL}}\nolimits(2,\R)}
\newcommand{\SO}{\mathop{\rm SO}\nolimits}
\newcommand{\SU}{\mathop{\rm SU}\nolimits}

\newcommand{\PSL}{\mathop{\rm PSL}\nolimits}

\newcommand{\aff}{\mathop{\rm aff}\nolimits}

\newcommand{\Diffeo}{\mathop{\rm Diffeo}\nolimits}
\newcommand{\id}{\mathop{\rm id}\nolimits}

\newcommand{\inter}{\mathop{\rm int}\nolimits}
\newcommand{\dom}{\mathop{\rm dom}\nolimits}
\newcommand{\ad}{\mathop{\rm ad}\nolimits}
\newcommand{\CoAd}{\mathop{\rm CoAd}\nolimits}

\newcommand{\arsinh}{\mathop{\rm arsinh}\nolimits}
\newcommand{\arcosh}{\mathop{\rm arcosh}\nolimits}
\newcommand{\artanh}{\mathop{\rm artanh}\nolimits}
\newcommand{\Ad}{\mathop{\rm Ad}\nolimits}

\newcommand{\sign}{\mathop{\rm sign}\nolimits}
\newcommand{\diag}{\mathop{\rm diag}\nolimits}

\newcommand{\trace}{\mathop{\rm trace}\nolimits}

\newcommand{\conv}{\mathop{\rm conv}\nolimits}

\newcommand{\norm}[1]{\left\|{#1}\right\|}

\newcommand{\qaq}{\quad\mbox{and}\quad}

\newtheorem{Theorem}{Theorem}[section]
\newtheorem{Proposition}[Theorem]{Proposition}
\newtheorem{Corollary}[Theorem]{Corollary}

      \newtheorem{rem}[Theorem]{Remark}
\newenvironment{Remark}{\begin{rem}\rm}{\end{rem}}
      \newtheorem{define}[Theorem]{Definition} 

      \newtheorem{exer}[Theorem]{Exercise}

      \newtheorem{ex}[Theorem]{Example}
\newenvironment{Example}{\begin{ex}\rm}{\end{ex}}
      \newtheorem{exs}[Theorem]{Examples}

      \newtheorem{conj}[Theorem]{Conjecture}
\newenvironment{Conjecture}{\begin{conj}\rm}{\end{conj}}
      \newtheorem{probl}[Theorem]{Problem}

      \newtheorem{prs}[Theorem]{Problems}

\newenvironment{Theorem*}[1]{\begin{trivlist}\item[\hskip%
    \labelsep{\bf Theorem #1.\quad}]\em}%
    {\rm\end{trivlist}}
\newenvironment{Proposition*}[1]{\begin{trivlist}\item[\hskip%
    \labelsep{\bf Proposition #1.\quad}]\em}%
    {\rm\end{trivlist}}
\newenvironment{Corollary*}[1]{\begin{trivlist}\item[\hskip%
    \labelsep{\bf Corollary #1.\quad}]\em}%
    {\rm\end{trivlist}}
\newenvironment{Lemma*}[1]{\begin{trivlist}\item[\hskip%
    \labelsep{\bf Lemma #1.\quad}]\em}%
    {\rm\end{trivlist}}
\newenvironment{Definition*}[1]{\begin{trivlist}\item[\hskip%
    \labelsep{\bf Definition #1.\quad}]}%
    {\end{trivlist}}
\newenvironment{Remark*}[1]{\begin{trivlist}\item[\hskip%
    \labelsep{\bf Remark #1.\quad}]}%
    {\end{trivlist}}
\newenvironment{Exercise*}[1]{\begin{trivlist}\item[\hskip%
    \labelsep{\bf Exercise #1.\quad}]}%
    {\end{trivlist}}
\newenvironment{Example*}[1]{\begin{trivlist}\item[\hskip%
    \labelsep{\bf Example #1.\quad}]}%
    {\end{trivlist}}
\newenvironment{Examples*}[1]{\begin{trivlist}\item[\hskip%
    \labelsep{\bf Examples #1.\quad}]}%
    {\end{trivlist}}
\newenvironment{Conjecture*}[1]{\begin{trivlist}\item[\hskip%
    \labelsep{\bf Conjecture #1.\quad}]}%
    {\end{trivlist}}
\newenvironment{Problem*}[1]{\begin{trivlist}\item[\hskip%
    \labelsep{\bf Problem #1.\quad}]}%
    {\end{trivlist}}
\newenvironment{Problems*}[1]{\begin{trivlist}\item[\hskip%
    \labelsep{\bf Problems #1.\quad}]}%
    {\end{trivlist}}

\newenvironment{Proof}{\par\noindent{\bf
Proof.\quad}}{\rule{.1cm}{0cm}\hfill$\Box$\par\bigskip} 
\newenvironment{Proof*}{\par\noindent{\bf Proof.\quad}}{\par\bigskip}

\newcommand{\smo}{\setminus\{0\}}
\newcommand{\qmb}[1]{\quad\mbox{#1}}
\newcommand{\qmbq}[1]{\quad\mbox{#1}\quad}

{  \end{center}
  \end{minipage}\\ \hline
 \end{tabular}
\end{center}
}

\newcommand{\nobox}[1]{#1}

\DeclareMathOperator{\mat}{mat}

\DeclareMathOperator{\Diss}{Diss}
\DeclareMathOperator{\Fix}{Fix}
\DeclareMathOperator{\opp}{op}

\newcommand{\tra}{\mathsf{T}}
\newcommand{\wh}{\widehat}
\newcommand{\wt}{\widetilde}
\newcommand{\tg}{{\tilde{g}}}

\newcommand{\seto}[1]{\{#1\}}
\newcommand{\Conj}[1]{\mathop{\rm I}\nolimits_{#1}}

\newcommand{\cover}{\mathop{f}\nolimits}      
\newcommand{\icover}{\mathop{f^{-1}}\nolimits}
\newcommand{\expt}{\mathop{\rm exp}\nolimits} 
\newcommand{\tcU}{\tilde{\mathcal{U}}}

\newcommand{\ggen}[1]{\left\langle \left\langle #1\right\rangle\right\rangle}
\newcommand{\cD}{\mathcal{D}}
\newcommand{\cE}{\mathcal{E}}
\newcommand{\cF}{\mathcal{F}}
\newcommand{\cH}{\mathcal{H}}
\newcommand{\cM}{\mathcal{M}}
\newcommand{\cQ}{\mathcal{Q}}
\newcommand{\cT}{\mathcal{T}}
\newcommand{\cU}{\mathcal{U}}

\newcommand{\sZ}{Z}

\newcommand{\bbD}{\mathbb{D}}
\newcommand{\bbF}{\mathbb{F}}
\newcommand{\bbbF}{\overline{\bbF}}

\newcommand{\bbM}{\mathbb{M}}

\newcommand{\bfF}{\mathbf{F}}
\newcommand{\bfFel}{\mathbf{F_{el}}}
\newcommand{\bfFpl}{\mathbf{F_{pl}}}
\newcommand{\bfR}{\mathbf{R}}

\newcommand{\bfg}{{\boldsymbol{g}}}
\newcommand{\bfz}{{\boldsymbol{z}}}

\newcommand{\bfgamma}{{\boldsymbol{\gamma}}}
\newcommand{\bfell}{{\boldsymbol{\ell}}}
\newcommand{\bfvarphi}{{\boldsymbol{\varphi}}}

\newcommand{\rmD}{\mathrm{D}}
\newcommand{\rmd}{\mathrm{d}}

\newcommand{\dst}{\displaystyle}

\newcommand{\mex}[1]{\mathop{\mathbb{M}}\left( #1 \right)}

\newcommand{\exP}[1]{
  \begin{pmatrix}
    1 & #1 \\ 0 & 1
  \end{pmatrix}
}

\newcommand{\exQ}[1]{
  \begin{pmatrix}
    1 & 0 \\ #1 & 1
  \end{pmatrix}
}

\newcommand{\exU}[1]{\begin{pmatrix} \cos({#1}) & \sin({#1}) \\
 -\sin({#1}) & \cos({#1}) \end{pmatrix}}

\newcommand{\matH}{\begin{pmatrix} 1 & 0 \\  0 & -1 \end{pmatrix}}
\newcommand{\matT}{\begin{pmatrix} 0 & 1 \\  1 &  0 \end{pmatrix}}
\newcommand{\matU}{\begin{pmatrix} 0 & 1 \\ -1 &  0 \end{pmatrix}}

\newcommand{\matI}{\begin{pmatrix} 1 & 0 \\  0 &  1 \end{pmatrix}}

\begin{document}

\title{The dissipation distance for a 2D single crystal with two symmetric slip systems}

\author{Dirk Mittenhuber\\
Math. Inst. A, Uni Stuttgart\\
 \texttt{Dirk.Mittenhuber@mathematik.uni-stuttgart.de}}


\maketitle

\begin{abstract}
  We solve a model problem from single crystal plasticity.
  We consider $4$ slip systems in the plane with orthogonal slip-directions
  and equal slip rates, forward as well as backwards.
  We compute the associated dissipation distance by solving
  an optimal control problem. It turns out that from a computational point of
  view computing the distance is inexpensive.
  We put special emphasis on visualization of the metric spheres and
  the associated length-minimizing curves.\\
  As a byproduct we also solve a related problem, optimal path planning for a car driving
  forwards and backwards with limited turning radius in the hyperbolic plane.
  This is a hyperbolic version of the Reeds-Shepp-Car-Problem first
  discussed in~\cite{reeds-shepp}.
\end{abstract}


\section{Introduction}

In this paper we provide the solution to an optimization problem which has
various interpretations. Although we will put special emphasis on the one
mentioned in the title, the underlying mathematical problem does not require
any knowledge from continuum mechanics and/or finite plasticity. Therefore we
start with a naive formulation as an optimal factorization problem in the group
$\SL(2,\R)$ of invertible $2$~by~$2$-matrices with determinant~$1$.

\subsection*{A factorization problem}

Let $\sL(2)$ denote the set of
$2$~by~$2$-matrices with zero trace, and let
$$
  P=\begin{pmatrix}0 & 1\cr 0 & 0 \end{pmatrix},\qaq
  Q=\begin{pmatrix}0 & 0\cr 1 & 0 \end{pmatrix}.
$$
 So $P$ and $Q$ are generators for shearings along the coordinate axes in~$\R^2$.
 Then it is well-known that every $g\in\SL(2,\R)$ may be written as a product
$$
  g=\exp(t_1A_1)\cdots\exp(t_kA_k)\quad
 \mbox{with $k\in\N$, $t_k\in\R$, and $A_k\in\{P,Q\}$.}
$$
 We want to find factorization(s) of a given~$g$
 such that $\sum_i |t_i|$ is minimal.
 Therefore we define the {\em factorization cost}
 $\mathcal{T}(g)$ as
$$
  \mathcal{T}(g)=\inf\setof{\sum_i|t_i| : g=\exp(t_1A_1)\dots\exp(t_kA_k),\ 
    \begin{array}[t]{l}
      k\in\N,\ t_k\in\R,\\  A_k\in\{P,Q\}
    \end{array}
  }.
$$
The factorization cost $\mathcal{T}(g)$ can be interpreted as the distance of~$g$ from the
identity matrix. It can also be used to measure distances
in the group $\SL(2)$:
$$
\wh D(g_0,g_1)\stackrel{\rm def}{=\!=}\mathcal{T}\left( g_0^{-1}g_1 \right), \qmbq{or}
\check{D}(g_0,g_1)\stackrel{\rm def}{=\!=}\mathcal{T}\left( g_0g_1^{-1}\right).
$$
It turns out that $\wh D$, $\check{D}$ are metrics on~$\SL(2)$, by construction
$\wh D$ is left-invariant while $\check{D}$ is right-invariant:
$$
\wh D(g_0,g_1)=\wh D(gg_0\,,gg_1),\quad \check{D}(g_0,g_1)=\check{D}(g_0g\,,g_1g)
\qmbq{for all} g,g_0,g_1\in\SL(2).
$$
In the sequel we will solve the problem of computing $\mathcal{T}(g)$
through an associated optimal control problem.  Our technique is kind of standard in
control theory in the sense that we use the Pontrjagin Maximum Principle (PMP)
as a necessary condition for optimality plus some adhoc comparison arguments.
We will also point out how the (PMP) relates to the \textbf{yield surface} and
\textbf{flow rule} used in the plasticity literature.

\subsection*{Sneak preview}

In order to give the reader an idea of the final outcome (and the
computational complexity) we state a few consequences of our final results.

\begin{Theorem}\label{six-factors-preview-thm}
   Every $g\in\SL(2)$ has an optimal factorization of the form
$$
  g=\exp(t_1A_1)\dots\exp(t_6A_6),\qmbq{with} A_k\in\setof{P, Q, P+Q},\ t_k\in\R.
$$
\end{Theorem}
So no more than $6$ factors are needed, but it is necessary to allow
factors of the form~$\exp(\pm t(P+Q))$, too.
Otherwise $\mathcal{T}(g)$ will be an infimum for some matrices $g\in\SL(2)$,
in particular for $g=\exp(t(P+Q))$ with $t\neq0$.

As a consequence, finding an optimal factorization is reduced to a finite
problem. We will actually obtain the following, much more detailed
information:

\begin{Theorem}\label{family-64-thm}
   There exists a sufficient family consisting of~$64$ maps, i.e., there exist
 functions $f_1,\dots,f_{64}\colon\R^3 \to\SL(2)$ with the following property:
\linebreak For every $g\in\SL(2)$
 there exist $k\in\setof{1,\dots,64}$ and $r,s,t\geq0$ such that $f_k(r,s,t)$
 provides an optimal factorization of~$g$.
\end{Theorem}

If one is really interested in computing an optimal factorization explicitly
one can exploit the symmetry of the problem and reduce the number of maps that
have to be inverted to~$13$ (rather than $64$).\label{64-13-notice}
For efficient computation of the function $\mathcal{T}$ one can even
reduce this to~$12$ maps.

To give a rough idea of the computational complexity we note that except for one map,
nothing worse than solving quadratic equations is required. In this
\textit{worst case} the challenge consists of solving a cubic equation $p(x)=y$,
and this  needs to be done only over an $x$-interval where the underlying cubic
polynomial $p$ is strictly increasing, convex, and, $p'(x)$ is bounded away
from~$0$.

\smallskip
Finally, we will show that for any  other pair of rank-$1$ matrices
$S^1,S^2\in\sL(2)$ the solution of the associated factorization problem
for $A_k\in\setof{\pm S^1,\pm S^2}$ can be obtained from $\mathcal{T}(g)$
in the following, very simple way:
\begin{Theorem}\label{other-2slip-thm}
  Let $S^1,S^2\in\sL(2)$ with $\det(S^1)=\det(S^2)=0$, $[S^1,S^2]\neq0$.
  Let $\mathcal{T}_S$ denote the factorization cost for $\setof{\pm S^1,\pm S^2}$.
  Then there exist $\lambda>0$ and an automorphism $\sigma\colon\SL(2)\to\SL(2)$
  such that $\mathcal{T}_S(g)=\lambda \mathcal{T}(\sigma(g))$.
\end{Theorem}

We will also show how $\lambda>0$ and $\sigma\in\Aut(\SL(2))$ are obtained, given
$S^1,S^2$. Thus we have determined the dissipation distance for 
\textbf{any $2$-slip system} with \textbf{symmetric dissipation functonial}.


\subsection*{A reader's guide}

This paper serves several purposes, therefore a few remarks concerning these
seem to be in order.

The main purpose is to illustrate the application of optimal control
techniques and Lie group methods to finite plasticity. So partly this paper is
intended as a tutorial for  non-specialists in optimal control on
Lie groups. Therefore we will discuss everything in great detail and provide
rigorous proofs. In this spirit this report is a successor of Sussmann's and Tang's
paper~\cite{suss-car} on the Reeds-Shepp-Car-Problem. As our factorization
problem is related to the Reeds-Shepp-Car-Problem in the hyperbolic plane, our
arguments and results will bear some strong resemblance with those
in~\cite{suss-car}. Therefore we would like to stress that in this paper we
put special emphasis on how to exploit the Lie group structure of
$\SL(2)$. The latter is instrumental in reducing the complexity and
streamlining the discussion. Moreover, it is indispensable if one's aim is to
treat similar problems in $\SL(2)$ and, eventually, in $\SL(3)$.

As the solution of the hyperbolic Reeds-Shepp-Car-Problem requires only little
extra effort, we will provide it in an appendix. Although the result
resembles that for the euclidean case, some aspects are different. For the
geometer these are properties that distinguish hyperbolic from
 euclidean geometry. The interpretation as a path planning problem in the
hyperbolic plane also provides a good visualization tool.
It is noteworthy to mention that even if one does not care about hyperbolic
geometry, one can benefit from it because some of the adhoc arguments suddenly
have a simple interpretation---they might seem perfectly obscure and
unmotivated, otherwise.

Visualization of the metric spheres (i.e., level sets of the factorization
cost $\mathcal{T}$) is another issue we deal with. Since the group $\SL(2)$
is three-dimensional, everything can be visualized in $\R^3$, but how? We will
use a parametrization coming from a polar decomposition, first proposed by
Hilgert and Hofmann in~\cite{hofmann-old-new}. As a set, $\SL(2)$
is identified with $\R^2\times[-\pi,\pi)\subseteq \R^3$, and $\R^3$ is
identified with the simply connected Lie group with Lie algebra~$\sL(2)$.  An
advantage of this parametrization is that it immediately allows to recognize
the symmetry inherent to the problem. A disadvantage is that the group
operation is more complicated than matrix multiplication. The purpose of
the first appendix is to collect information about this parametrization which
is scattered around in the literature. This information is not necessary to
understand and interprete the pictures of the metric spheres, but it is
indispensible for generating them.

The remainder of the paper is organized as follows:
\begin{description}
\item[2. From finite plasticity to Lie groups]
\ \\  Brief outline how plasticity leads to consider metrics on Lie groups.
\item[3. The underlying optimal control problem]
\item[4. Symmetries and isometries]
\item[5. The structure of extremals]
\ \\  Discussion of the (PMP), description of yield surface and flow rule.
\item[6. A sufficent family for \boldmath{$\SL(2)$}]
\ \\  Summary and short discussion of how to find best factorizations.
\item[7. Comparison arguments]
\ \\  Rigorous proofs for the sufficiency of the familiy described in Section~6.
\item[8. Conclusion]
\ \\ Brief outlook on future work and how to treat similar problems.
\item[Appendix A: Parametrizing the simply connected group \boldmath{$\SLt$}]
  All information necessary to generate the graphics.
\item[Appendix B: The hyperbolic Reeds-Shepp-Car]
\ \\  Missing arguments and comparison with the results in~\cite{suss-car}.
\item[Appendix C: More details for \boldmath{$\SLt$}]
\ \\  Additional information clarifying some of the arguments given in Section~7.
\end{description}

\noindent\textbf{Notation.} As we will have to write products of exponentials
repeatedly, we need  a shorthand notation. Let $G$ be a Lie group with Lie
algebra~$\g$ and exponential function $\exp\colon\g\to G$. Then we define
$$
\mex{X_1,\dots,X_k}:=\exp(X_1)\cdots\exp(X_k),\quad k\in\N,\ X_1,\dots,X_k\in\g.
$$
Thus $\mex{\cdot}$ is a map from $\bigcup_{k\in\N}\g^k$ to $G$.
The map $\mex{\cdot}$ depends, of course, on the group $G$,
so one should write $\mathbb{M}_G(X_1,\dots)$. But except for a few situations
in the appendix it will always be clear from the context in which group we are
working, so we omit the subscript~$G$ most of the time.

\medskip
\noindent\textbf{Acknowledgement.}  This research was supported by DFG within
\textit{SFB 404 Multifield Problems.} The author is grateful to Alexander
Mielke and Klaus Hackl for stimulating discussions.


\section{From finite plasticity to dissipation distances on Lie groups}
\nocite{mielke-geodesics-sld,mielke-sfb-dissipation}

The idea to use left-invariant metrics on Lie groups 
in the modelling of elasto-plastic material behavior is due to  Mielke,
cf.~\cite{mielke-geodesics-sld}. For a detailed overview of this approach
we refer to~\cite{mielke-sfb-dissipation}.
Let us quickly outline some of the main ideas of this approach.

\subsection*{A global formulation of elastoplasticity}

Consider a body $\Omega\subseteq\R^d$ that undergoes a deformation
$\bfvarphi:\Omega\to\R^d$. Let $\bfF=\rm\rmD\bfvarphi$ denote the deformation
gradient. Inelastic material behavior is described by an internal state~$\bfz$
from some set $\sZ$. The whole material model is based on two scalar
constitutive functions,
the elastic potential~$\wh\psi$ and the dissipation potential~$\wh\Delta$.
These give rise to an elastic storage energy and a dissipation functional.
Considering the evolution $(\bfvarphi(t),\bfz(t))$ under the influence
 of some time-varying external forces,
the total elastic and potential energy (or Gibb's energy) at time~$t$ is
$$
\cE(t,\bfvarphi,\bfz)=\int_\Omega
\wh\psi(x,\rmD\bfvarphi,\bfz)\,\rmd x-\left\langle\bfell(t),
 \bfvarphi\right\rangle,
$$
the second term corresponding to the work by external forces.
The dissipation~$\wh\Delta$ is supposed to depend only on the evolution of the internal
state $\bfz(t)$, i.e., $\wh\Delta=\wh\Delta(x,\bfz,\dot{\bfz})$. One defines 
the dissipation distance $\wh D$ as
$$
\wh D(x,\bfz_0,\bfz_1)=
\inf\setof{ \int_0^1 \wh\Delta(x,\bfz(s),\dot{\bfz}(s))\,\rmd s \mid
\begin{array}[t]{l}
\bfz(\cdot)\in C^1([0,1],\sZ),\\ \bfz(0)=\bfz_0,\ \bfz(1)=\bfz_1
\end{array}}.
$$
Integrating over~$\Omega$ one defines $\cD(\bfz_0,\bfz_1)=\int_\Omega\wh
D(x,\bfz_0,\bfz_1)\,\rmd x$. Finally  the total dissipation along a
path $\bfz(t)$ is defined as
$$
\Diss(\bfz;[t_1,t_2])=\sup\setof{ 
\sum_{j=1}^n \cD(\bfz(\tau_j),\bfz(\tau_{j-1})) \mid
t_1=\tau_0<\cdots<\tau_n=t_2
}
$$
With these functionals one obtains a notion of solution processes without making any
differentiability assumptions. A process $(\bfvarphi(t),\bfz(t))$ is called
a solution process over~$[0,T]$ if it satisfies the following two conditions:
\\[2mm]
\textbf{(S) Stability:} $\cE(t,\bfvarphi,\bfz)\leq
\cE(t,\wt\bfvarphi,\wt\bfz)+\cD(\bfz,\wt\bfz)$ for all $t\in[0,T]$ and all
 comparison states $(\wt\bfvarphi,\wt\bfz)$;
\\[2mm]
\textbf{(E) Energy inequality:}
$$\cE(t_1,\bfvarphi(t_1),\bfz(t_1))+\Diss(\bfz;[t_1,t_2])\leq
\cE(t_2,\bfvarphi(t_2),\bfz(t_2))
-\int_{t_1}^{t_2} \left\langle \dot{\bfell}(s),\bfvarphi(s)\right\rangle\,\rmd s.
 $$
This formulation does not involve any derivatives, neither of $\wh\psi,\wh\Delta$
nor of $\rmD\bfvarphi,\bfz$. As is shown in~\cite{mielke-sfb-dissipation} this
formulation is consistent with classical elasto-plastic flow rules if the
solution is sufficiently smooth.
 A particular advantage of this global formulation is that it allows to derive
incremental time-stepping algorithms which are robust.

\subsection*{Multiplicative elastoplasticity: constitutive laws}

So far we outlined the general approach without making any assumptions
on the internal state space~$\sZ$. Multiplicative elastoplasticity uses
the splitting $\rmD\bfvarphi=\bfF=\bfFel\,\bfFpl$ and considers $\bfFpl$ as an
internal variable while the elastic potential $\wh\psi$ is supposed to depend
only on $\bfFel=\bfF\,\bfFpl^{-1}$.
Actually  $\bfz=\bfFpl^{-1}$ is used as internal state, and the
internal state space~$\sZ$ is a connected Lie subgroup, say~$G$, of $\GL(d)$. Typically,
$G=\sZ=\SL(d)$, but other groups may be considered, too.
The following constitutive laws are postulated:
\begin{description}
\item[(Sy1) Objectivity:] (frame indifference)
 $\wh\psi(x,\bfR\,\bfF,\bfz)=\wh\psi(x,\bfF,\bfz)$ for all $\bfR\in\SO(d)$;
\item[(Sy2) Plastic indifference:]
$\wh\psi(x,\bfF\bfg^{-1},\bfg\bfz)=\wh\psi(x,\bfF,\bfz)$,\\
$\wh\Delta(x,\bfg\bfz,\bfg\dot{\bfz})=\wh\Delta(x,\bfz,\dot{\bfz})$ for all $\bfg\in G$;
\item[(Sy3) Rate independence:]
 $\wh\Delta(x,\bfz,\alpha\dot{\bfz})=\alpha\wh\Delta(x,\bfz,\dot{\bfz})$ for $\alpha\geq0$;
\end{description}
Material symmetries may be captured, for example, 
in a group $S\subseteq\mathrm{O}(d)\cap G$.
Following the notation in~\cite{mielke-sfb-dissipation} we postulate this as
constitutive law, too:
\begin{description}
\item[(Sy4) Material symmetry:]
$\wh\psi(x,\bfF,\bfz\bfgamma)=\wh\psi(x,\bfF,\bfz)$,\\
$\wh\Delta(x,\bfz\bfgamma,\dot{\bfz}\bfgamma)=\wh\Delta(x,\bfz,\dot{\bfz})$
for all $\bfgamma\in S$.
\end{description}

Property~(Sy2) implies that the dissipation distance $\wh D$ defined in
the previous subsection is invariant under left-multiplication with elements
from~$G$, hence
$$
\wh D(\bfz_0,\bfz_1)=\wh D(\1,\bfz_0^{-1}\bfz_1)=:\wt D(\bfz_0^{-1}\bfz_1).
$$
These metrics are the objects we want to study. Here we dropped the material
point~$x$ for sake of simplicity. In the sequel we will never consider
dependency on~$x$. This does not necessarily mean that our considerations are
limited to homogeneous media. There are suitable formulations where all
considerations are first limited to a fixed  material point~$x$, and the final result
is obtained by integration over~$\Omega$, cf.~\cite{mielke-sfb-dissipation}

\subsection*{Dissipation distances on Lie groups and time optimal control problems}

From now on we will assume that the internal state space is a connected Lie group $G$
with Lie algebra~$\g$. Therefore we slightly change the notation.
Henceforth, we write $g\in G$ (instead of~$\bfz\in\sZ$).
Our next goal is to discuss the consequences of the constitutive laws~(Sy2)
and~(Sy3). By~(Sy2) the distance function $\wh D\colon G\times G\to[0,\infty]$
is left-invariant. For the dissipation potential $\wh\Delta$ this means that
$\wh\Delta(g,\dot g)=\wt\Delta(g^{-1}\dot g)$ for some $\wt\Delta\colon\g\to[0,\infty]$.
Thus, given $\wt\Delta$ our goal is to analyze the function $\wt D\colon G\to[0,\infty]$,
$$
\wt D(g_0)=\inf\setof{ \int_0^1 \wt\Delta(g^{-1}(t)\dot g(t))\,\rmd t : 
g\in C^1([0,1],G),\ g(0)=\1,\ g(1)=g_0
}.
$$
Rate independence~(Sy3) implies $\wt\Delta(\alpha X)=\alpha\wt\Delta(X)$ for
all $\alpha\geq0$, $X\in\g$. Therefore the definition of~$\wt D$ still
contains redundancy. We can reparametrize curves by their $\wt\Delta$-length.

Indeed, assume~(Sy3) and $\wt\Delta(X)>0$ for all
$X\neq0$.  Now suppose that  $g\in C^1([0,1],G)$ is given.
Let $L(t)=\int_0^t\wt\Delta(g^{-1}(t)\dot g(t))\,dt$ and set $L_1=L(1)$.
We only need to consider finite length, so $L_1<\infty$. 
Now $L(t)$ is differentiable, and $L'(t)\geq0$ in~$[0,1]$.

If $\dot g(t)\neq0$ in~$[0,1]$, then $L'(t)>0$, and $L$ has a differentiable
inverse~$L^{-1}$. An elementary computation shows that
$g\circ L^{-1}:[0,L_1]\to G$ is parametrized by $\wt\Delta$-length,
hence for $\gamma(t)=g(L^{-1}(t L_1))$ we obtain $\gamma(0)=g(0)$,
$\gamma(1)=g(1)$, and $\wt\Delta(\gamma^{-1}\dot\gamma)\equiv L_1$.

In the general case ($\dot g(t)=0$ is possible), $L(t)$ is only monotone
 increasing. In that case one uses 
$\check L(s)=\sup\setof{t:L(t)\leq s}$ and shows that $g\circ\check L$ is
 differentiable (although $\check L$ need  not be differentiable).
Hence $\wt D$ can be characterized in the following ways
\begin{eqnarray*}
\lefteqn{\wt D(g_0)=\inf\setof{\int_0^1 \wt\Delta(g^{-1}\dot g)\,\rmd t:
  \begin{array}[t]{l}
g\in C^1([0,1];G),\ \wt\Delta(g^{-1}\dot g)\equiv\mathrm{const}
\\
 g(0)=\1,\ g(1)=g_0
  \end{array}
}.}
\\
& = & \inf\setof{T: (\exists g\in C^1([0,T];G))
\ g(0)=\1,\ g(T)=g_0,\ \wt\Delta(g^{-1}\dot g)\equiv1
  }
\\
& = & \inf\setof{T: (\exists g\in C^1([0,T];G))
\ g(0)=\1,\ g(T)=g_0,\ \wt\Delta(g^{-1}\dot g)\leq1
  }.
\end{eqnarray*}
Now let $\cU=\seto{X\in\g:\wt\Delta(X)\leq 1}$. Then the last statement says
we must
look for solutions of the differential inclusion $g^{-1}\dot g\in\cU$, with
boundary data $g(0)=\1$, $g(T)=g_0$ such that $T$ is minimal.
Thus computing the dissipation distance $\wt D$ is equivalent to solving
a time-optimal left-invariant control problem on the Lie group~$G$.
Such problems are well-studied and standard results are available, 
cf.~\cite{jurd-elastica,jurd-book,dimi-oberwolf}
\begin{Theorem}
  Assume that $\cU$ is compact and convex and  $\wt D(g_0)<\infty$.
  Then there exists an absolutely continuous $g\colon[0,\wt D(g_0)]\to G$
  such that $g^{-1}\dot g\in\cU$ a.e., $g(0)=\1$, and $g(\wt D(g_0) )=g_0$.
\end{Theorem}
Thus length minimizing arcs always exist within the class of absolutely continuous
paths, provided there is some path with finite length and
the set $\cU\subseteq\g$ is compact and convex.
If $\cU$ is a $0$-neighborhood then it is clear that $\wt D(g_0)<\infty$ for
all $g_0\in G$. But this is still true under much weaker hypotheses.
For example, let $\ggen{\cU}$ denote the
smallest subalgebra of~$\g$ containing $\cU$.
Then we have:
\begin{Theorem}
  Assume that $\cU=-\cU$ and $\ggen{\cU}=\g$. 
  Then $\wt D(g_0)<\infty$ for all $g_0\in G$.
\end{Theorem}
In control theory language: every $g_0\in G$ is reachable from the group
identity~$\1$ along a trajectory of $g^{-1}\dot g\in\cU$, see~\cite[Theorem~1]{brockett-system-theory}
or~\cite[Thm.~5.1]{susscontrol},
for example. The condition $\ggen{\cU}=\g$ is also necessary because
the set $\seto{g_0\in G: \wt D(g_0)<\infty}$ is nothing but the reachable set (from~$\1$)
of the system $g^{-1}\dot g\in\cU$. And if $\ggen{\cU}\neq\g$, this
reachable set is contained in a proper subgroup of~$G$.

Still $\wt D(g_0)<\infty$ for all $g_0\in G$ may hold true under much weaker
hypotheses.
In fact, let  $S(\cU):=\left\langle \exp\R^+\cU \right\rangle$ denote the
subsemigroup of~$G$ generated by $\exp(\R^+\cU)$.
Then $\wt D(g_0)<\infty$ for all $g_0\in S(\cU)$, and
$$
\wt D(g_0)<\infty\mbox{ for all $g_0\in G$ } \iff S(\cU)=G.
$$
Actually, $S(\cU)=G$  may hold under extremely weak assumptions.
To give just one more example, consider $\g=\sL(2)$ and
$\cU=\conv(0,P,-Q)=[0,1]\conv(P,-Q)$. Then $S(\cU)=\SL(2)$ holds true.

\subsection*{Control systems on Lie groups and automorphisms}

Given a set $\cU\subseteq\g$ we now consider the left-invariant control system
given by the differential inclusion $g^{-1}\dot g\in\cU$ a.e., and analyze
some of its properties. Left-invariance means that for a trajectory $g(t)$ and
an arbitrary $g_0\in G$ the path $\tilde g(t):=g_0g(t)$ is a trajectory, too.

We now write $\cT_\cU(g_0)$ instead of $\wt D(g_0)$ because we want to consider
various possibilities for~$\cU$. Our first observation is that
\begin{equation}\label{eq:t-mu-U}
\cT_{\mu\cU}(g_0)=\frac1\mu\,\, \cT_{\cU}(g_0)
\mbox{ for all $\cU\subseteq\g$, $\mu>0$, and $g_0\in G$.}
\end{equation}
In fact, this is obtained simply by reparametrization.
Next we observe
\begin{Proposition}\label{symmetric-propo}
  If $\cU=-\cU$ then $\cT_\cU(g_0)=\cT_\cU(g_0^{-1})$ for all $g_0\in G$.
\end{Proposition}
\begin{Proof}
  Take $g:[0,t^*]\to G$ with $g(0)=\1$, $g(t^*)=g_0$, and let $\tilde
  g(t)=g_0^{-1}\,g(t^*-t)$. Then $\tilde g(0)=\1$, $\tilde g(t^*)=g_0^{-1}$,
 and
$$
\tg^{-1}(t)\dot{\tg}(t)=g(t^*-t)^{-1}g_0g_0^{-1}\dot g(t^*-t)\,(-1)
=-g(t^*-t)^{-1}\dot g(t^*-t)\in-\cU.
$$
Since $\cU=-\cU$, the claim follows.
\end{Proof}
For the metric $\wh D$ this means that $\wh D(g_0,g_1)$ is symmetric if
$\cU=-\cU$, for then
$\wh D(g_1,g_0)=\wh D(\1,g_1^{-1}g_0)=\wh D(\1,g_0^{-1}g_1)=\wh D(g_0,g_1)$.
Finally we observe that group automorphisms interact well with the ODE
$\dot g(t)=g(t)u(t)$. For $g\in G$ we denote left-multiplication with~$g$ by
$\lambda_g\colon G\to G$, $\lambda_g(g_0)=gg_0$. As $\lambda_g$ is
differentiable, we denote its differential by $d\lambda_g$.

\begin{Proposition}
  Let $\sigma\in\Aut(G)$ and let $\sigma'=d\sigma(\1)$. Then
$$
d\sigma(g)=d\lambda_{\sigma(g)}\,\sigma'\,d\lambda_{g^{-1}}(g)
\qmb{for all $g\in G$.}
$$
In particular, if $g(t)$ is such that
  $g^{-1}\dot g\in\cU$, then $\tg(t):=\sigma(\tg(t))$ satisfies
   $\tg^{-1}\dot{\tg}\in\sigma'(\cU)$.
\end{Proposition}

As an immediate consequence we obtain the following estimate:

\begin{Proposition}\label{aut-inequality-propo}
  Let $\sigma\in\Aut(G)$ and $\cU,\tcU\subseteq\g$
  such that $\sigma'\,\cU\subseteq\tcU$.
  Then
$$
\cT_{\tcU}(g_0) \leq \cT_{\cU}(\sigma^{-1}(g_0))\qmb{for all $g_0\in G$.}
$$
\end{Proposition}
\begin{Proof}
  Let $g\colon[0,t^*]\to G$ be absolutely continuous with
  $g^{-1}\dot g\in\mathcal{U}$ a.e.,
 $g(0)=\1$, $g(t^*)=\sigma^{-1}(g_0)$.
 Then
  $\tilde g(t):=\sigma(g(t))$ satisfies
$$
\tilde g^{-1}\,\dot{\tilde g} \in\sigma'(\cU)\subseteq\tcU,\quad \tg(0)=\1,\ \tg(t^*)=g_0.
$$
Hence $\cT_{\tcU}(g_0)\leq t^*$ follows.
Since we may choose $g$ such that $t^*$ is arbitrarily close to 
$\cT_{\cU}(\sigma^{-1}(g_0))$, 
our claim follows.
\end{Proof}

The material symmetry axiom~(Sy4) can be re-interpreted now in the following way:
right multiplication with $g\in G$ leaves $\wh D$ invariant, iff
the inner automorphism $\Conj{g}=(g_0\mapsto gg_0g^{-1})\colon G\to G$
leaves~$\cU$ invariant, i.e., $\Ad(g)\cU=\cU$. This implies that
$\Conj{g}$ leaves $\cT_\cU$ invariant: $\cT_\cU=\cT_\cU\circ\Conj{g}$.
In other words: $\Conj{g}$ is an isometry for the distance $\wh D$.

\subsection*{Single-crystal plasticity}

In single-crystal plasticity the plastic flow occurs through plastic slips
induced by movements of dislocations. These movements are generated by
shearings or \textbf{slip systems}, say, $S^\alpha= x_\alpha y_\alpha^\tra$,
$\alpha=1\dots m$, where
$x_\alpha,y_\alpha\in\R^d$, $\norm{x_\alpha}=\norm{y_\alpha}=1$, and $x_\alpha\perp y_\alpha$.
Geometrically, $x_\alpha$ is the \textbf{slip direction} and $y_\alpha$ is the unit normal
of the \textbf{slip plane}.
All plastic flow has the form $\dot g=g\,\sum_\alpha \nu_\alpha S^\alpha$
with $\nu_\alpha\geq0$. Formally one distinguishes between $S^\alpha$ and
$-S^\alpha$ because mechanically the slip strains in these directions
must be distinguished, cf.~\cite{gurtin-single-crystals}.

The associated Lie algebra is
$\g=\ggen{\setof{S^\alpha:\alpha=1,\dots,m}}\subseteq\sL(d)$
because of $\trace(S^\alpha)=0$, $\alpha=1,\dots,m$.
In this case the dissipation functional has the form
$$
\wt\Delta(X)=\min\setof{\sum_\alpha \kappa_\alpha\gamma_\alpha:
  \gamma_\alpha\geq0,\ X=\sum_\alpha\gamma_\alpha S^\alpha}
$$
where $\kappa_\alpha>0$ are threshhold parameters.  The set
 $\cU=\seto{X:\wt\Delta(X)\leq1}$ is a convex polytope:
$$
\seto{X:\wt\Delta(X)\leq 1}
=\conv\left(\seto 0\cup\seto{\kappa_\alpha^{-1}\,S^\alpha:\alpha=1,\dots,m} \right).
$$
Indeed, since $\wt\Delta(S^\alpha)\leq \kappa_\alpha$,
the inclusion $\cU\supseteq\conv(\{0\}\cup\setof{\kappa_\alpha^{-1}\,S^\alpha})$
 obviously holds true. Conversely, if $\wt\Delta(X)\leq1$, we find
$\gamma_1,\dots,\gamma_m\geq0$ such that $X=\sum_\alpha \gamma_\alpha
S^\alpha$ and $\sum_\alpha \kappa_\alpha\gamma_\alpha\leq1$.
 Hence $X=\sum_\alpha \lambda_\alpha(\kappa_\alpha^{-1}S^\alpha)$
with $\lambda_\alpha:=\kappa_\alpha\gamma_\alpha\geq0$,
$\sum_\alpha\lambda_\alpha\leq1$.
Whence $X\in\conv(\{0\}\cup\setof{\kappa_\alpha^{-1}\,S^\alpha})$.

\medskip
Thus the factorization problem described in the introduction
can be interpreted as the problem of finding dissipation minimizing
paths for a $2$-dimensional single-crystal with four slip systems
$S^1=P$, $S^2=Q$, $S^3=-P$, $S^4=-Q$ and equal slip rates
$\kappa_\alpha\equiv1$. Equivalently, we can say that
we have  two slip systems $S^1$, $S^2$ and a symmetric
dissipation functional:  $\wt\Delta(-X)=\wt\Delta(X)$ for all $X\in\sL(2)$.
 The factorization cost $\cT(g)$
is nothing but the dissipation distance from the identity:
 $\cT(g)=\wt D(g)$.


\section{The associated optimal control problem}

Before we state the control problem we fix some more notation.
The Lie algebra $\sL(2)$ is the set of $2\times2$-matrices of zero trace.
The bracket is the commutator: $[X,Y]=XY-YX$.
The following matrices form a basis of $\sL(2)$:
$$
H=\matH,\quad T=P+Q=\matT,\qaq U=P-Q=\matU.
$$
We observe that $P=\frac12(T+U)$ and $Q=\frac12(T-U)$.
\begin{figure}[htbp]
  \centering
 \nobox{\includegraphics{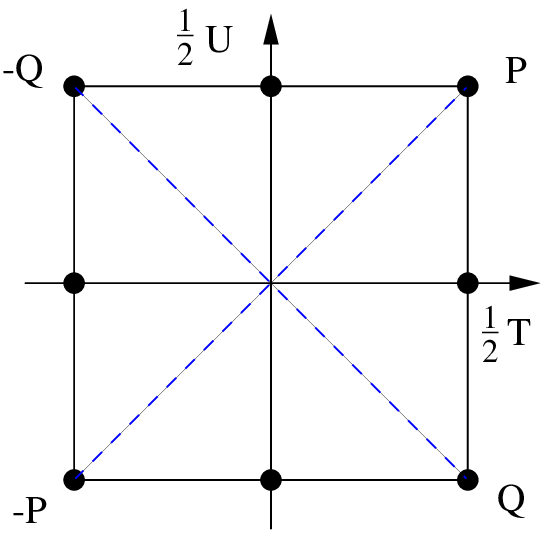}}\quad
\includegraphics[height=60mm]{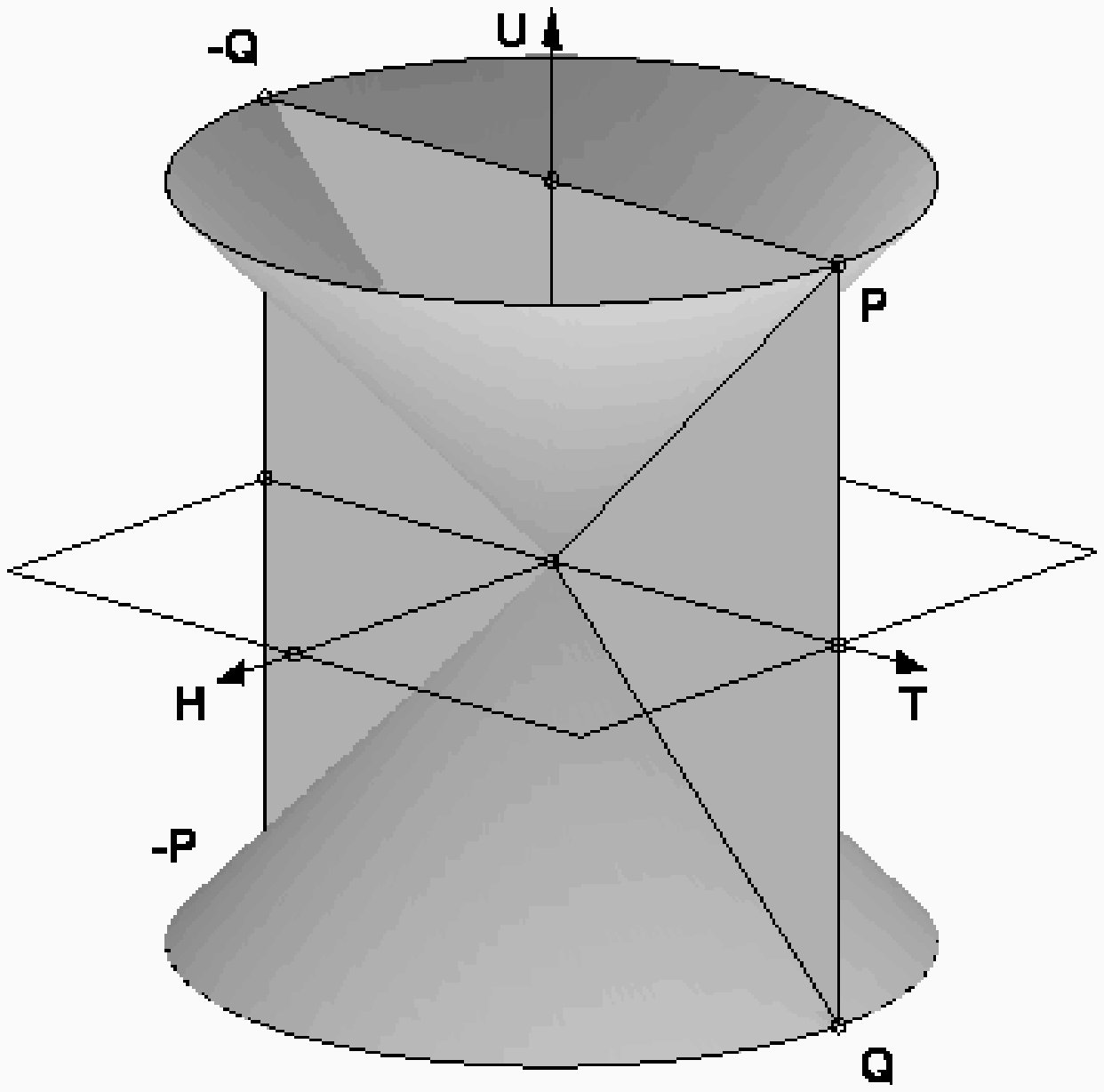}
\caption{The set $\cU=\conv(\pm P,\pm Q)\subseteq\sL(2)$ and
  the Lorentzian double cone $\setof{X\in\sL(2):\det(X)=0}$.}
  \label{square-pic}
\end{figure}

Let $\cU=\conv(\pm P,\pm Q)\subseteq\sL(2)$. This set is simply a square in the
plane $\R T+\R U\subseteq\sL(2)$. Figure~\ref{square-pic} shows how $\cU$
is situated in~$\sL(2)$. The Lorentzian double cone consists of all
matrices $X\in\sL(2)$ such that $\det(X)=0$. It is the set of all
possible two-dimensional slip systems (plus the origin),
cf. the discussion at the end of the previous section.
The elements in the
interior of the double cone all have purely imaginary spectrum, so they
are generators for compact (circle) subgroups of~$\SL(2)$.
The plane $\R H+\R T\subseteq\sL(2)$ is the set of symmetric matrices
in~$\sL(2)$.  All elements outside the double cone are 
diagonalizable (over~$\R$), in fact they are
conjugate to $\lambda H$ for some $\lambda>0$.

\smallskip
Now we  consider the
left-invariant control system (ODE on $\SL(2)$):
$$
\dot{g}(t)= g(t)\,u(t),\quad g(t)\in\SL(2),\ u(\cdot)\in
L^\infty(\R;\mathcal{U}).
\eqno{\rm(LICS)}
$$
Admissible control functions are measurable, essentially bounded functions.
The factorizations we are looking for are in one-to-one correspondence to 
the trajectories of~(LICS) generated by
piecewise constant controls. Given $t_1,\dots,t_k >0$ and
$A_1,\dots,A_k\in\mathcal{U}$, we let $\tau_j=\sum_{i=1}^jt_i$, ($j=0,\dots,k$) and define the
control $u\colon[0,\tau_k]\to \mathcal{U}$ by $u(t)=\sum_j A_j\,
\chi_{[\tau_{j-1},\tau_j)}$.
\begin{center}
  \nobox{\includegraphics{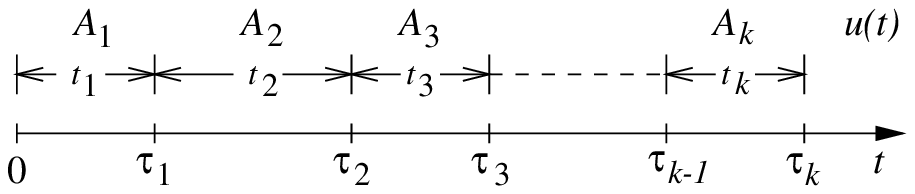}}
\end{center}
Let $g(t)$ denote the associated trajectory of~(LICS) with 
$g(0)=\1$. Then
$$
\label{trajectory-for-mex}
g(t)=
\begin{cases}
   \exp(t A_1), &\quad t\in[0,\tau_1) \\
   \mex{t_1 A_1,\, (t-\tau_1)A_2},  &\quad  t\in[\tau_1,\tau_2),\\
   \vdots & \vdots \\
   \mex{t_1 A_1,\, t_2 A_2,\dots,t_{j-1}A_{j-1},\, (t-\tau_{j-1})A_j},  &\quad t\in[\tau_{j-1},\tau_j),\\
\end{cases}
$$
and $g(\tau_k)=\mex{t_1A_1,\dots, t_kA_k}$.

\bigskip The system (LICS) is controllable because $P,Q$ generate $\sL(2)$ as
a Lie algebra: $[P,Q]=H$, $\sL(2)=\R H+\R P+\R Q$.  Therefore every $g_0\in
\SL(2)$ is reachable from the group identity. Our Optimal Control Problem
(OCP) consists of finding a time-minimal trajectory from the group
identity~$\1$ to~$g_0$:
$$
\int_0^{t^*} \rmd t \longrightarrow\mbox{min, subj. to}\ 
\dot{g}=g\,u,\ u\in \mathcal{U},\ g(0)=\1,\ g(t^*)=g_0.
\eqno{\rm(OCP)}
$$
Since $\mathcal{U}$ is compact and convex, a standard result in control
theory states that time-optimal arcs (and controls) always exist. Using
standard methods we will show that the corresponding controls are piecewise
constant.  This will tell us which factorizations have a chance to be optimal
and in which sense the original factorization problem has to be modified in
order to have solutions.  Eventually our goal is
to classify the optimal arcs of~(OCP).


\subsection*{Related problems}

There may be several non-isomorphic (connected) Lie groups having the
same---i.e., isomorphic---Lie algebras. For example, the Lie algebras $\so(3)$
(real, skew-symmetric, $3\times3$) and $\su(2)$ (complex, skew-hermitian,
$2\times2$) are isomorphic,
but the groups $\SO(3)$ and $\SU(2)$ are not ($\SO(3)$ has trivial center
while the center of $\SU(2)$ is $\{\pm\id_2\}$).

The local structure of the control system (LICS) only depends on the Lie algebra
structure. Therefore, if we consider (OCP) on any other group~$G$ with Lie
algebra isomorphic to $\sL(2)$, the whole discussion of the (PMP) will
apply equally to any such~$G$.

For example, if $B=\diag(1,1,-1)$ denotes a bilinear form of Lorentzian signature
on~$\R^3$, then we can
consider (OCP) on the group
$$
\SO_0(2,1)=\setof{g\in \mat(3,\R): g^\tra Bg=B,\ \det(g)=1,\ g_{33}>0}
$$
because the Lie algebras $\so(2,1)=\setof{X: X^\tra B+BX=0}$ and $\sL(2)$ are indeed
isomorphic. Note that $\SO_0(2,1)\not\cong\SL(2)$ (the centers consist of
one, resp. two, elements).

Therefore Problem~(OCP) is strongly related to the Hyperbolic Dubins' Problem
(HDP) and the Hyperbolic Reeds-Shepp-Car-Problem (HRSCP). We wish to emphasize
this connection because it allows nice geometric interpretations of
several of the arguments to come.

Explaining Dubins' Problem~(DP) is simple. Given two points $x_0,x_1$ and tangent
directions $v_0,v_1$ in the plane $\R^2$, the goal is to find a $C^1$-curve
$\gamma\colon[0,L]\to\R^2$
such that
\begin{itemize}
\item[(1)] $\gamma(0)=x_0$, $\dot\gamma(0)=v_0$, $\gamma(L)=x_1$, $\dot\gamma(L)=v_1$,
\item[(2)] $\gamma(s)$ has curvature $\kappa(s)$ almost everywhere, and $|\kappa(s)|\leq1$,
\item[(3)] $\gamma$ is parametrized by arc-length, and $L$ is minimal.
\end{itemize}
One may interpret this as follows: Imagine driving a car in the plane. The car
moves forward at constant speed~$1$ and its turning radius is limited.
At time $t=0$ the car is located at $x_0$, pointing into the $v_0$-direction.
The goal is to drive the car in minimal time to position~$x_1$, pointing into
the $v_1$-direction. If we allow the car to move backward as well as forward,
we obtain the so-called Reeds-Shepp-Car-Problem (RSCP).

The original sources for these problems are~\cite{dubins}
and~\cite{reeds-shepp}.
Variational problems with cost functional depending only on the
curvature can be generalized directly to manifolds with curvature, in particular
the sphere $\mathbb{S}^2$ and the hyperbolic plane~$\HH^2$ (constant
curvature~$1$, resp.~$-1$). In the latter case this leads to control systems
on the Lie groups $\SO(3)$, resp. $\SO_0(2,1)$, an observation first made by Jurdjevic in
his paper~\cite{jurd-elastica} on non-euclidean elastica.

The hyperbolic Dubins' Problem (HDP) has been investigated  by Monroy 
in~\cite{monroy,monroy-dubins}. 
We observed in~\cite{dubins-controllable} that (HDP) is equivalent to
finding time-optimal paths for the control system:\label{hdp-psl2-system}
$$
\dot g=d\lambda_g(\1) u,\quad u \in\conv(P,Q)\mbox{ a.e., }g\in\PSL(2,\R),
\eqno{\rm(HDP)}
$$
where $\PSL(2,\R)=\SL(2,\R)/\setof{\pm\1}\cong\SO_0(2,1)$. The controls $P,Q$ correspond to
left, resp., right turns, and the control $\frac12(P+Q)=\frac12 T$
corresponds to a geodesic arc.

Thus a product of the form $\mex{rP,sT,tQ,\dots}$ corresponds to a path in the
hyperbolic plane consisting of circular arcs and geodesic segments. Equality
of two such products means that two (seemingly different) paths have the same
initial and terminal positions and tangents.

In order to visualize such paths we use
the so-called \textbf{conformal disc model} of~$\HH^2$. Here $\HH^2$ is identified (as
a set) with the open unit disc
$\mathbb{D}=\setof{z\in\C: |z|<1}$. The geodesics are either diameters of
$\mathbb{D}$ or circular arcs perpendicular to the unit circle. The circular
arcs are parts of so-called \textbf{horocycles}, i.e., ordinary (euclidean) circles
touching the unit circle from inside. 


\section{Symmetries and isometries}

Apparently the set $\mathcal{U}=\conv(\pm P,\pm Q)\subseteq\sL(2)$ has
some symmetries which we would like to exploit.
 Let $\imath$ denote inversion, i.e., $\imath(g)=g^{-1}$.
Then we already observed that $\mathcal{T}(\imath(g))=\mathcal{T}(g)$ for all
$g\in\SL(2)$ because of~$\mathcal{U}=-\mathcal{U}$.
The appropriate strategy for finding more symmetries is
to look for group automorphisms preserving the factorization cost
$\mathcal{T}$.

Consider $\sigma_H,\sigma_T,\sigma_U:\SL(2)\to\SL(2)$,
$$
\sigma_H(g)=HgH,\ \sigma_T(g)=TgT,\ \sigma_U(g)=UgU^\tra=-UgU.
$$
We observe that $\sigma_H,\sigma_T,\sigma_U\in\Aut(\SL(2))$. For example,
$\sigma_H(g_0g_1)=Hg_0g_1H=Hg_0H\,Hg_1H=\sigma_H(g_0)\sigma_H(g_1)$.
Also, $\sigma_H^2=\sigma_T^2=\id_{\SL(2)}$, and
$\sigma_U=\sigma_H\sigma_T=\sigma_T\sigma_H$.
Hence
$\Gamma:=\setof{\id_{\SL(2)},\sigma_H,\sigma_T,\sigma_U}\subseteq\Aut(\SL(2))$ is a
group, actually $\Gamma\cong\Z_2\times\Z_2$.

For any Lie group automorphism $\sigma$ its derivative at the identity
 $\sigma':=d\sigma(\1)$ is a Lie algebra
automorphism. In the present situation we have 
$$
\sigma_H',\sigma_T',\sigma_U':\sL(2)\to\sL(2),\
\sigma_H'(X)=HXH,\ \sigma_T'(X)=TXT,\ \sigma_U'(X)=-UXU.
$$
Apparently there seems to be no difference between $\sigma_H$, and
$\sigma_H'$, for example. Nevertheless the distinction makes sense because
these maps have different domains and ranges.
Next we observe that 
$$
\begin{array}{lclcl}
\sigma_H'(P) = -P & \ & \sigma_T'(P)=Q, &\ & \sigma_U'(P)=-Q,\\
\sigma_H'(Q) = -Q & \ & \sigma_T'(Q)=P, &  & \sigma_U'(Q)=-P.
\end{array}
$$
Thus $\sigma'(\mathcal{U})=\mathcal{U}$ for all $\sigma\in\Gamma$.
These are the symmetries we have been looking for because the following proposition
holds true for any Lie group $G$:
\begin{Proposition}
  Assume that $\sigma\in\Aut(G)$ and let $\sigma'=d\sigma(\1)\in\Aut(\g)$.
  If $\sigma'(\mathcal{U})=\mathcal{U}$,
  then $\mathcal{T}(\sigma(g))=\mathcal{T}(g)$ for all $g\in G$.
  Moreover, $D(g_0,\,g_1)=D(\sigma(g_0),\,\sigma(g_1))$ for all $g_0,g_1\in G$, i.e.,
  $\sigma$ is an isometry of the metric $D$.
\end{Proposition}

\begin{Remark}
  Although inversion $\imath$ preserves $\mathcal{T}$, it is \textbf{not}
  an isometry of the metric $D$, in general. In fact, it is an easy exercise to show
  that for a left-invariant metric $D$ on $G$ we have:
$$
  \mbox{$\imath$ is an isometry}\iff\mbox{$D$ is right-invariant}
  \iff\mbox{$D$ is bi-invariant}.
$$
 The left-invariant metric defined by~$\mathcal{U}$ will be bi-invariant,
 for example, if $G$ is abelian or if $\Ad(G)\mathcal{U}=\mathcal{U}$.
 But one cannot expect bi-invariance  otherwise.
\end{Remark}
\medskip
A general proof is easily obtained using Proposition~\ref{aut-inequality-propo}.
In the special problem that we consider here, an elementary computation allows to
verify that every $\sigma\in\Gamma$ maps trajectories of (LICS) onto
trajectories. Assume $\dot\gamma(t)=\gamma(t) u(t)$ and
let, for example, $\eta(t)=\sigma_H(\gamma(t))=H\gamma(t)H$. Then
$$
\dot\eta(t)=H\dot\gamma(t)H=H\gamma(t)u(t)H=H\gamma(t)H\,Hu(t)H=\eta(t)\,\sigma_H'(u(t)).
$$
Hence $\mathcal{T}(\sigma_H(g))\leq\mathcal{T}(g)$ for all $g$ follows. As the
same is true for $\sigma_H^{-1}$, $\mathcal{T}(\sigma_H(g))=\mathcal{T}(g)$ follows.

\smallskip
In terms of the basis $\{H,T,U\}$ of $\sL(2)$ the maps 
$\sigma_H',\sigma_T',\sigma_U'$ are nothing but 180 degree rotations around
the $H$-, $T$-, and $U$-axis. One can actually show that 
$\Gamma=\setof{\sigma\in\Aut(\SL(2)): \sigma'(\mathcal{U})=\mathcal{U}}$.
Finally, let $\tilde\Gamma\subseteq\Diffeo(\SL(2))$
 denote the group generated by
$\Gamma\cup\{\imath\}$, then
$$
\tilde\Gamma
=\setof{\id_{\SL(2)}, \sigma_H,\sigma_T,\sigma_U,
\imath,\imath\sigma_H,\imath\sigma_T,\imath\sigma_U}
.
$$
In particular, $\tilde\Gamma$ is abelian (isomorphic to
$\Z_2\times\Z_2\times\Z_2$), and $\mathcal{T}(\phi(g))=\mathcal{T}(g)$ for all
$g\in\SL(2),\phi\in\tilde\Gamma$.


We will use the group $\tilde\Gamma$ in various ways. First it will allow us to
streamline the discussion of the (PMP) and the comparison arguments.
Practically it will allow us to reduce the number of factorization maps that
have to be inverted from $64$ down to~$13$: instead of solving $64$
 (systems of three) equations for the same right hand side, say $g$, one may
solve $13$ systems for various righthand sides from
$\setof{\phi(g):\phi\in\tilde\Gamma}$. Although this does not really affect the
overall computational cost, it is a great help for programming, testing, and debugging.

It is also worthwhile to mention that the fixed point sets
$\Fix(\phi):=\setof{g:\phi(g)=g}$,($\phi\in\tilde\Gamma$)
provide information about the cut-locus. In many cases geodesics lose their
global optimality once they hit some $\Fix(\phi)$.



\section{The structure of extremals}
\label{sec:pmp-discussion}

We already observed that for  $\mathcal{U}=\conv(\pm P,\pm Q)\subseteq\sL(2)$ 
the left-invariant control system 
$$
  \dot{\gamma}(t)={\gamma(t)}u(t),\quad u(t)\in \mathcal{U}\quad a.e.,
$$
is controllable, and that every $g$ is reachable from the
group identity in minimal time.

Next we want to apply the Pontrjagin Maxmum Principle (PMP) to obtain necessary
conditions for optimality.  For invariant systems on Lie groups the (PMP)
takes a particularly simple form, cf.~\cite{jurd-book,dimi-oberwolf}.
A proper statement requires some extra terminology: the Lie algebra dual and
the adjoint and coadjoint action.

\subsection*{Adjoint and coadjoint action}

Let $G$ be an arbitrary Lie group with Lie algebra $\g$. Let $\g^*$ denote the
vector space dual of $\g$. Then there is a natural action of $G$ on~$\g$ and, by
duality, also on $\g^*$. These are the so-called \textbf{adjoint} and
\textbf{coadjoint action}. Both actions come from conjugation on the group.
For $g_0\in G$ let $\Conj{g_0}=(g\mapsto g_0gg_0^{-1}):G\to G$. Then 
$\Conj{g_0}\in\Aut(G)$. The adjoint action is obtained by differentiating
$\Conj{g_0}\colon G\to G$ at the group identity~$\1$,
$\Ad(g_0)=d\Conj{g_0}(\1)$.
The coadjoint acton is obtained via duality.
 Instead of a coordinate-free discussion we now give explicit
representations (in coordinates) for $G=\SL(2)$ because 
these will be needed in the subsequent discussion  of the structure of extremals.

So consider $G=\SL(2)$ and $\g=\sL(2)$.  Fix the basis
$\{H,T,U\}$ of~$\g$. So we can write $X\in\sL(2)$ as $X=hH+tT+uU$ with $(h,t,u)\in\R^3$.
Next we fix a dual basis in $\g^*$, so we may write $p\in\g^*$ as a row vector
$p=(p_H,p_T,p_U)$, and $\left\langle p,hH+tT+uU\right\rangle=p_Hh+p_Tt+p_Uu$.

The adjoint action is conjugation $\Ad(g)X=gXg^{-1}$ for $g\in\SL(2)$,
$X\in\sL(2)$. In terms of the basis $\{H,T,U\}$ $\Ad(g)$ is a
$3\times3$-matrix (actually $\Ad(g)\in\SO_0(2,1)$).
Although we will not need the explicit expression, we state it just for sake of
completeness:
\begin{equation}
  \label{eq:Adg}
\left.\begin{array}{c}
g=\twomatrix{a}{b}{c}{d}\\
 \rule{0mm}{14pt}\det(g)=1
\end{array}\right\}\Rightarrow
\Ad(g)=
  \begin{pmatrix}
     ad+bc & -ac+bd                  & -ac-bd                  \\
    -ab+cd & \frac{a^2-b^2-c^2+d^2}2 & \frac{a^2+b^2-c^2-d^2}2 \\
    -ab-cd & \frac{a^2-b^2+c^2-d^2}2 & \frac{a^2+b^2+c^2+d^2}2
  \end{pmatrix}.
\end{equation}
For $X\in\sL(2)$ let $\ad(X)\colon\sL(2)\to\sL(2)$, $\ad(X)Y=[X,Y]$. If
$X=hH+tT+uU$, then in terms of the basis $\{H,T,U\}$ we obtain
$$
\ad(X)=
\begin{pmatrix}
  0 & 2u & -2t \\ -2u & 0 & 2h \\ -2t & 2h & 0
\end{pmatrix}.
$$
We note that $\ad\colon\sL(2)\to\so(2,1)$ is a Lie algebra isomorphism.
For $P=\frac12(T+U)$ and $Q=\frac12(T-U)$ we obtain:
\begin{equation}
  \label{eq:adpq}
  \ad(P)=\begin{pmatrix}
  0 & 1 & -1 \\ -1 & 0 & 0 \\ -1 & 0 & 0
\end{pmatrix},\quad
  \ad(Q)=\begin{pmatrix}
  0 & -1 & -1 \\  1 & 0 & 0 \\ -1 & 0 & 0
\end{pmatrix}.
\end{equation}
Since $\ad(P)^3=\ad(Q)^3=0$, we obtain
\begin{equation}
  \label{eq:e2adpq}
  e^{\tau\ad(P)}=
\begin{pmatrix}
     1 &             \tau & -\tau \\
 -\tau & 1-\frac{\tau^2}2 &   \frac{\tau^2}2 \\
 -\tau &  -\frac{\tau^2}2 & 1+\frac{\tau^2}2
\end{pmatrix},
\ %
  e^{\tau\ad(Q)}=
\begin{pmatrix}
     1 &            -\tau & -\tau \\
  \tau & 1-\frac{\tau^2}2 &  -\frac{\tau^2}2 \\
 -\tau &   \frac{\tau^2}2 & 1+\frac{\tau^2}2
\end{pmatrix}.
\end{equation}
Using the duality between $\sL(2)$ and $\sL(2)^*$ one obtains
$\Ad(g)^*p=p\circ\Ad(g)$, $\ad(X)^*p=p\circ\ad(X)$,
and with our choice of coordinates the latter is nothing but
left-multiplication
of the row vector~$p$ with the matrix $\Ad(g)$, resp. $\ad(X)$.
The \textbf{coadjoint action} is defined as 
$$
\CoAd(g)p=\Ad(g^{-1})^*p,\quad g\in G,\ p\in\g^*. 
$$
We are more in favor of this notation (rather than the frequently used
$\Ad^*(g)$ because it prevents confusion between $\Ad(g)^*$ and $\Ad^*(g)=\Ad(g^{-1})^*$.

\smallskip
The \textbf{Killing form} is defined as $\kappa(X,Y)=\trace(\ad(X)\ad(Y))$.
For $\sL(2)$ it is a nondegenerate symmetric bilinear form of
signature~$(+,+,-)$. The adjoint action leaves $\kappa$ invariant,
in particular the quadratic form $q_\kappa(X):=\kappa(X,X)=\trace(\ad(X)^2)$
is $\Ad$-invariant, $q_\kappa(\Ad(g)X)=q_\kappa(X)$ for all $X\in\sL(2)$, $g\in\SL(2)$.
With our choice of coordinates $q_\kappa(hH+tT+uU)=8(h^2+t^2-u^2)$.

By duality, the quadratic form $\mathcal{C}(p)=p_H^2+p_T^2-p_U^2$ is
$\Ad^*$-invariant, i.e., $\mathcal{C}(\Ad(g)^*p)=\mathcal{C}(p)$
for all $p\in\sL(2)^*$, $g\in\SL(2)$.
As we will see soon, $\mathcal{C}(p)$ will appear as a \textbf{first integral} for
\textbf{any} optimal control problem on \textbf{any} Lie group with
Lie algebra isomorphic to~$\sL(2)$.

For $p\in\sL(2)^*$ let $\mathcal{O}_p=\setof{\Ad(g)^*p: g\in\SL(2)}$ denote
its orbit under the coadjoint action. Roughly speaking, the coadjoint orbits
are level sets of $\mathcal{C}$, more precisely we have four different types
of orbits:
\begin{description}
\item[Hyp1:] If $\mathcal{C}(p)>0$ the orbit $\mathcal{O}_p$ is a one-sheeted
  hyperboloid.
\item[Hyp2:] If $\mathcal{C}(p)<0$ the orbit $\mathcal{O}_p$ is the upper
 ($p_U>0$) or  lower ($p_U<0$) part of a two-sheeted hyperboloid.

\item[Cone:] If $\mathcal{C}(p)=0$, then either $p=0$ and $\mathcal{O}_p=\{0\}$ is
  singleton, or $p\neq0$ and $\mathcal{O}_p$ is the upper ($p_U>0$) or  lower
  ($p_U<0$) part of the (boundary of the) Lorentzian double cone, cf. Figure~\ref{square-pic}.
\end{description}

For $p\in\sL(2)^*$ we denote $G_p=\setof{g\in\SL(2):\Ad(g)^*p=p}$ the
\textbf{stabilizer group of $p$}, and $\g_p=\setof{X\in\sL(2):\ad(X)^*p=0}$
the \textbf{stabilizer algebra of~$p$}. One quickly verifies that for
$p=(p_H,p_T,p_U)\neq0$ one has $\g_p=\R (p_H H+p_T T-p_U U)$.

\subsection*{The Pontrjagin Maximum Principle on Lie groups\\ (for time-optimal problems)}

\begin{Theorem}
  Let $G$ be a Lie group with Lie algebra $\g$, $\mathcal{U}\subseteq\g$, and
  consider the \textsc{(LICS)} $\dot\gamma(t)=\gamma(t)u(t)$, $u(t)\in\mathcal{U}$.
  Assume that $g(t)$ is a trajectory which is time-optimal on some
  interval~$I$. Let $u(t)=g(t)^{-1}\dot g(t)\in L^\infty(I;\mathcal{U})$.  Then
  there exists an absolutely continuous covector function $p:I\to\g^*$
  such that
\begin{itemize}
 \item[\textsc{(0)}] $p(t)\neq0$ for some {\rm(}hence all\/{\rm)} $t\in I$,
 \item[\textsc{(1)}] $\Ad(g(t)^{-1})^*p(t)$ is constant,
 \item[\textsc{(2)}] $\langle p(t),u(t)\rangle=\min_{v\in \mathcal{U}} \langle
 p(t),v\rangle$, a.e. in~$I$,
 \item[\textsc{(3)}] $\langle p(t),u(t)\rangle$ is constant {\rm(}a.e.{\rm)},  either~$-1$, or~$0$.
\end{itemize}
\end{Theorem}

A triple $(g(t),p(t),u(t))$ consisting of a trajectory~$g(t)$, a control
$u(t)$ generating $g(t)$ and a covector~$p(t)$ satisfying conditions~(0)--(3)
is called an \textbf{extremal} of the optimal control problem.
In the present situation $u=g^{-1}\dot g$, so specifying $u$ is actually redundant.

Let ${\mathcal H}(p)=\min_{v\in \mathcal{U}} \langle p(t),v\rangle$ denote the
\textbf{optimal Hamiltonian.} By~(3) $\mathcal{H}$ is an integral of motion.
An extremal for which $\mathcal{H}(p)\equiv 0$ is called an \textbf{abnormal}
extremal. The other extremals are called \textbf{normal} or \textbf{regular}
extremals.  Geometrically an extremal being abnormal means it satisfies the
first order necessary conditions for \textbf{any} cost functional, not only
the one under consideration.

If $\mathcal{U}=-\mathcal{U}$ optimal abnormal
extremals will appear only in  exceptional, degenerate cases.
They will \textbf{never} appear, for example, if $\mathcal{U}$ is a
zero-neighborhood, for in that case $\mathcal{H}(p)=0\implies p=0$, and the
latter is impossible by~(0).

The level set $\setof{\mathcal{H}=-1}$ is the
\textbf{yield surface} found in the plasiticity literature,
and the \textbf{flow rule} is encoded in~(1) and~(2).
Differentiating~(1) yields $\dot p(t)=p(t)\circ\ad(u(t))$, while the
minimizing condition~(2) implies that $u(t)$ lies in the subgradient
of~$\mathcal{H}$.
As usual the derivative-free formulation has the advantage that we need not
make smoothness assumptions.

In view of our discussion of symmetries and isometries we observe:
\begin{Proposition}\label{symmetries-for-pmp}
  If $\sigma\in\Aut(G)$ satisfies $\sigma'(\mathcal{U})=\mathcal{U}$, 
  then $\sigma$ maps extremals onto extremals.
  Similarly, if $\mathcal{U}=-\mathcal{U}$ and $(g,p,u)$ is an extremal
 {\rm(} in $[0,t^*]$\rm{)}, then
  $(g(t^*-t), -u(t^*-t),\, -p(t^*-t) )$ is an extremal, too.
\end{Proposition}
The proof is an easy exercise. If $(g,p,u)$ is an extremal and $\sigma$ as
above, the image extremal is $(\sigma(g(t)),\, ((\sigma')^{-1})^*p(t),\, \sigma'(u(t)))$.

\subsection*{Yield surface and flow rule for the model problem}

Using the basis $\setof{H,T,U}$ of $\g=\sL(2)$ and fixing 
a dual basis (of~$\setof{H,T,U}$) in~$\g^*$, we write $p=(p_H,p_T,p_U)$.
With $\mathcal{U}=\conv(\pm P,\pm Q)$ the optimal Hamiltonian is
$$
{\mathcal H}(p)=-\frac12\, (|p_T|+|p_U|).
$$
We already observed that there is another integral of motion:
$$
  \mathcal{C}(p)=p_H^2+p_T^2-p_U^2
$$
because $\mathcal{C}$ is constant along coadjoint orbits.

\begin{Proposition}\label{no-abnormals-propo}
  Abnormal extremals are not optimal.
\end{Proposition}
\begin{Proof}
  Let $(g(t),p(t),u(t))$ be an abnormal extremal. Then
  $\mathcal{H}(p)\equiv0$ implies $p(t)=(p_H(t),0,0)$.
  Hence $\mathcal{C}(p)=p_H(t)^2$, whence $p_H$ is constant.
  As $p\not\equiv0$, $p_H\neq0$ must hold. As $p$ is constant,
  $p$ is differentiable. Let $u(t)=u_1(t)P+u_2(t)Q$.
  As $0\equiv \dot p=p\circ\ad(u)$, $u(t)\in\mathcal{U}\cap\g_p$ must hold.
  Since $\g_p=\R H$, $u(t)\equiv0$ follows.
  Or, elementary,
$$
0=\dot p=p\,\ad(u_1P+u_2 Q)=(p_H,0,0)
\begin{pmatrix}
        0 & u_1-u_2 & -u_1-u_2 \\
 -u_1+u_2 & 0 & 0 \\
 -u_1-u_2 & 0 & 0
\end{pmatrix},
$$
so $u_1-u_2=0$ and $u_1+u_2=0$, hence $u_1=u_2=0$.
Thus $\dot g\equiv 0$,
i.e., $g(t)\equiv g_0$ is a constant path, hence not optimal. 
\end{Proof}

\begin{figure}
\centering
\includegraphics[width=100mm]{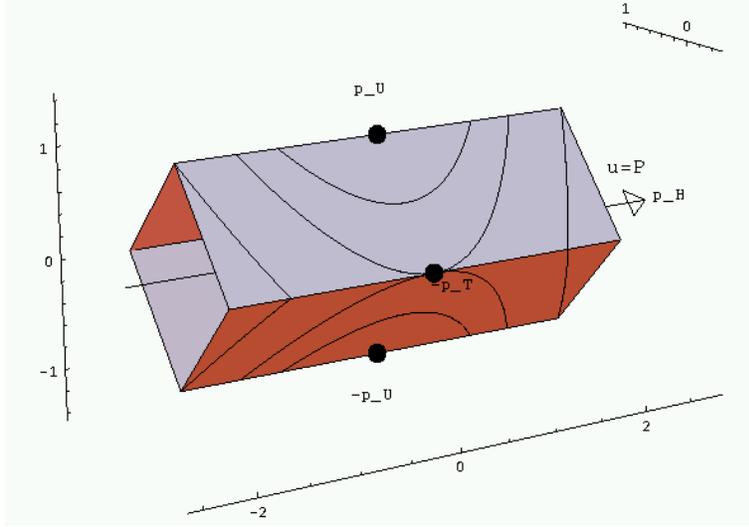}
\caption{The  yield surface $\setof{\mathcal{H}=-1}$.}
\label{pmp-flowpic}
\end{figure}

For regular extremals the covector $p(t)$ evolves on  the level set (yield surface)
$\{\mathcal{H}=-1\}$. The latter is a cylinder over a square,
cf. Figure~\ref{pmp-flowpic} which also shows some flow lines.
On each open face the minimizing control is uniquely determined. Switches
(may) occur when $p(t)$ hits one of the four edges, i.e., when $p_U(t)=0$ or
$p_T(t)=0$. At such a point the minimizing condition $\left\langle
  p(t),u(t)\right\rangle=\min_{v\in\mathcal{U}}\left\langle
  p(t),v\right\rangle$ does not suffice to characterize the control~$u(t)$
uniquely.  But we also have the geometric information
$p(t)\in\setof{\mathcal{H}=-1}$, and $\dot p(t)=p(t)\circ\ad(u(t))$.
Carefully exploiting all this information, one obtains that optimal controls
are piecewise constant and that there are $4$ types of extremals:
\begin{description}
\item[(ALT)] The optimal control $u(t)$ follows an {\em alternating switching pattern}.
$$
P\vdash (-Q)\vdash P\vdash(-Q)\dots,\qmbq{resp.} 
Q\vdash(-P)\vdash Q\vdash(-P)\dots,
$$
and the time $s$ between successive switches is a constant, $0<s<2\sqrt{2}$.
The corresponding path $g(t)$ is a subarc of 
$$
\mex{sP,-sQ,sP,-sQ,\dots}\qmbq{resp.,} \mex{sQ,-sP,sQ,-sP,\dots}.
$$

\item[(CSP)] The {\em circular  switching pattern} (CSP), here the control $u(t)$
switches in either clockwise or counterclockwise order from vertex to vertex:
$$
P \vdash Q \vdash -P\vdash -Q\vdash P\dots,\quad\mbox{resp.}\quad
P \vdash -Q \vdash -P\vdash Q\vdash P \dots,
$$
the time~$s$ between successive switches is constant, $s\in(0,\sqrt{2})$.
The corresponding path $g(t)$ is a subarc of 
$$
\mex{sP, sQ,-sP,-sQ,\dots}\qmbq{resp.,}
\mex{sP,-sQ,-sP, sQ,\dots}.
$$

\item[(SSP)] The {\em singular switching pattern(s)} (SSP).
 In this case singular arcs may occur, the corresponding  control is
 {\em not} bang-bang (i.e., in $\{\pm P,\pm Q\}$).
 The singular controls are constant, $\pm\frac12 T$, they may
 be applied on an interval of arbitrary length.
 The switching time for an intermediate bang arc is
 exactly~$\sqrt{2}$.
 Describing all possible switching sequences in general is complicated.
 All possibilities can be obtained as paths in the directed graph shown in Figure~\ref{ssp-graph}.
 \begin{figure}[hbtph]
   \centering
   \includegraphics{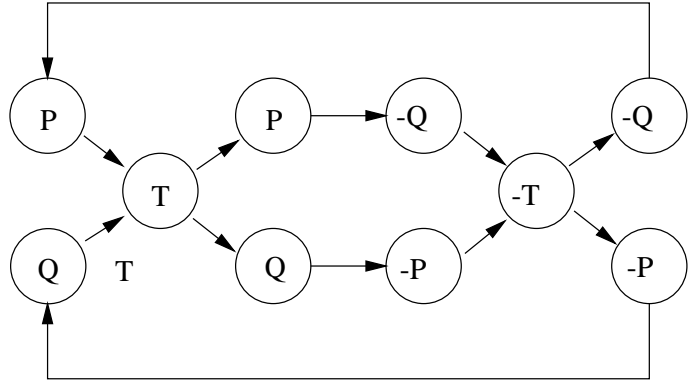}
   \caption{The directed graph describing the singular switching patterns.}
   \label{ssp-graph}
 \end{figure}
 If the $S$-arc between 2 $B$-arcs has zero length, the value of the
 control function $u(t)$ need not change. We call 
these~\textbf{virtual switches.} Example:
$$
\mex{\sqrt2\, P, s_1T, \sqrt2\, P, -\sqrt2\,Q,-\sqrt2 Q,\sqrt2
  P,s_2\,T,\sqrt2\,P,\dots},
\ s_1,s_2\geq0.
$$

\item[(U/2)] The constant controls $u_\pm \equiv\pm\frac12(P-Q)=\pm\frac12 U$.
  The corresponding path is $g(t)=\exp(\pm t\cdot\frac12 U)$, $t\geq0$.
  These are singular controls, too.

\end{description}

From this characterization we can already deduce that every $g\in\SL(2)$ has
an optimal factorization with factors from the set $\exp(\R
P)\cup\exp(\R Q)\cup\exp(\R T)\cup\exp(\R U)$.


\subsection*{How to find the switching patterns and times}
\label{switching-pattern-howto}

In this subsection we prove that indeed, our list consists of extremals only,
and we explain how this list is obtained.
We also indicate why there are no other extremals, but a
rigorous proof has to be carried out in a different way.
Therefore it is provided separately in the subseqent
subsection.

The yield surface has four faces, on each of the open faces we have
$p_Tp_U\neq0$, and the minimizing condition~(2) determines the control $u$ uniquely
(the subgradient $\partial\mathcal{H}$ is singleton).

Due to the symmetries from $\Gamma=\setof{\id,\sigma_H,\sigma_T,\sigma_U}$
it suffices to consider only one of these faces, say,
$$
\bbF=\setof{p:p_T+p_U=2,\ p_T,p_U>0},\qaq\bbbF
$$
because any other face is mapped onto $\bbF$ by
$(\sigma')^{-*}$ for some $\sigma\in\Gamma$.
The maps $(\sigma')^{-*}$ ($\sigma\in\Gamma$) are
180-degree rotations around the coordinate axes in $\sL(2)^*$. 

So let us consider $p\in\bbF$.  Then $\partial\mathcal{H}(p)=\{-P\}$.
It is crucial to analyze what happens when $p(t)$ hits the
two boundary lines $(0,2,0)+\R(1,0,0)$, resp. $(0,0,2)+\R(1,0,0)$.
Therefore our next step is to consider the flow $(p,\tau)\longmapsto p\,e^{-\tau\ad(P)}$ for
$p\in\bbbF$, $\tau\in\R$. We compute
$$
e^{-\tau\ad(P)}=
\begin{pmatrix}
    1 & -\tau            & \tau             \\
 \tau & 1-\frac{\tau^2}2 &   \frac{\tau^2}2 \\
 \tau &  -\frac{\tau^2}2 & 1+\frac{\tau^2}2
\end{pmatrix},\quad \tau\in\R.
$$
So $(p_H,0,2)e^{-\tau\ad(P)}=(p_H+2\tau,\, -\tau(p_H+\tau),\,
\tau^2+p_H\,\tau+2)$.
\begin{figure}[htbp]
  \centering
  \includegraphics{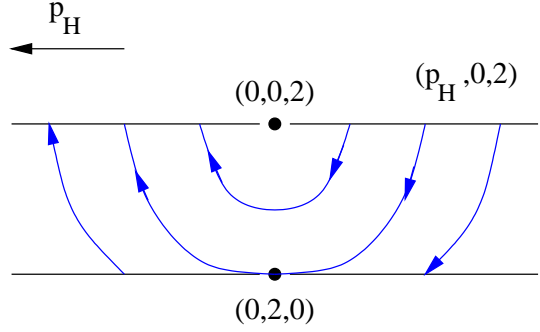}  
  \caption{The flow on the face $\bbF=\setof{p_U+p_T=2, p_U>0, p_T>0}$.}
  \label{face-flowpic}
\end{figure}
Hence we have the following cases:
\begin{description}
\item[\boldmath{$p_H<-2\sqrt2$}:] $p\,e^{-\tau\ad(P)}\in\bbF$ for
  $\tau\in(0,\tau^*)$ where $\tau^*=\frac{|p_H|-\sqrt{p_H^2-8}}2$, and
  $p\,e^{-\tau^*\ad(P)}=(-\sqrt{p_H^2-8},2,0)$.

\item[\boldmath{$p_H=-2\sqrt2$}:] $p\,e^{-\tau\ad(P)}\in\bbF$ for
  $\tau\in(0,2\sqrt2)\setminus\{\sqrt2\}$, and
  $p,e^{-\sqrt2\,\ad(P)}=(0,2,0)$, $p,e^{-2\sqrt2\,\ad(P)}=(2\sqrt2,0,2)$.

\item[\boldmath{$p_H\in(-2\sqrt2,0)$}:]  $p\,e^{-\tau\ad(P)}\in\bbF$ for
  $\tau\in(0,\tau^*)$ where $\tau^*=|p_H|>0$, and $p\,e^{-\tau^*\ad(P)}=(|p_H|,0,2)$.

\item[\boldmath{$p_H\in(0,2\sqrt2)$}:] $p\,e^{-\tau\ad(P)}\in\bbF$ for
  $\tau\in(-\tau^*,0)$ with $\tau^*=p_H$, cf. the previous case.

\item[\boldmath{$p_H=2\sqrt2$}:] $p\,e^{-\tau\ad(P)}\in\bbF$ for
  $\tau\in(-2\sqrt2,0)\setminus\{-\sqrt2\}$, cf. the second case.

\item[\boldmath{$p_H>2\sqrt2$}:] $p\,e^{-\tau\ad(P)}\in\bbF$ for
  $\tau\in(-\tau^*,0)$ where $\tau^*=\frac{p_H-\sqrt{p_H^2-8}}2$,
  $p\,e^{\tau^*\ad(P)}=(\sqrt{p_H^2-8},2,0)$.
\end{description}
Applying symmetries we obtain the flow lines on the other open faces, too.
Figure~\ref{flowgraph} shows 
the yield surface  ``unfolded'' and the various possibilities for $p(t)$.
\begin{figure}
  \centering
  \includegraphics{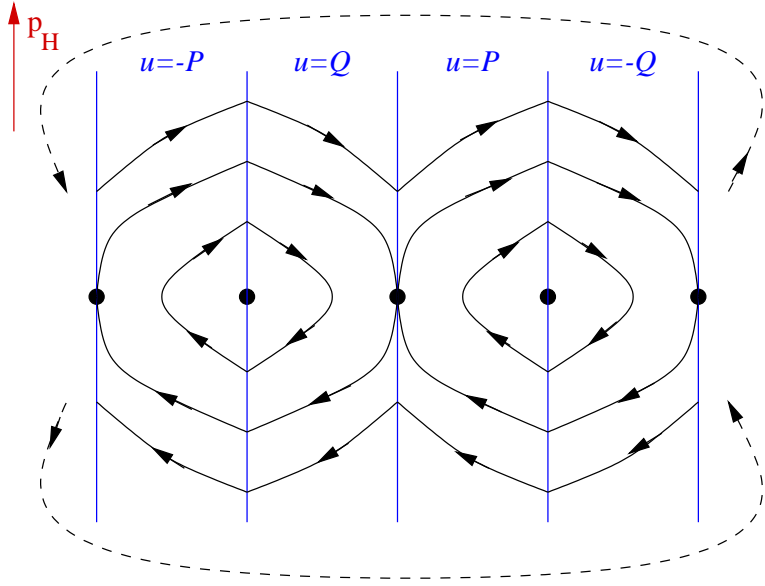}
  \caption{The Hamiltonian flow on $\setof{\mathcal{H}=-1}$.}
  \label{flowgraph}
\end{figure}
Recalling that $\mathcal{C}(p)=p_H^2+p_T^2-p_U^2$ is constant along extremals
we can state  the following
\begin{Proposition}
  Let $(g,p,u)$ be an extremal such that $p(t_0)\in\bbF$ for some~$t_0$.
  Let $C=\mathcal{C}(p)$ and
  \begin{eqnarray*}
    \alpha &=& \inf\setof{t<t_0:(\forall \tau\in(t,t_0))\ p(\tau)\in\bbF}\\    
    \beta  &=& \sup\setof{t>t_0:(\forall \tau\in(t_0,t))\ p(\tau)\in\bbF}
  \end{eqnarray*}
  Then one of the following three cases occurs:
  \begin{description}
  \item[CSP:]  $C>4$, $p_H(t_0)\neq0$,
    $\beta-\alpha=\frac{\sqrt{C+4}-\sqrt{C-4}}2\in(0,\sqrt2)$, and
    $$
    \begin{array}{lcll}
      p_H(t_0)>0 & \Rightarrow & p(\alpha)=(\sqrt{C-4},2,0), & p(\beta)=(\sqrt{C+4},0,2),\\
      p_H(t_0)<0 & \Rightarrow & p(\alpha)=(-\sqrt{C+4},0,2), & p(\beta)=(-\sqrt{C-4},2,0).
    \end{array}
    $$

  \item[SSP:] $C=4$, $p_H(t_0)\neq0$, $\beta-\alpha=\sqrt2$, and
    $$
    \begin{array}{lcll}
      p_H(t_0)>0 & \Rightarrow & p(\alpha)=(0,2,0), & p(\beta)=(2\sqrt2,0,2), \\
      p_H(t_0)<0 & \Rightarrow & p(\alpha)=(-2\sqrt2,0,2), &p(\beta)=(0,2,0).
    \end{array}
    $$

  \item[ALT:] $C\in(-4,4)$, $\beta-\alpha=\sqrt{C+4}\in(0,2\sqrt2)$, and
    $p(\alpha)=(-\sqrt{C+4},0,2)$, $p(\beta)=(\sqrt{C+4},0,2)$.
  \end{description}
\end{Proposition}
It should be obvious why we used the labels (ALT),(SSP),(CSP).
It is also obvious that the following triples $(g,p,u)$ are extremals ($t\in\R$):
$$
\left( \exp\left(\pm t\,\frac12 T\right),\ (0,\mp2,0),\ \pm\frac12 T \right),\quad
\left( \exp\left(\pm t\,\frac12 U\right),\ (0,0,\mp2),\ \pm\frac12 U \right).
$$
We refer to these as the \textbf{singular arcs} because the control $u(t)$ is
not an extreme point of $\mathcal{U}$. The first one has $\mathcal{C}(p)=4$
while the second one has $\mathcal{C}(p)=-4$.
The geometric reason for the existence of these singular arcs is that
for $p\in\{(0,0,\pm2),(0,\pm2,0)\}$ the yield surface and the coadjoint
orbit $\mathcal{O}_p$ have first order contact: 
$\partial\mathcal{H}(p)\cap \R\partial\mathcal{C}(p)\nonvoid$.

We observe that an extremal with $\mathcal{C}(p)=4$ inevitably hits one of the
points~$(0,\pm2,0)$. There we may ``glue'' it together with a singular arc.
Hence we obtain that our so-called \textbf{list of extremals} really 
consists of extremals.


In order to prove that our list of extremals is complete we must prove that
switches must occur whenever $p(t)$ hits one of the edges at some
$p\neq(0,\pm2,0)$. We first give a geometric explanation that would actually
prove our claim if we knew that $p(t)$ is piecewise differentiable.
Unfortunately the (PMP) only provides a $p(t)$ which is a.e. differentiable.
Therefore a rigorous proof has to be provided. This will be done in the next
subsection.

\begin{figure}
  \centering
  \includegraphics{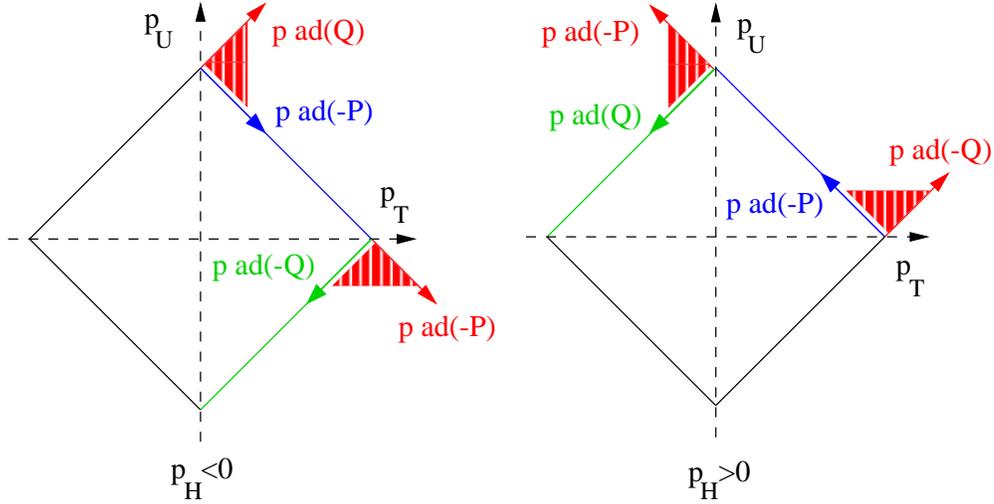}
  \caption{The reason why switches must occur on the edges $p_U=2$ ($p_T=2$), if
    $p_H\neq0$}
  \label{switch-pic1}
\end{figure}

Let us consider $p=(p_H,0,2)$ first. The minimizing condition~(2) implies
$u\in\conv(-P,Q)$. Next we compute
\begin{eqnarray*}
p\,\ad(-P)&=& (p_H,0,2)\begin{pmatrix}
  0 & -1 & 1 \\  1 & 0 & 0 \\ 1 & 0 & 0
\end{pmatrix}
=(2,-p_H,p_H)
\\
p\,\ad(Q) &=& (p_H,0,2)\begin{pmatrix} 0 & -1 & -1 \\  1 & 0 & 0 \\ -1 &
  0 & 0
\end{pmatrix}
=(-2, -p_H, -p_H).
\end{eqnarray*}
Using the fact that $p(t)$ evolves on the yield surface, $p\,\ad(u)$ must be
subtangent to the yield surface. A glance at Figure~\ref{switch-pic1}
convinces us that for $p_H\neq0$ this determines $u$ uniquely.
The same thing happens on the other edge for $p=(p_H,2,0)$ with $p_H\neq0$.
%


%
Proposition~\ref{no-abnormals-propo} gives a flavor of how to attack
a rigorous proof.
The remainder of this section is devoted to providing one,
it is not mandatory for understanding the rest of this paper.

\subsection*{Rigorous proofs}

In the previous discussion we already observed that
the cases  to distinguish are $\mathcal{C}(p)=-4$,
$\mathcal{C}(p)\in(-4,4)$, $\mathcal{C}(p)=4$, and $\mathcal{C}(p)>4$.
They correspond (in this order) to (U/2)-, (ALT)-, (SSP)-, and (CSP)-extremals, respectively.
W.l.o.g. we only consider extremals $(g,p,u)$ starting at $g(0)=\1$.

\begin{Proposition}
  Let $(g,p,u)$ be an extremal such that $p_U=0$ on an open interval~$I$. Then
  $\mathcal{C}(p)=4$, $p(t)$ is constant in $I$, $p\equiv(0,\pm2,0)$, and
  $u(t)\equiv\mp\frac12 T$ in~$I$.  
\end{Proposition}
\begin{Proof}
  Suppose $p_U\equiv0\in I$, then $\mathcal{H}(p)\equiv-1$ implies
  $|p_T|\equiv2$ in~$I$. Hence $\mathcal{C}(p)=p_H^2+p_T^2-p_U^2=p_H^2+4$, whence
  $p(t)$ is constant in $I$. 
  Thus $0=\dot p=p\,\ad(u)$, so $u\in \g_p\cap -\sign(p_T)\conv(P,Q)$.
  Since $\g_p=\R(p_H H+p_T T)$ this intersection is nonempty iff $p_H=0$.
  Thus $p\equiv(0,\pm2,0)$ and  $u\equiv \mp\frac12\,T$ in~$I$. 
  So $\mathcal{C}(p)=0+4-0$ in~$I$ and hence for all $t$ as it is a first integral.
\end{Proof}

\begin{Proposition}
  Let $(g,p,u)$ be an extremal such that $p_T=0$ on an open interval~$I$.
  Then $\mathcal{C}(p)=-4$, $p(t)$ is constant for all~$t$, either
  $(0,0,\pm2)$, $u(t)$ is constant (either $\mp\frac12 U$), and $g(t)=\exp(\mp\frac12 U)$.
\end{Proposition}
\begin{Proof}
  Suppose $p_T\equiv0$ in some interval~$I$. Then $|p_U|\equiv2$ in~$I$, and
  as in the previous proof we deduce that $p_H$ is constant in~$I$, hence
  $p(t)\equiv(p_H,0,p_U)$ is
  constant in~$I$. Thus $0=\dot p=p\,\ad(u)$ yields $u(t)\in\g_p\cap-\sign(p_U)\conv(P,-Q)$.
  As $\g_p=\R(p_H H-p_U U)$, the intersection is nonempty iff $p_H=0$.
  Hence $p(t)\equiv(0,0,\pm2)$ and $u(t)\equiv \mp\frac12 U$ in~$I$.
  Also $\mathcal{C}(p)=0^2+0^2-4=-4$.

  Since $|p_U|\leq2$ and $p_H^2+p_T^2-p_U^2=-4$ hold for \textbf{all} $t$ (not
  just in $I$),
  we deduce $p_H^2+p_T^2=0$ and $|p_U|=2$ for all $t$.
  So $p\equiv(0,0,\pm2)$ and $u=\mp\frac12 U$ follows for all~$t$.
  Hence $g(t)=\exp(\mp\frac12 U)$ follows, too.
\end{Proof}
The \textbf{switching surface} for our problem is
$\Sigma:=\setof{p:p_U=0}\cup\setof{p:p_T=0}$.
Up to symmetry ($\tilde\Gamma$) we have to consider only two cases:
$p_H<0,p_T=0,p_U=2$ and $p_H<0,p_T=2,p_U=0$.

\begin{Proposition}\label{switch-propo-u2}
  Let $(g,p,u)$ be an extremal 
  on some interval~$I$ and $t_0\in I$ such that
  $p(t_0)=(p_H^0,0,2)$ with $p_H^0<0$.
  Then there exists an $\varepsilon>0$ such that
  $p_T(t)<0$ in $(t_0-\varepsilon,t_0)$ and $p_T(t)>0$ in $(t_0,t_0+\varepsilon)$.
  In particular, $u(t)=Q$ a.e. in $(t_0-\varepsilon,t_0)$ and
  $u(t)=-P$ a.e. in $(t_0,t_0+\varepsilon)$.
  Hence a switch occurs.
\end{Proposition}
\begin{Proof}
  Since $p_U(t_0)=2$, take $\epsilon$ such that $p_U(t)>0$ for
  $t\in I_\epsilon:=I\cap(t_0-\epsilon,t_0+\epsilon)$.
  Then $u(t)\in\conv(-P,Q)$, so $u(t)=-\frac12\,U+\lambda(t)\cdot\frac12 T$
  for some measurable function $\lambda:I_\epsilon\to[-1,1]$.
  Since $p(t)$ is absolutely continuous, we have $\dot p=p\,\ad(u)$ a.e. in
  $I_\epsilon$. Hence
$$
  \dot p=p\ad\left(-\frac12\,U+\lambda\frac12 T\right)
  =(p_H,p_T,p_U)
  \begin{pmatrix}
    0 & -1 & -\lambda\\ 1 & 0 & 0\\ -\lambda & 0 & 0
  \end{pmatrix}
  =({*}, -p_H, \, {*}).
$$
  Thus $\dot p_T(t)=-p_H(t)$ a.e. in $I_\epsilon$. As $p_H$ is continuous, $p_T$ is
  differentiable in $I_\epsilon$, and since $p_T(t_0)=0$, $\dot p_T(t_0)=-p_H(t_0)>0$
  we find $\varepsilon$ such that $p_T(t)<0$ in $(t_0-\varepsilon,t_0)$ and
  $p_T(t)>0$ in $(t_0,t_0+\varepsilon)$.
  Hence $u(t)=Q$ a.e. in $(t_0-\varepsilon,t_0)$, $u(t)=-P$ a.e. in $(t_0,t_0+\varepsilon)$.
\end{Proof}

The other case is treated similarly.

\begin{Proposition}\label{switch-propo-t2}
  Let $(g,p,u)$ be an extremal 
  on some interval~$I$ and $t_0\in I$ such that
  $p(t_0)=(p_H^0,2,0)$ with $p_H^0<0$.
  Then there exists an $\epsilon>0$ such that
  $p_U(t)>0$ in $(t_0-\epsilon,t_0)$ and $p_T(t)<0$ in $(t_0,t_0+\epsilon)$.
  In particular, $u(t)=-P$ a.e. in $(t_0-\epsilon,t_0)$ and
  $u(t)=-Q$ a.e. in $(t_0,t_0+\epsilon)$.
  Hence a switch occurs.
\end{Proposition}
\begin{Proof}
  Since $p_T(t_0)=2$, we find  $\epsilon>0$ such that $p_T(t)>0$ for
  $t\in I_\epsilon:=I\cap(t_0-\epsilon,t_0+\epsilon)$.
  Then $u(t)\in\conv(-P,-Q)$, so $u(t)=-\frac12\,T+\lambda(t)\cdot\frac12 U$
  for some measurable function $\lambda:I_\epsilon\to[-1,1]$.
  Since $p(t)$ is absolutely continuous, we have $\dot p=p,\ad(u)$ a.e. in
  $I_\epsilon$. Hence
$$
  \dot p=p\ad\left(-\frac12\,T+\lambda\frac12 U\right)
  =(p_H,p_T,p_U)
  \begin{pmatrix}
    0 & \lambda & 1\\ -\lambda & 0 & 0\\ 1 & 0 & 0
  \end{pmatrix}
  =({*},{*}, p_H).
$$
  Thus $\dot p_U(t)=p_H$ a.e.,  so $p_U$ is
  differentiable in $I_\epsilon$ because $p_H$ is continuous.
  Since $p_U(t_0)=0$ and  $\dot p_U(t_0)=p_H(t_0)<0$,
  we find $\epsilon_1>0$ such that $p_U(t)>0$ in $(t_0-\epsilon_1,t_0)$ and
  $p_U(t)<0$ in $(t_0,t_0+\epsilon_1)$.
  Hence $u(t)=-P$ a.e. in $(t_0-\epsilon_1,t_0)$ and
  $u(t)=-Q$ a.e. in $(t_0,t_0+\epsilon_1)$.
\end{Proof}

This finishes the proof that our list of extremals is complete.
It is noteworthy to mention that the whole discussion involves only the
Lie algebra, its dual, and the adjoint, resp. coadjoint action.
 In particular our results are valid for
\textbf{any} group with Lie algebra isomorphic to $\sL(2)$.


\section{A sufficient family for $\SL(2)$}

Although the (PMP) provides very detailed information, it is only a first
order necessary condition. We need two types of additional arguments:
global arguments (no local condition could replace  them)
and higher order variations. Our goal is to provide a \textbf{sufficient family}
of paths, i.e., a finite set of maps $f_k:\R^3\to\SL(2)$
such that for every $g\in\SL(2)$ some map $f_k$ provides
an optimal path $\1\leadsto g$.

Our classification of extremals provides us with candidates. Looking at
alternating extremals we come up with
\begin{equation}
  \label{alt-maps}
  \begin{array}[b]{ccll}
    A3(r,s,t) & = & \mex{rP,-sQ,tP}, & 0\leq r,t\leq  s< 2\sqrt2,\\
    A4(r,s,t) & = & \mex{rP,-sQ,sP,-tQ},\\
    A5(r,s,t) & = & \mex{rP,-sQ,sP,-sQ,tP},
  \end{array}
\end{equation}
and $A6$, $A7$, \dots are defined similarly. Due to symmetry it suffices to
consider only (ALT)-extremals starting with an $\exp(r P)$-arc because the
group $\Gamma$ acts transitively on $\{\pm P,\pm Q\}$.

Similarly, when considering (CSP)-extremals, there are up to $8$ possibilities
($4$ for the first vertex, $2$ for the orientation: clockwise/ anticlockwise).
The group $\tilde\Gamma$ has $8$ elements, and it turns out that
for an odd number of factors it suffices to consider $1$ case, whereas
for an even number of factors there are $2$ distinguished cases.
We may always assume that the first factor is $\exp(r P)$, but we must distinguish
between clockwise and anticlockwise switching patterns if the total
number of factors is even. We obtain a sequence of maps starting with:
\begin{equation}
  \label{csp-maps}
  \begin{array}[b]{lcll}
    C3(r,s,t)  & = & \mex{rP,sQ,-tP}, & 0\leq r,t\leq s<\sqrt2,\\
    C4a(r,s,t) & = & \mex{rP,-sQ,-sP, tQ},\\
    C4c(r,s,t) & = & \mex{rP,sQ,-sP,-tQ}.
  \end{array}
\end{equation}
For (SSP)-extremals we can always
assume that (up to symmetry)  the first singular ($S$-) arc is $\exp(\frac s2\, T)$, and that the
preceding bang ($B$-) arc is $\exp(rP)$. For example, for $3$ or $4$ factors with one
singular arc we obtain the maps:
\begin{equation}
  \label{ssp-maps34}
  \begin{array}{lcll}
    S3P(r,s,t) & = & \mex{rP,\frac s2\,T, tP}, & s\geq0,\ r,t\in[0,\sqrt2], \\
    S3Q(r,s,t) & = & \mex{rP,\frac s2\,T, tQ},  \\
    S4P(r,s,t) & = & \mex{rP,\frac s2\,T,\sqrt2\,P,\,-tQ},  \\
    S4Q(r,s,t) & = & \mex{rP,\frac s2\,T,\sqrt2\,Q,\,-tP}.
  \end{array}
\end{equation}
Enumerating all possible patterns for a large number, say $n$, of factors
is an unpractical task since the number of possibilities grows exponentially in
$n$. After a singular ($S$-) arc we always have two
choices for the next bang ($B$-) arc.
But fortunately enough, it turns out that it suffices to consider
(SSP)-extremals with one $S$-arc only.
The proof is basically a verification of the identity:
$$
\mex{sT,\,\sqrt2\, P,-\sqrt2\, Q}=\mex{\sqrt2\, P,-\sqrt2\, Q,-s T}.\eqno{(*)}
$$
Let us abbreviate $w:=\sqrt2$. Then $(*)$ implies
$$
\mex{rP,s_1T,wP,-wQ,-s_2T}=\mex{rP,(s_1+s_2)T,wP,-wQ}.
$$
Occurence of this identity is not a miracle but it is kind of locally detected by the
(PMP): $e^{r\ad Q}e^{-r\ad P}T=-T$ $\iff$ $r^2=2$. So the singular
switching time is in a certain  sense geometrically distinguished. The hyperbolic
Reeds-Shepp-Car problem allows to visualize this neatly,
(cf. Fig.~\ref{singular-switching-time-pic},
p.~\pageref{singular-switching-time-pic}).
 This is completely analogous to
the euclidean case, cf.~\cite[Fig.~18, p.59]{suss-car}.

Eventually it turns out that one only needs three more maps:
\begin{equation}
  \label{ssp-maps57}
  \begin{array}{lcl}
    S5P(r,s,t) & = & \mex{-rQ,\,wP,\frac s2\,T, wP,-tQ},
 \\
    S5Q(r,s,t) & = & \mex{-rQ,\,wP,\frac s2\,T, wQ,-tP},  \\
    S7b(r,s,t) & = & \mex{-rQ,wP,\frac s2\,T,wP,-wQ,-wQ,tP}
  \end{array}
\end{equation}
In addition, it is convenient to consider also the following maps which are derived
from the previous ones:
\begin{equation}
  \label{ssp-maps-extra}
  \begin{array}[b]{lcl@{=}l}
    S5a(r,s,t) & = & \mex{ rP,\frac s2\,T, wP,-wQ,-tQ}     &   S4P(r,s,t+w)   \\
    S6 (r,s,t) & = & \mex{-rQ,\,wP,\frac s2\,T, wP,-wQ,-tQ}&   S5P(r,s,t+w) \\
    S7a(r,s,t) & = & \multicolumn{2}{l}{\mex{rP,\frac s2\,T,wP,-wQ,-wQ,wP,tP}}\\
    & = & \multicolumn{2}{l}{\sigma_U(S5P(r+w,s,t+w)).} 
  \end{array}
\end{equation}
\begin{table}[htbp]
  \centering
  \begin{tabular}{|c|c|c|c|c|l|l|l|}\hline
    \multicolumn{7}{|c|}{\bf Factorization maps with relevant domain 
      and cost\rule{0pt}{14pt}}\\ \hline\hline
     Type & Map &  \multicolumn{3}{c|}{Symmetry} & Domain  & Cost\\ \hline
     \rule{0mm}{14pt}ALT
     & $A3$ & $\Gamma$ & $4$ & $\imath\sigma_H$ & $0\leq r,t\leq s\leq
     2\sqrt2$
     & $r+s+t\in[0,6\sqrt2]$
    \\
     & $A4$ & $\Gamma$ & $4$ & $\imath\sigma_T$ & and $s\geq1$ & $r+2s+t\in[2,8\sqrt2]$
    \\
     & $A5$  & $\Gamma$ & $4$ & $\imath\sigma_H$ & and $s\in[\sqrt2,\sqrt3]$ &
     $r+3s+t\in[3\sqrt2,5\sqrt3]$
    \\ \hline
    \rule{0mm}{14pt}CSP
    & $C3$ & $\tilde\Gamma$ & $8$ & & $0\leq r,t\leq s \leq\sqrt2$ & $r+s+t\in[0,3\sqrt2]$
     \\ \cline{2-5}
    & $C4a$ & $\Gamma$ & $4$ & $\imath\sigma_U$ & and $s\leq1$ &$r+2s+t\in[0,4]$ \\
    & $C4c$ & $\Gamma$ & $4$ & $\imath\sigma_T$ & &$r+2s+t\in[0,4]$ \\ \hline
    \rule{0mm}{14pt}SSP
    & $S3P$ & $\Gamma$ & $4$ & $\imath\sigma_H$ & $s\geq0$,
     $r,t\in[0,\sqrt2]$ & $r+s+t\geq0$ \\
     & $S3Q$ & $\Gamma$ & $4$ & $\imath\sigma_U$ & &  \\ \cline{2-5}\cline{7-7}
     \rule{0mm}{14pt}
     & $S4P$ & $\tilde\Gamma$ & $8$ &  &  & $r+s+t+\sqrt2\geq\sqrt2$   \\
     & $S4Q$ & $\tilde\Gamma$ & $8$ &  &  &    \\ \cline{2-5}\cline{7-7}
     \rule{0mm}{14pt}
     & $S5P$ & $\Gamma$ & $4$ & $\imath\sigma_H$ & & $r+s+t+2\sqrt2\geq 2\sqrt2$    \\
     & $S5Q$ & $\Gamma$ & $4$ & $\imath\sigma_U$ & & \\
     & $S5a$ & $\Gamma$ & $4$ & $\imath\sigma_T$ & &
      \\  \cline{2-5}\cline{7-7}
     \rule{0mm}{14pt}
     & $S6$  & $\tilde\Gamma$ & $8$ &            & & $r+s+t+3\sqrt2\geq 3\sqrt2$ 
   \\  \cline{2-5}\cline{7-7}
%
%
     & $S7b$ & $\Gamma$ & $4$ & $\imath\sigma_T$ & &  \\
     \hline
     $\mathcal{F}$ & $16$ & & $76$ & \multicolumn{3}{l|}{\ }
 \\ \hline\hline
     \multicolumn{7}{|c|}{Use bigger domains and drop some maps}\\ \hline
     \rule{0mm}{14pt}drop 
     & $S5a$ & $\Gamma$ & $4$ &  \multicolumn{3}{l|}{%
$S4P:[0,\sqrt2]\times\R_+\times[0,2\sqrt2]\to\SL(2)$ } \\ \cline{2-7}
     \rule{0mm}{14pt}drop
     & $S6$  & $\tilde\Gamma$ & $8$ &
     \multicolumn{3}{l|}{$S5P:[0,2\sqrt2]\times\R_+\times[0,2\sqrt2]\to\SL(2)$} 
 \\ \hline
    $\mathcal{F}_1$ & $13$ & & $64$ &  \multicolumn{3}{l|}{\ }
 \\ \hline
     \rule{0mm}{14pt}drop
     & $A3$  & $\Gamma$ & $4$ &\multicolumn{3}{l|}{} \\
     drop & $C3$  & $\tilde\Gamma$ & $8$ &\multicolumn{3}{l|}{} \\ 
     add  & $B3$ & $\sigma_T$ & $2$ &
     \multicolumn{3}{l|}{$
\mex{rP,sQ,tP}$,\quad $|r|,|s|,|t|\leq2\sqrt2$}
\\ \hline
      $\mathcal{F}_2$ & $12$ & & $54$ &  \multicolumn{3}{l|}{\ }
\\ \hline
\end{tabular}
  \caption{A sufficient family of extremals for $\SL(2)$ and smaller families
that suffice for computing optimal factorizations, resp., $\mathcal{T}(g)$}
  \label{sl2-family-tab}
\end{table}
Table~\ref{sl2-family-tab} specifies a family $\mathcal{F}$ consisting
entirely of extremals. We obtain a list of $|\mathcal{F}|=76$ maps,
but since $\mathcal{F}$ is $\tilde\Gamma$-invariant, it suffices to specify
one map from each $\tilde\Gamma$-orbit. There are~$16$ such orbits.
We list a representative, a subgroup generating the orbit and, if nontrivial,
``the'' stabilizer in columns~2--4. Now we can state the main result for $\SL(2)$.

\begin{Theorem}\label{sl2-family-thm}
  Let $\mathcal{F}$ denote the family of maps specified in
  Table~{\rm\ref{sl2-family-tab}}. Then $\mathcal{F}$ is a sufficient family for
  \textsc{(OCP)} in $\SL(2)$, i.e., for every $g\in\SL(2)$ there exist
  $f\in\mathcal{F}$ and $(r,s,t)\in\dom(f)$ such that $f(r,s,t)=g$, and
  the associated path from~$\1$ to~$g$ has minimal length $\mathcal{T}(g)$.
\end{Theorem}

The proof will be given in the next section. Observing
 $S7b(r,s,t)=S5P(r,s,2\sqrt2)\exp(tP)$, we immediately obtain
\begin{Corollary}
  Every $g\in\SL(2)$ has an optimal factorization with at most~$6$ factors
  from $\setof{\pm P,\pm Q, \pm\frac12\,T}$. There always exists a geodesic
  from $\1$ to~$g$ with at most~$5$ switches.
\end{Corollary}
As is indicated in
Table~{\rm\ref{sl2-family-tab}}, there exist smaller families (fewer maps but
with bigger domains) which are sufficient, too. For example, we find a family
$\mathcal{F}_1$ consisting of~$64$ maps, resp., $13$ $\tilde\Gamma$-orbits,
and another family $\mathcal{F}_2$ consisting of~$54$ maps, resp., $12$
orbits.

\begin{Corollary}
  The families $\mathcal{F}_1$, $\mathcal{F}_2$ specified in
  Table~{\rm\ref{sl2-family-tab}} are sufficient for \textsc{(OCP)}, too.
\end{Corollary}

The proof is an immediate consequence of the preceding theorem.
The family~$\mathcal{F}$ is the appropriate one for visualizing the metric
spheres. That's the reason why we wanted to keep the domains as small as
possible.
The family $\mathcal{F}_2$ is appropriate for computing $\mathcal{T}(g)$.
In that case we want to keep the number of maps that have to be inverted
as small as possible. Finally, $\mathcal{F}_1$ is appropriate if we want to
find optimal factorizations because in that case we have to keep track of the
type of factorization (ALT or CSP). Thus it makes sense not to merge $A3,C3$
into $B3$.

Since computation of $\mathcal{T}(g)$ as well as finding optimal
factorizations requires inverting the maps from $\mathcal{F}_2$, resp.,
$\mathcal{F}_1$, it is important to realize that this is an easy task
which can be carried out efficiently at low computational costs.
Therefore we observe:
\begin{Remark}
 We have $\exp(tP)=\id+tP$, and $\exp(tQ)=\id+tQ$.
 In particular $f(r,s,t)$ is linear in~$r$ and~$t$  for every
 $f\in\mathcal{F}$. The dependence on~$s$ is as follows:
\begin{center}
 \begin{tabular}{|c|c|c|l|} \hline
  $A3,C3,B3$ & $A4,C4a,C4c$ & $A5$ & SSP \\ \hline
  linear & quadratic & cubic & linear in
  $\xi,\xi^{-1}$ where $\xi=e^s>0$.\\
  in~$s$ & in~$s$ & in~$s$ & $\leadsto$ quadratic equation in~$\xi$.
\\ \hline
 \end{tabular}
\end{center}
 Hence inverting the maps in $\mathcal{F}$ is nothing but a
 (time consuming) exercise in college algebra, except for $A5$ where we have
 $p(s)=e_2^\tra A5(r,s,t)e_1=s^3-2s$.
\\
 It suffices to invert one map from each $\tilde\Gamma$-orbit
 because $\sigma(f(r,s,t))=g$ iff $f(r,s,t)=\sigma^{-1}(g)$.
 So instead of inverting, say, all $54$ maps in $\mathcal{F}_2$, it
 suffices to invert $12$ representatives and apply them to several
 right hand sides, according to Table~\ref{sl2-family-tab}. An efficient
 implementation can also use the information on the domains to ``discard''
 solutions of $f(r,s,t)=g$ as soon as one of the parameters $r,s,t$ is
 not in the appropriate range, i.e., $\dom(f)$. A smart implementation
 also uses the fact that some extremals inevitably generate large cost.
 For example, an optimal $A5$-extremal has at least cost~$3\sqrt2$.
 Hence it is not necessary to check $A5$-type factorizations if some other
 factorization has already given cost less than $3\sqrt2$.
\end{Remark}


\section{Comparison arguments}

In this section we provide all the arguments necessary to prove Theorem~\ref{sl2-family-thm}.
The proof is based on comparison arguments.
We will show that certain factorizations are {\em not} optimal.
Hence they can never appear as a subarc of an optimal arc. 

We must treat the four classes of extremals separately.  The singular
(U/2)-extremals and the (CSP)-extremals are most pleasant in the sense that
the proofs are purely Lie algebraic. Therefore they apply to \textbf{any}
group with Lie algebra~$\sL(2)$. The arguments for (ALT)- and (SSP)-extremals
(partly) make very explicit use of the underlying group.
We start with the elimination of (U/2)-extremals.
\begin{Proposition}
  For all $\alpha\in\R$ the ``factorization'' $\exp\left(2\alpha\,\frac12
    (P-Q)\right)$ is not optimal.  In particular, optimal controls take values
   only in $\setof{\pm P,\pm Q, \pm\frac12 T}$, and optimal factorizations only
  use factors from $\exp(\R P\cup\R Q\cup \R T)$.
\end{Proposition}
\begin{Proof}
 W.l.o.g.  let us consider $\alpha>0$ small. With $r=\tan(\alpha/2)$ and
 $s=\sin(\alpha)$
 we obtain
$$
\exp\left(2\alpha\,\frac12   (P-Q)\right)
=\begin{pmatrix}
  \cos\alpha & \sin\alpha \\ -\sin\alpha & \cos\alpha
\end{pmatrix}
=\mex{r P,-s Q,r P}=A3(r,s,r).
$$
The cost of the lefthand side is $2\alpha$, on the righthand side  we have cost
$\sin(\alpha)+2\tan(\alpha/2)$.
Consider the difference $f(\alpha)=2\alpha-\sin(\alpha)-2\tan(\alpha/2)$.
Then a simple Taylor expansion yields
$$
0=f(0)=f'(0)=f''(0),\qaq f'''(0)=\frac12.
$$
Hence for small $\alpha>0$ the factorization on the righthand side is better
than the factorization on the lefthand side.
\end{Proof}

\begin{Remark}
  For the Hyperbolic Reeds-Shepp-Car Problem the previous result has an
  interesting consequence. It says that a ``rotation on the spot'' is never
  optimal. So if one wants to drive and return to the starting point but heading
  into a different direction, it is better to move forward (or backward)!
  This is in contrast to the euclidean case where it doesn't matter
  if one moves forward or turns on
  the spot: both paths have equal cost. 
  Of course this is important, because
  a rotation on the spot is not feasible for the Reeds-Shepp-Problem, it is only feasible
  for the convexified problem!

  This property reflects the fact that in a euclidean triangle the sum of the angles
  equals~$\pi$ whereas in the hyperbolic plane this sum is strictly less
  than~$\pi$.
\end{Remark}


\subsection{Circular (CSP)-extremals}

Our next goal is to prove:
\begin{Proposition}\label{csp-propo}
  Optimal \textsc{(CSP)}-extremals have at most~$4$ factors.
  This is actually true for any Lie group $G$ with Lie algebra $\sL(2)$.
\end{Proposition}

The proof of this result requires two different types of arguments.
The first one is elementary  and involves only the adjoint action.
It is a matter of elementary computations to verify that
$$
  e^{r\ad Q}(-P)= e^{1/r \ad P}(r^2 Q),\qaq
  e^{r\ad P}(-Q)= e^{1/r \ad Q}(r^2 P),\quad r\neq0.
$$

\begin{Proposition}\label{csp-s-ge-1}
  The factorization $\mex{r P,-s Q,-r P}$ is not optimal for
  $s\in(1,\sqrt{2})$ and $r>1$.
   The factorization $\mex{P,-Q,-P,t Q}$ is not optimal for $t>0$.  In
  particular, optimal \textsc{(CSP)}-extremals with switching time $s\geq1$ have at
  most~$4$ factors.
\end{Proposition}
\begin{Proof}
  Let $s\in(1,\sqrt2)$ and $r>1$, assume w.l.o.g. $r<s$.  Then we obtain
  $\mex{r Q,-s P,-r Q}=\mex{\frac1r P, r^2s\, Q,\,-\frac1r P}$ from the above
  equation.  Comparing the factorization costs on both sides we obtain
$$
(s+2r)-\left(\frac2r -r^2s\right) = \frac1r\,(r^2-1)(2-r s)
>0\quad\mbox{ if }  s\in(1,\sqrt{2}),\ r\in(1,s).
$$
Thus $\mex{r P,-s Q,-r P}$ is not optimal for $s\in(1,\sqrt{2})$ and $r>1$.
\\[2mm]
To prove the second claim we consider $r=s=1$ and obtain 
$$
\mex{P,-Q,-P,t Q} =\mex{Q,-P,-Q,tQ}=\mex{Q,-P,-(1-t)Q}.
$$
The lefthand side has cost $3s+t$ as opposed to $3s-t$ on the righthand side.
Thus our claim is proved.
\end{Proof}

While this was quite elementary, proving that optimal (CSP)-extremals with
switching time $s\in(0,1)$ have at most~$4$ factors requires some more
sophisticated arguments.  Therefore we first supply another, relatively simple
argument, that allows to reduce the number of factors to~5:
\begin{Proposition}
  For all $s>0$ the factorization $\mex{s P,-s Q,-s P,s Q}$ is not optimal.  
\end{Proposition}
\begin{Proof}
  For $s>0$ let $r=\frac{s}{1+s^2}\in(0,s)$. Then an elementary computation yields:
$$
\mex{r P,-s Q,-s P,r Q}=\mex{-r P,s Q,s P,-r Q}.
$$
Thus 
\begin{eqnarray*}
\mex{s P,-s Q,-s P,s Q} & = & \mex{(s-r)P,r P,-s Q,-s P,r Q,(s-r)Q}
\\ & = & \mex{(s-2r)P,-s Q,s P,(s-2r)Q}.
\end{eqnarray*}
The cost of LHS and RHS are $4s$, resp., $2s+2|s-2r|$. As $0<r<s$, we deduce
$|s-2r|<s$,
i.e.,  RHS is better. 
\end{Proof}
\begin{Remark}
  The preceding argument has a nice geometric interpretation.  Both
  factorizations (or arcs) actually generate a parallel translation along the
  geodesic perpendicular to the initial tangent. So initial and terminal
  tangent are both perpendicular to the geodesic joining the initial and
  terminal point.
  \begin{center}
    \nobox{\includegraphics[height=50mm,bb=90 540 320 760]{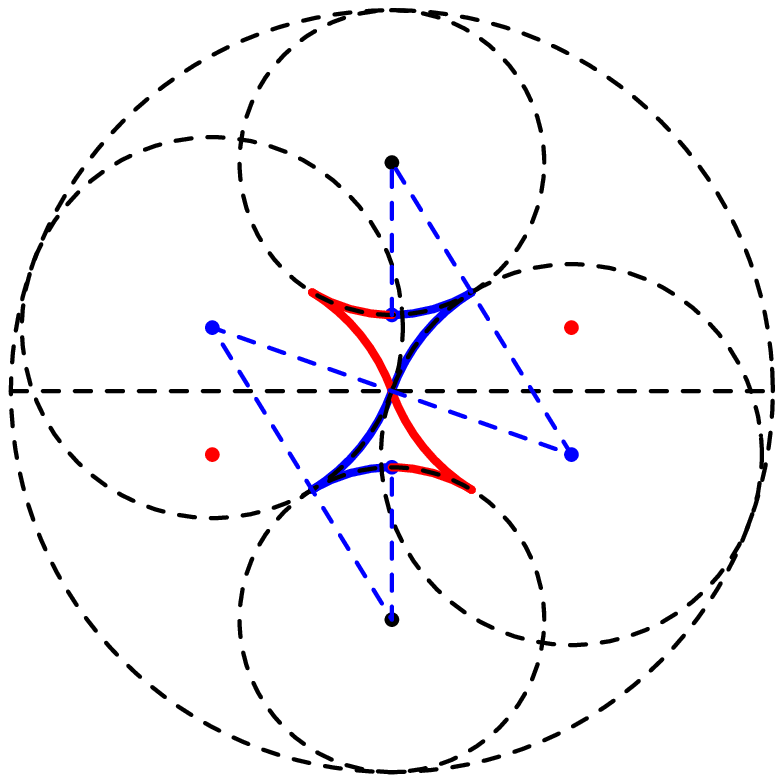}}
\quad 
\begin{minipage}[b]{5cm}\begin{center}
 $\mex{rP,-sQ,-sP,rQ}$\\[4pt] versus\\[4pt]
 $\mex{-r P,sQ,sP,-rQ}$.
\end{center}
\vspace*{1.5cm}
\end{minipage}
  \end{center}
\end{Remark}

A kind of second order variational argument allows us to prove the following
stronger result:
\begin{Proposition}\label{csp-variational-propo}
  Let $G$ be an arbitrary Lie group with Lie algebra~$\sL(2)$.
  Then for $s\in(0,1)$ and $t>0$ the factorization $\mex{s Q,-s P,-s Q,t P}$ is not optimal.
\end{Proposition}
\begin{Proof}
Fix $s\in(0,1)$ and consider $\gamma(t)=\mex{s Q,-s P,-s Q,t P}$.
Let $F:\R^3\to G$, $F(x)=\mex{x_1 P,x_2 Q,-x_2 P,-x_3 Q}$.
  Then $\gamma(0)=F(0,s,s)$. We claim that the following statements hold true:
  \begin{itemize}
  \item[(1)] The differential $dF(0,s,s)$ is invertible. Therefore there
    exist $t_0>0$ and a smooth (actually analytic) curve
    $x\colon(-t_0,t_0)\to\R^3$ such that $F(x(t))=\gamma(t)$ for all
    $|t|<t_0$.

    \item[(2)] We have $x_k(t)>0$ for $t\in(0,t_0)$, $k=1,2,3$.
      \item[(3)] Let $\delta(t)=3s+t-(x_1(t)+2x_2(t)+x_3(t))$.\\ Then
        $0=\delta(0)=\delta'(0)$, and $\delta''(0)>0$.
  \end{itemize}
  Once we have proved (1)--(3) it is clear that 
  $\mex{s Q,-s
  P,-s Q,t P}$ cannot be optimal because for small $t>0$ the map $F(x(t))$
  provides a better factorization.
  So let us prove (1)--(3). First we must compute the differential $dF(x)$.
  We observe:
  \begin{eqnarray*}
    dF(x)e_1&=& {F(x)}\, e^{x_3 \ad Q} e^{x_2\ad P}\,e^{-x_2\ad Q} P,\\
    dF(x)e_2&=& {F(x)}\, e^{x_3 \ad Q} e^{x_2\ad P}\, (Q-P) ,\\
    dF(x)e_3&=& {F(x)}\,(-Q)\\
    &=&d\lambda_{F(x)}\, e^{x_3 \ad Q} e^{x_2\ad P}\,e^{-x_2\ad P}(-Q).
  \end{eqnarray*}
 Thus $dF(x)={F(x)}\, e^{x_3 \ad Q} e^{x_2\ad P} M(x_2)$ with
$$
M(x_2)=\left[e^{-x_2\ad Q} P,\ Q-P,\ e^{-x_2\ad P} (-Q)\right].
$$
In terms of the basis $\setof{H,P,Q}$ we compute
$$
M(s)=
\begin{pmatrix}
     s &  0 & s  \\
     1 & -1 & s^2\\
  -s^2 &  1 & -1
\end{pmatrix},
\quad \det M(s)=-2s(s^2-1).
$$
Since $0<s<1$, Claim~(1) is proved. Next we use the Implicit Function Theorem
(IFT) to compute the derivate $\dot x(0)$. Differentiating $F(x(t))=\gamma(t)$
we obtain $dF(x(t))\dot x(t) = \dot\gamma(t)=d\lambda_{\gamma(t)}(\1) P$. Thus
\begin{displaymath}
  e^{x_3(t) \ad Q} e^{x_2(t)\ad P} M(x_2(t))\, \dot x(t) \equiv P.
\end{displaymath}
Since $x(0)=(0,s,s)$, we deduce  $\dot x(0)=M(s)^{-1} e^{-s\ad P} e^{-s\ad Q}
P$.
Simplification with \textsc{Mathematica} yields:
$$
\dot x(0)=
\left(
  \frac{s^2+2}{2(1-s^2)},\quad  -\frac{s^2}{2}, \quad \frac{s^2+2s^4}{2(-1+s^2)}
\right).
$$
Since $x_2(0)=x_3(0)=s>0$, continuity of $x(t)$ yields $x_2(t),x_3(t)>0$ for
sufficiently small $t>0$. Now $x_1(0)=0$ and $\dot x_1(0)>0$ because of
$1-s^2>0$,
hence $x_1(t)>0$ for small positive~$t$ follows, too. Thus we may assume
w.l.o.g. that $x_k(t)>0$ for $t\in(0,t_0)$, $k=1,2,3$ which proves Claim~(2).

Verifying $\delta(0)=0=\delta'(0)$ is straight forward. Finally, we must use
the (IFT) again to compute $\delta''(0)=-(\ddot x_1(0)+2 \ddot x_2(0)+\ddot
x_3(0))$.
Differentiating the identity
$$ 
M(x_2(t))\, \dot x(t) = e^{-x_2(t) \ad P} e^{-x_3(t) \ad Q} P
$$
we obtain
\begin{eqnarray*}
  M(x_2) \ddot x & = & - \dot x_2\, M'(x_2) \dot x
 -\dot x_2 e^{-x_2(t) \ad P} \ad(P) e^{-x_3(t) \ad Q} P
\\ &&{}
 -\dot x_3 e^{-x_2(t) \ad P}  e^{-x_3(t) \ad Q} \ad(Q) P.
\end{eqnarray*}
Simplification with \textsc{Mathematica} yields
$$
\ddot x(0) = 
\begin{pmatrix}
  -\frac{3 s^3(3+s^2)}{4(-1+s^2)^2}\\[1mm]
  -\frac{s^3}{4} \\[1mm]
  \frac{s^3(-1+5s^2+8s^4)}{4(-1+s^2)^2}
\end{pmatrix},\qaq
\delta''(0) = \frac{3s^3(2+s^2)}{2(1-s^2)}.
$$
For $0<s<1$ the denominator is positive, hence $\delta''(0)>0$ follows.
This proves Claim~(3) and finishes the proof of the proposition.
\end{Proof}

Now we are ready to prove Proposition~\ref{csp-propo}.

\begin{Proof}
  Consider a (CSP)-extremal with $5$ factors. Up to symmetry ($\Gamma$) we may
  assume that the first factor is  $\exp(rP)$.
  The switching pattern is clockwise or counterclockwise:
$$
\mex{r P,s Q,-s P,-s Q,t P},\quad\mbox{or}\quad \mex{r P,-s Q,-s P,s Q,t P},
$$
with $0< r,t < s < \sqrt{2}$. 
 Recalling that $\imath{\sigma_H}$ simply reverts the factors,
we see that the second pattern is transformed into the first pattern.

For $s\geq1$ Proposition~\ref{csp-s-ge-1} applies, showing  that this arc contains
a non-optimal subarc. For $s\in(0,1)$ we apply the previous proposition and obtain that 
$\mex{s Q,-s P,-s Q,t P}$ is not optimal.
This finishes the proof.
\end{Proof}

We conclude our discussion of (CSP)-extremals with one more observation:
\begin{Proposition}\label{csp-lagniappe-propo}
 Let $r>0$  and $s\in(1,\sqrt2)$. Then $\mex{rP,s Q,-sP}$ is not optimal.
 In particular, an optimal \textsc{(CSP)}-extremal with switching time~$>1$ has at
 most~$3$ factors. Conversely, if $C4a(r,s,t)$ or $C4c(r,s,t)$ is optimal,
 then $s\leq1$ and $r+2s+t\leq4$.
\end{Proposition}
\begin{Proof}
  We compute
 $\mex{rP,s Q,-sP}=\mex{x_1Q, x_2  P, x_3 Q}$ with
$$
x_2=-(s+r(s^2-1)),\quad x_1=-\frac{s^2}{x_2},\qaq x_3=\frac{rs}{x_2}.
$$
For $s\geq1$ and $r\geq0$ we have $x_2\leq-1$, $x_1\ge0$, and $x_3\le0$. Hence
$|x_1|+|x_2|+|x_3|=x_1-x_2-x_3$, and
$$
(2s+r)-(x_1-x_2-x_3) = \frac{r^2(s^2-1)(2-s^2)}{s+r(s^2-1)}.
$$
For $r>0$ and $s\in(1,\sqrt2)$ the righthand side is positive, hence
$\mex{rP,sQ,-sP}$ is not optimal.
\end{Proof}

\subsection{Alternating (ALT)-extremals}

While the results for optimal (CSP)-extremals are valid in any group $G$
with Lie algebra $\sL(2)$, some of the results for (ALT)-extremals
depend explicitly on~$G$. The crucial parameter is actually $|Z(G)|=:n$, the
cardinality of the center of $G$. We start with results that hold for any $G$.
Whenever we make explicit use of $G$ (resp.~$n$) we will indicate this clearly.

\begin{Proposition}\label{alt3-s3-propo}
  For $s\in(0,1)$ and $r>0$ the factorization $\mex{r P,-s Q,s P}$ is not
  optimal.  In particular, an optimal \textsc{(ALT)}-extremal with switching time
  $s\in(0,1)$ has at most $3$ factors. Conversely, if $A4(r,s,t)$ is optimal,
  then $s\geq1$.
\end{Proposition}
\begin{Proof}
  We compute
$\mex{(r P,-s Q,s P)} =\mex{x_1 Q,x_2 P, x_3 Q}$ with
$$
x_1=-\frac{s^2}{s+r(1-s^2)},\quad 
x_2=s+r(1-s^2),\quad
x_3=-\frac{rs}{s+r(1-s^2)}.
$$
For $s\in(0,1)$ and $r>0$ we deduce $x_2>0$ and $x_1,x_3<0$, so the second
factorization has cost $-x_1+x_2-x_3$. The cost difference is
$$
r+2s-(-x_1+x_2-x_3) = \frac{r^2\,s^2\,\left( 1 - s^2 \right) }
  {s + r\,\left( 1 - s^2 \right) }.
$$
For $s\in(0,1)$ and $r>0$ numerator and denominator are both positive, 
so $r+2s-(-x_1+x_2-x_3)>0$. Hence $\mex{r P,-s Q,s P}$ is not
  optimal.
\end{Proof}

\begin{Proposition}
  The factorizations 
$$
\mex{P,-Q,P,-tQ}\qaq\mex{\sqrt2 P,-\sqrt2 Q,\sqrt2 P,-\sqrt2 Q,tP}
$$ 
are not optimal for $t>0$.
In particular, an optimal \textsc{(ALT)}-extremal with $s=1$ has at most $4$
factors, and with $s=\sqrt2$ it has at most~$5$ factors.
\end{Proposition}
\begin{Proof}
 Let $w=\sqrt2$. An elementary computation yields
$$
\mex{P,-Q,P}=\matU,\quad
\mex{w P,-w Q,w P,-w Q}=-\matI.
$$
 Since $\sigma_U(U)=U$ and $\sigma_U(-\1)=-\1$, we get
 $\mex{-Q,P,-Q}=\mex{P,-Q,P}$ and 
 $\mex{w P,-w Q,w P,-w Q}=\mex{-w Q,w P,-w Q,w P}$.
 Hence
\begin{eqnarray*}
\mex{P,-Q,P,-t Q} 
&=&\mex{-Q,P,-(t+1)Q}, \\
\mex{w P,-wQ,w P,-w Q,tP}&=&\mex{-w Q,w P,-w Q,(t+w) P}.
\end{eqnarray*}
 In each case LHS and RHS  have equal cost, but for $t>0$ the righthand side is not an
 extremal, hence it cannot be optimal.
 We also observe that for all $r\in\R$ the identities
 \begin{eqnarray*}
&&\mex{(1-r)P,-Q,P,-rQ}\equiv\matU,\\
&&\mex{(w-r)P,-wQ,wP,-wQ,rP}\equiv-\matI
\end{eqnarray*}
hold true. Therefore an optimal (ALT)-extremal with $s=1$ and $4$ factors may
be replaced (at equal cost) by another extremal with~$3$ factors.
Similarly, for switching time~$s=\sqrt2$ it suffices to consider at most~$4$ factors.
\end{Proof}

The same technique as in the proof of Proposition~\ref{csp-variational-propo}
can be used to obtain:
\begin{Proposition}\label{alt4-s4-propo}
  For $1<s<\sqrt{2}$ the factorization $\mex{-s Q,s P,-s Q,t P}$ is not
  optimal.
  In particular,  optimal \textsc{(ALT)}-extremals with switching time
  $s\in(1,\sqrt{2})$ have at most $4$ factors. Conversely, if $A5(r,s,t)$ is
  optimal,
  then $s\geq\sqrt2$.
\end{Proposition}
\begin{Proof}
 Let $\gamma(t)= \mex{-s Q,s P,-s Q,t P}$ and consider $F\colon\R^3\to G$,
 $F(x)=\mex{x_1 P,-x_2 Q,x_2 P,-x_3 Q}$.  Then $\gamma(0)=F(0,s,s)$.
 We claim that $dF(0,s,s)$ is invertible and that there exists a $t_0>0$ and
 a smooth curve  $x\colon(-t_0,t_0)\to\R^3$ such that $x_k(t)>0$ for
 $t\in(0,t_0)$, $k=1,2,3$, and for $\delta(t)=3s+t-(x_1(t)+2 x_2(t)+x_3(t))$
 we have  $0=\delta(0)=\delta'(0)$, $\delta''(0)>0$.
Computing the differential of $F$ yields:
$$
dF(x)={F(x)}\, e^{x_3\ad Q} e^{-x_2\ad P} M(x_2)
$$
with
$$
M(x_2)=\left[ e^{x_2\ad Q} P,\quad -Q+ P,\quad e^{x_2\ad P}(-Q)\right].
$$
In terms of the basis $\setof{H,P,Q}$ we obtain 
$$
M(s)=
\begin{pmatrix}
  -s   &  0 & -s  \\
  1    &  1 & s^2 \\
  -s^2 & -1 & -1 
\end{pmatrix},
\quad \det M(s) = -2s(s^2-1).
$$
Differentiating $F(x(t))=\gamma(t)$ yields
$$
e^{x_3\ad Q} e^{-x_2\ad P} M(x_2) \dot x = P,\quad
M(x_2) \dot x = e^{x_2\ad P} e^{-x_3\ad Q} P.
$$
For $x(0)=(0,s,s)$ we obtain
$$
\dot x(0)=\left(\frac{s^2}{-2+2s^2},\ 1-\frac{s^2}{2},\ \frac{2-5s^2+2s^4}{-2+2s^2}\right).
$$
As $s>1$, $\dot x_1(0)>0$. Therefore $x_k(t)>0$ for $t>0$ small follows.
With $\delta(t)=3s+t-(x_1+2x_2+x_3)$ we have $\delta(0)=0=\delta'(0)$.
Differentiating once more we deduce
$$
M(x_2) \ddot x = -\dot x_2\, M'(x_2)\dot x 
+\dot x_2 e^{x_2\ad P} \ad(P) e^{-x_3\ad Q} P
-\dot x_3 e^{x_2\ad P} e^{-x_3\ad Q} \ad(Q) P.
$$
\textsc{Mathematica} yields:
$$
\ddot x(0)=\left(\frac{4 - 8\, s^2 - 3\, s^4 + 3\, s^6}{4\, 
    s\, {\left( -1 + 
              s^2 \right) }^2}, \frac{-{\left( -2 + s^2 \right) }^2}{4\, 
    s}, \frac{{\left( -2 + s^2 \right) }^2\, \left( 1 - 5\, s^2 + 8\, 
        s^4 \right) }{4\, s\, {\left( -1 + s^2 \right) }^2}
\right),
$$
and
$$
\delta''(0)=\frac{s\, \left( 4 - 8\, s^2 + 3\, s^4 \right) }{2 - 2\, s^2}=
\frac{\left( 
        s\, \left( 2 - s^2 \right) \, \left( -2 + 3\, 
            s^2 \right)  \right) }{2\, \left( -1 + s \right) \, \left( 
        1 + s \right) }.
$$
For $1<s<\sqrt{2}$ the numerator is positive and so is the denominator,
hence
$\delta''(0)>0$ as we claimed. This finishes the proof.
\end{Proof}

Up to this point all statements hold true in \textbf{any} group with Lie
algebra~$\sL(2)$. The last proposition in this subsection
makes very explicit use of $\SL(2)$. 

\begin{Proposition}\label{alt-sl2-a5-propo}
  Let $s>\sqrt{2}$ and $r,t\geq 0$. If $r+t > \frac{2s}{s^2-1}$, then
  the factorization $\mex{r P,-s Q,s P,-t Q}$ is not optimal in $\SL(2)$.
  In particular,
$$
  \begin{array}{lcll}
 s>\sqrt{3}       & \implies & \mex{-s Q,s P,-s Q}    & \mbox{not optimal  in $\SL(2)$,}\\
 s=\sqrt{3},\ t>0 & \implies & \mex{-s Q,s P,-s Q,tP} & \mbox{not optimal  in $\SL(2)$,}\\
 s\in(\sqrt{2},\sqrt{3}) & \implies &
  \mex{-s Q,s P,-s Q,s P} & \mbox{not optimal  in $\SL(2)$.}
  \end{array}
$$
 In particular, optimal \textsc{(ALT)}-extremals in $\SL(2)$ have at most~$5$
 factors, and if $A5(r,s,t)$ is optimal, then $s\in[\sqrt2,\sqrt3]$

\end{Proposition}
\begin{Proof}
  For $s\neq\pm1$ let $\displaystyle\mu=\frac{2s}{s^2-1}$, then
$$
 \mex{ \mu P,\,-s Q,s P} = 
 \mex{s Q,-s P,\,\mu Q}.
$$
Now assume w.l.o.g. that $r,t \leq \mu$. Then
\begin{eqnarray*}
\mex{r P,-s Q,s P,-t Q} &=&\mex{-(\mu-r) P,\mu P,-s Q,s P,-\mu Q,(\mu-t)Q}
\\
   & = & \mex{-(\mu-r)P,s Q,-s P,(\mu-t) Q}.
\end{eqnarray*}
The LHS has cost $2s+r+t$, the RHS has cost $2s+2\mu-(r+t)$. As $r+t>\mu$,
$2s+r+t>2s+\mu>2s+2\mu-(r+t)$, so LHS is not optimal.

For $s>\sqrt{3}$ we can apply the previous result with $r=0$ and $t=s>\mu$.
For $s=\sqrt{3}$ we observe $\mu=s=\sqrt{3}$, and
$$
\mex{s P,-s Q,s P,-t Q}=\mex{s Q,-s P,s Q,-t Q}=\mex{s Q,-s P,-(s-t) Q}.
$$
Comparing costs we obtain $3s+t$ (LHS) and $3s-t$ (RHS), so LHS is not
optimal for $t>0$.
\\[1mm]
Finally, for $s\in(\sqrt{2},\sqrt{3})$ we have $s^2-1>1$, hence $2s >
\frac{2s}{s^2-1}$.  Taking $r=t=s$ we therefore have $r+t=2s > 2\mu$. Hence
$\mex{s P,-s Q,s P,-s Q}$ is not optimal in $\SL(2)$, whence
optimal (ALT)-extremals cannot have $4+2=6$ factors.
\end{Proof}
The following picture in the hyperbolic plane shows geometrically where
the mysterious identity comes from.
\begin{center}
  \nobox{\includegraphics[height=50mm,bb=90 490 320 710]{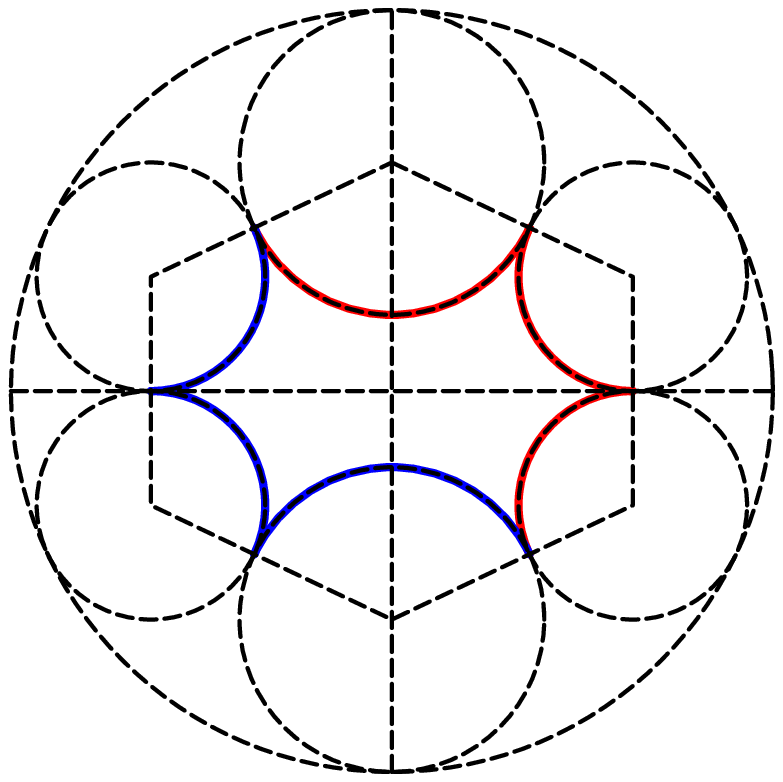}}
  \quad 
  \begin{minipage}[b]{7cm}
    \begin{center}
      $\mex{\mu P,-sQ,sP}=\mex{sQ,-sP,\mu Q}$.
    \end{center}
    \vspace*{1.8cm}
  \end{minipage}
\end{center}
It is noteworthy to mention that the bound on the number of factors of optimal
(ALT)-extremals always depends on the group. In particular we woul like to
stress that for the simply connected group $\SLt$ there is \textbf{no a
priori bound} on the number of factors.


\subsection{Singular (SSP)-extremals}

There are two reasons for non-optimality of (SSP)-extremals.
One of them is purely local, the other one is global in nature.
Looking at the graph from Figure~\ref{ssp-graph} we can say this as follows:
For (SSP)-extremals that switch between top and bottom line (of the graph) we get a bound
on the number of factors using purely local arguments.
An (SSP)-extremal that does not switch between top and bottom basically looks
like an (ALT)-extremal with switching time $2\sqrt2$ plus an interspersed $S$-arc.
A global argument (depending explicitly on the group)
is indispensable to bound the number of
factors for  these (SSP)-extremals.

 To save some typing throughout this subsection we let $w=\sqrt{2}$.
We recall that this is the ``switching'' time of the bang-arcs of an
(SSP)-extremal. The quotation marks indicate that a switch might be
\textbf{virtual}
because the singular arc between two identical bang-arcs may have zero length.

\begin{Proposition}\label{ssp-eqns-sl2}
Let $w=\sqrt{2}$. Then the following identities hold true in $\SL(2)$:
\begin{eqnarray}
    \mex{sP,-\frac2s\,Q}&=&-\mex{\frac2s\,Q,-sP},\qmb{if $s\neq0$,}
    \label{sl2:eq-1}
\\
    \mex{w P,-w Q} & = & -\mex{wQ,-wP},\label{sl2:eq-1w}
 \\
    \mex{s T,w P,-w Q} & = & \mex{w P,-w Q,-s T}, \quad s\in\R,\label{sl2:eq-2}
\\
\label{sl2:eq-4}%
\mex{r P, w Q,-w P} & = &\mex{ \frac{2}{r+w}Q,\, -(r+w)P,\, -\frac{r\,w}{r+w}
  Q},
 r\in\R,
\\
\mex{w P,w Q,-w P} & = & \mex{\frac w2 Q,\,-2 wP,\,-\frac w2 Q}.\label{sl2:eq-5}
  \end{eqnarray}
\end{Proposition}

The proof consists of nothing but elementary computations.
The next result is valid for any group. Basically it means that the crucial
(SSP)-extremals are those that look like (ALT)-extremals with switching
time~$2\sqrt2$ with an interspersed $S$-arc.

\begin{Proposition}\label{ssp-3w-propo}
  The factorization $\mex{w P,w Q,-w P}$ is not optimal in any group.
  In particular, if an optimal \textsc{(SSP)}-extremal contains $5$ $B$-arcs,
  then for each $S$-arc the two adjacent $B$-arcs must be equal.
\end{Proposition}
\begin{Proof}
  In view of the last item of the previous proposition we observe that
$$
\mex{w P,w Q,-w P}=\mex{\frac w2 Q,-2w P,-\frac w2 Q}.
$$
Both factorizations have equal cost $3w$, but the righthand side cannot be
optimal because it does not come from an extremal! Indeed, (ALT)-extremals are
impossible because the switching pattern is circular. (CSP)-extremals are
impossible because they have switching time $s\in(0,\sqrt2)$ whereas here the
switching time is $2\sqrt2$. Finally, (SSP)-extremals are impossible, too,
because for an (SSP)-extremal with middle arc~$\exp(-2 wP)$ the third arc would
have to be $\exp(t Q)$ with $t\geq0$.
\end{Proof}

A similar argument allows to prove the following stronger result:
\begin{Proposition}\label{ssp-top-bottom-bound}
  For $r,s>0$ the factorization $\mex{rP,sT,wQ,-wP}$ is not optimal.
  In particular, an optimal \textsc{(SSP)}-extremal that switches
  between top and bottom has at most $5$ factors.
\end{Proposition}
\begin{Proof}
  Let $r,s>0$. Using Proposition~\ref{ssp-eqns-sl2}(\ref{sl2:eq-2}) and
  (\ref{sl2:eq-4}) we obtain
$$
\mex{rP,sT,wQ,-wP}=\mex{\frac{2}{r+w}Q,\, -(r+w)P,\, -\frac{r\,w}{r+w} Q,\,-sT}.
$$
Both sides have equal cost $r+s+2w$, but the righthand side is not an
extremal if $r>0$ and $s>0$, hence it cannot be optimal.
\end{Proof}

One more argument is needed.

\begin{Proposition}\label{s7a-sl2-propo}
 In $\SL(2)$ $\mex{rP,sT,wP,-wQ,-wQ,wP}$ is not optimal.
 In particular, optimal \textsc{(SSP)}-extremals in $\SL(2)$ have at most~$7$ factors.
\end{Proposition}
\begin{Proof}
 Due to Eqn.~(\ref{eq-2}) we may
 shift the $S$-arc  to the right end,  and it suffices to show that
 $\mex{rP,wP,-wQ,-wQ,wP}$
 is not optimal for $r>0$.

Next we use the  seemingly obscure Eqn.~(\ref{sl2:eq-1})
whose geometric meaning becomes clear in the discussion of (HRSCP) (cf. proof
of Proposition~\ref{alt-psl2-propo}), namely
$$
\mex{s P,-\frac2s\, Q}=-\mex{\frac2s\, Q,-s P},
$$ which is
valid for all $s\neq0$. We consider $s=r+w>w$. Then $\frac2s < \frac2w = w$.
Let $\epsilon=w-\frac2{w+r}>0$. 
Then we obtain:
$$
\mex{(r+w)P,-wQ,-wQ,wP}=\mex{\frac2{r+w} Q,-(r+w)P,-\epsilon Q,-wP,wQ}.
$$
Both factorizations   have equal
cost. But the second factorization  cannot be optimal because it
does not come from an extremal! Indeed, the switching pattern is neither
alternating nor circular. Besides the switching times $r+w,\epsilon,w$ are not
equal. Hence $\mex{(r+w) P,-w Q,-wQ,wP}$ cannot be optimal in $\SL(2)$.
An (SSP)-extremal with $8$ factors  will always contain a subarc
of the above (non-optimal) form.
\end{Proof}

Our last task in this subsection is to show that the previous propositions
imply that the singular extremals listed in Table~\ref{sl2-family-tab} are
sufficient.
Up to five factors we obtain the factorization maps $S3P,S4Q,S4P,S4Q,S5P,S5Q$,
and $S5a$. Due to Proposition~\ref{ssp-top-bottom-bound} we know that optimal
(SSP)-extremals with more than five factors look like (ALT)-extremals with
switching time~$2\sqrt2$ and an interspersed $S$-arc.
In view of Eqn.~(\ref{sl2:eq-2}) the $S$-arc may be the second or third arc in
the product.  But if the number of factors
is even, we may revert the order of the factors (i.e., apply $\imath\sigma_H$), if
necessary, so that we can also assume that the $S$-arc \textit{is} the second
arc. Therefore it suffices to consider a single factorization map (say,
$S6,S8,\dots$) for an even
number of factors, but it is necessary to distinguish between, say,
$S7a=\mex{rP,\frac s2\,T,\dots}$ and 
$S7b=\mex{-rQ,wP,\frac s2\,T,\dots}$, resp. $S9a,S9b,\dots$, for an odd number
of factors. As $S7a$-extremals are not optimal in~$\SL(2)$, this naming
convention may look surprising, but it perfectly makes sense if one also wants
to consider the simply connected group~$\SLt$. 
Since $S8$-extremals cannot be optimal, we see that the
family $\mathcal{F}$ given in Table~\ref{sl2-family-tab} exhausts all possibilities.
Thus  the proof of Theorem~\ref{sl2-family-thm} is finished.


\section{Conclusion}

\begin{figure}[htbp]
  \centering
  \begin{tabular}{cc}
    \includegraphics[height=58mm]{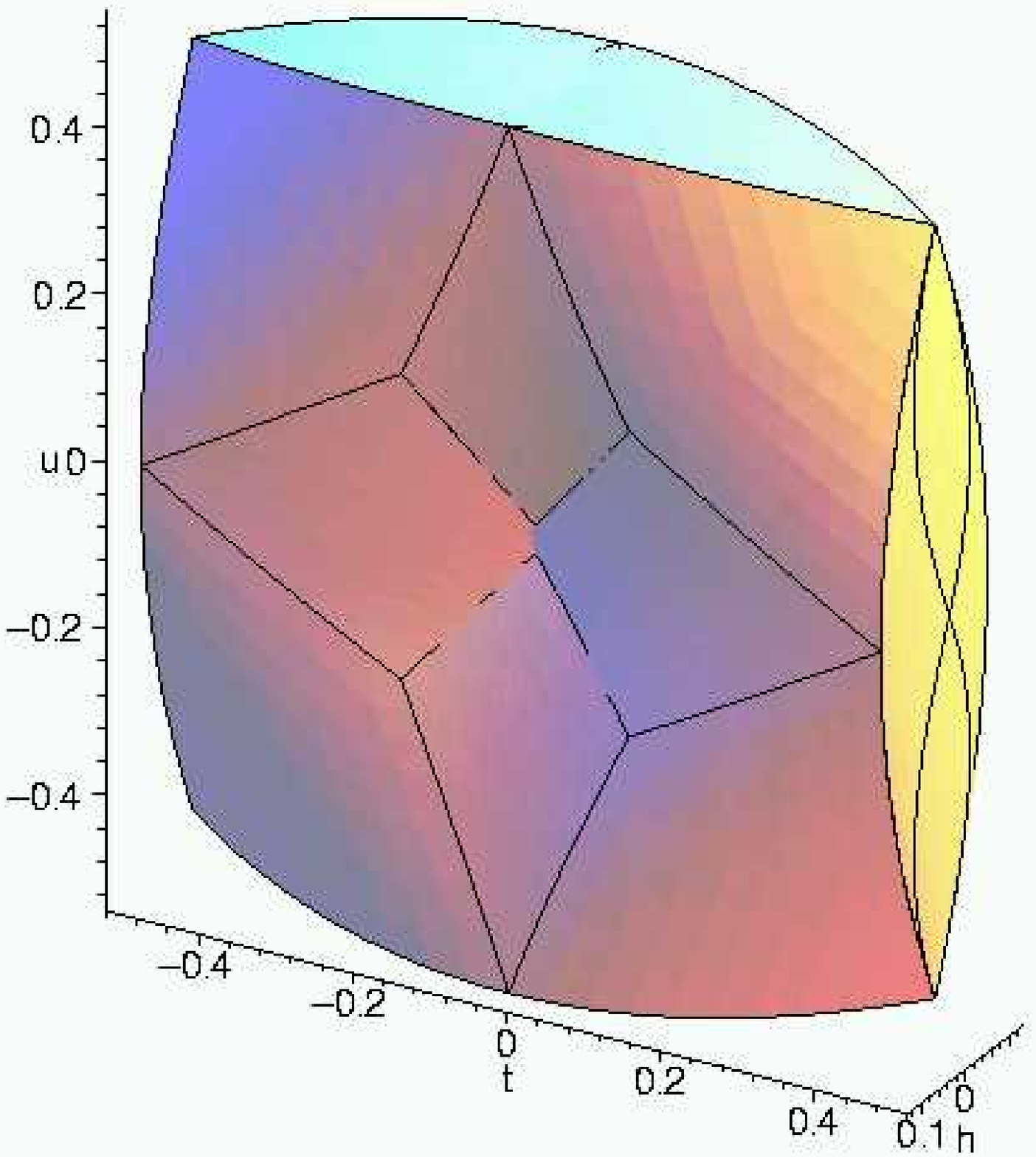} &
    \includegraphics[height=58mm]{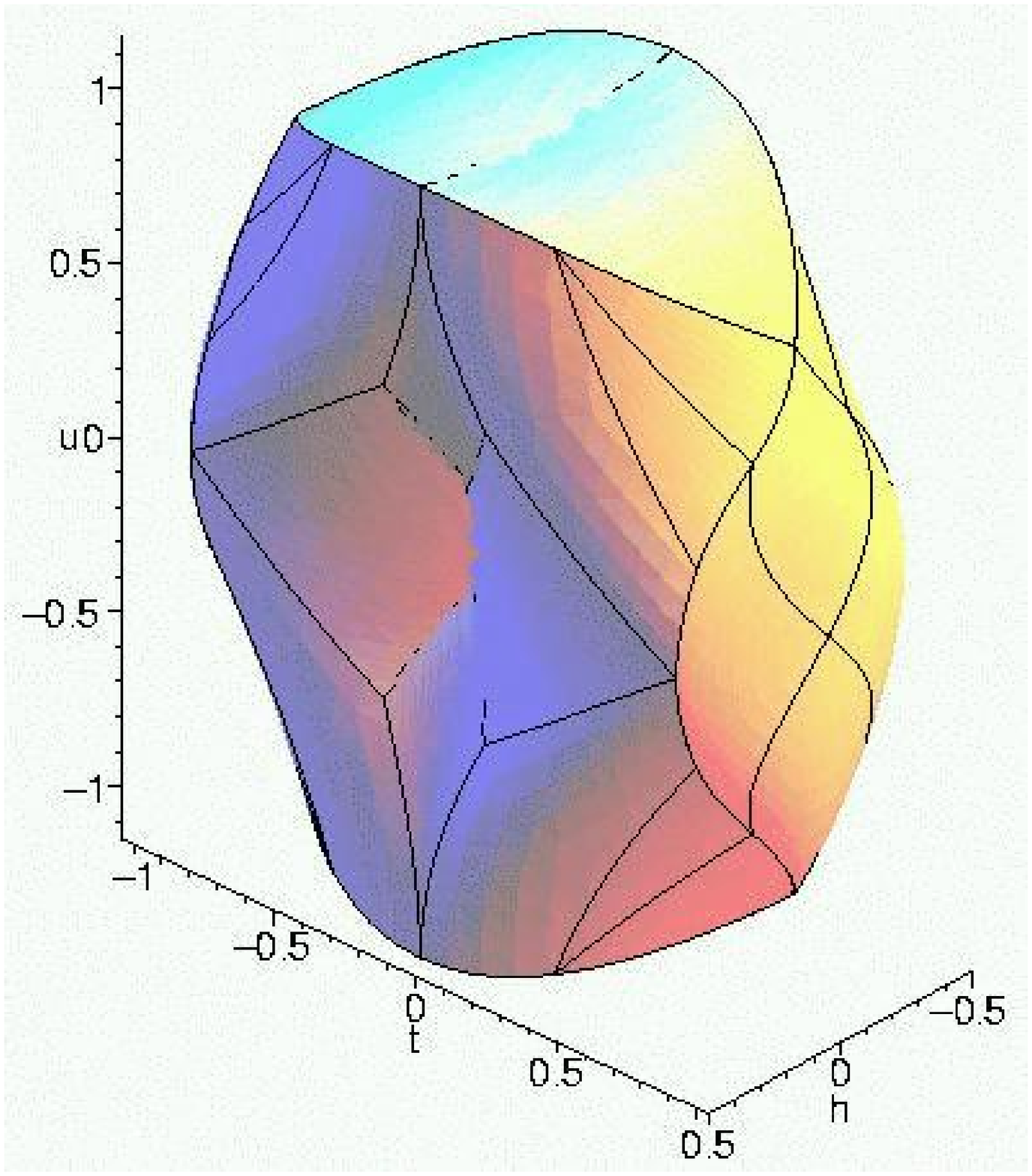} 
    \\ $\bbS(1)$ & $\bbS(2)$ \\[2mm]
    \includegraphics[height=58mm]{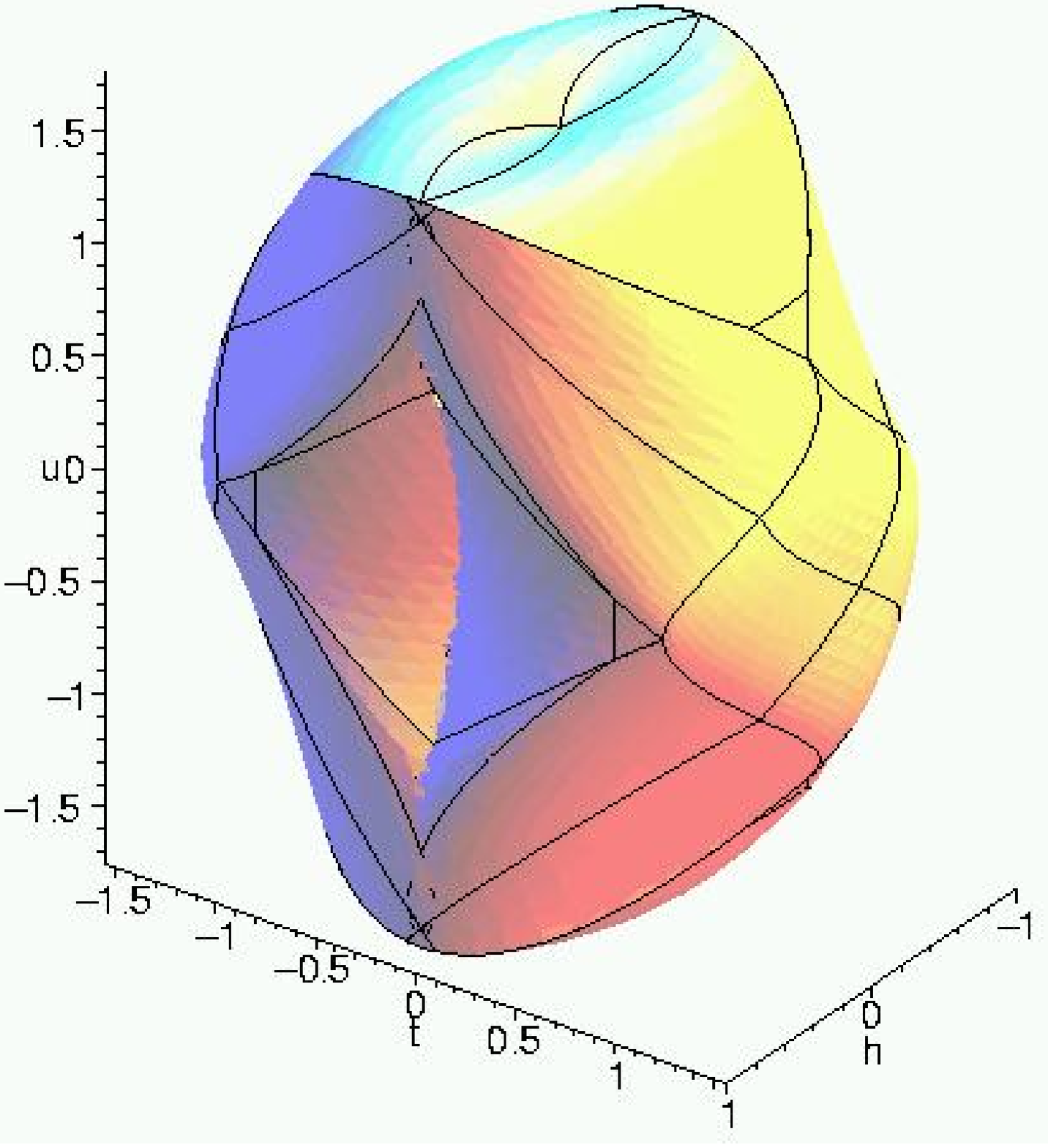} &
    \includegraphics[height=58mm]{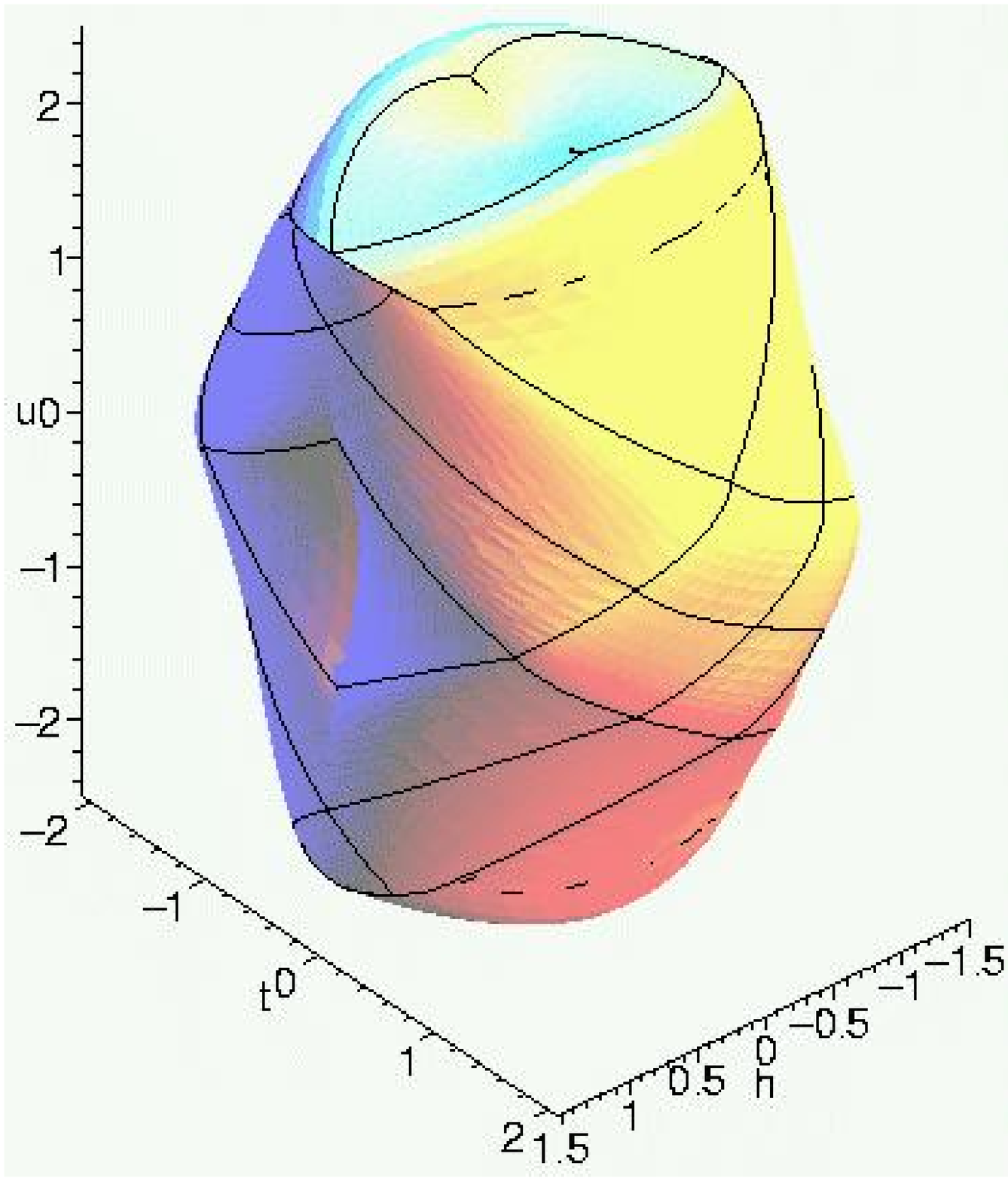} 
    \\ $\bbS(3)$ & $\bbS(4)$ \\[2mm]
    \includegraphics[height=58mm]{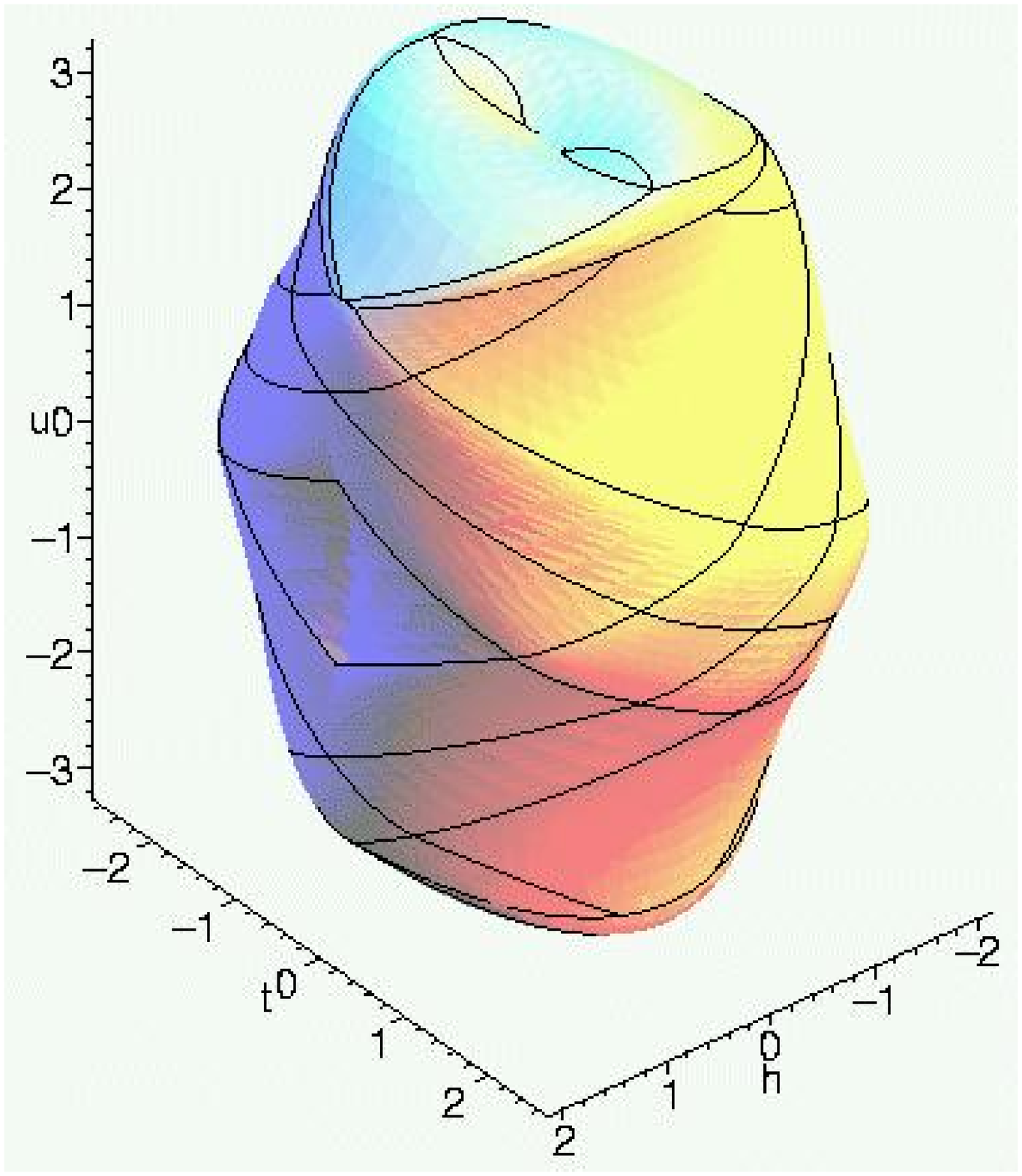} &
    \includegraphics[height=58mm]{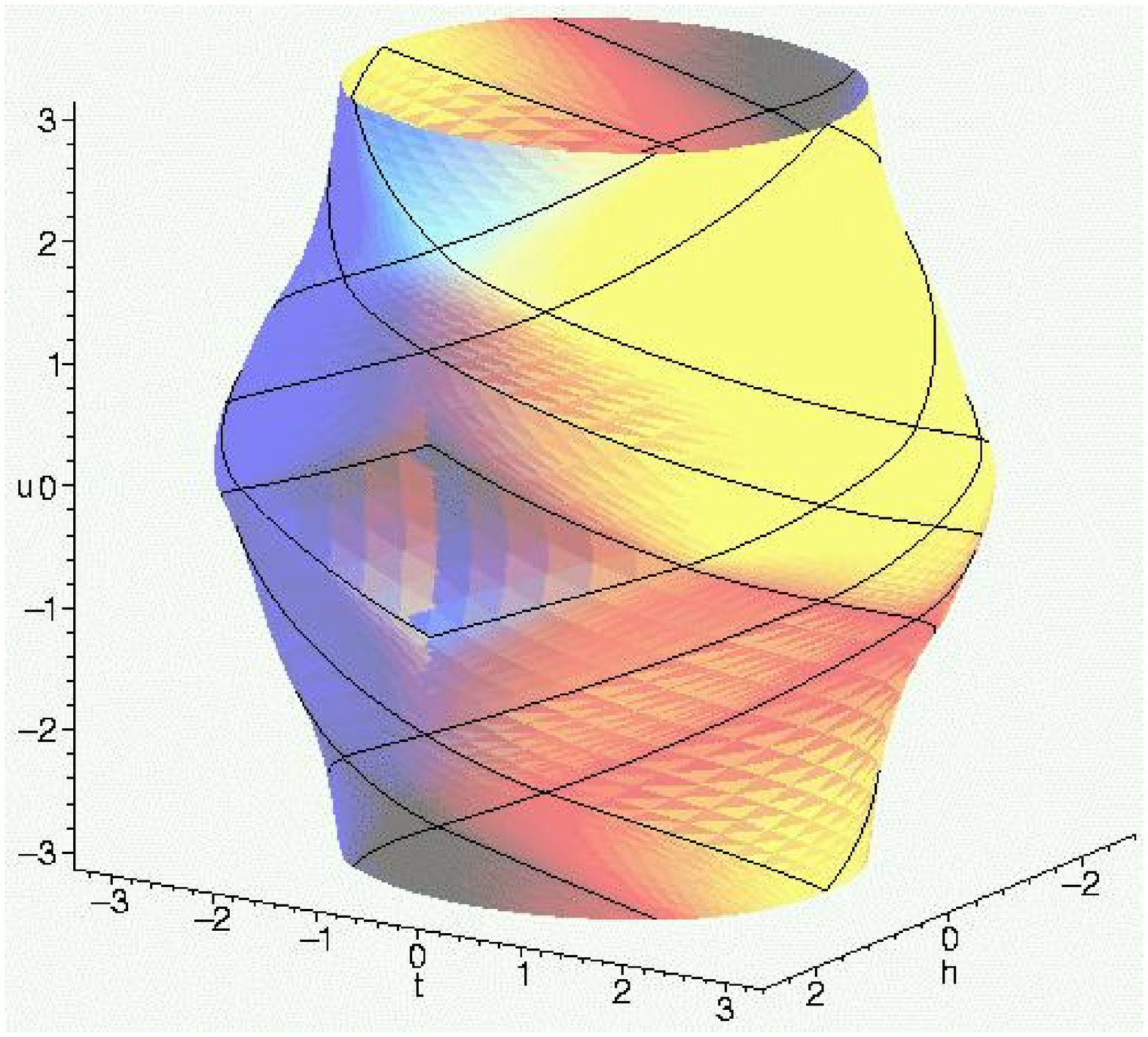} 
    \\ $\bbS(c)$ for $c=3\sqrt3\approx5.196$& $\bbS(c)$ for $c=\frac{32}7\sqrt2\approx 6.465$
  \end{tabular}
\caption{Metric spheres $\bbS(c)$ for some values of $c$.}
  \label{fig:metric-spheres}
\end{figure}

Having proved that the families $\cF,\cF_1,\cF_2$ \textbf{are} sufficient,
it is natural to ask if they are minimal with this property. The answer is
affirmative
in the sense that for every $f\in\cF$ there exist $(r,s,t)\in\dom(f)$ such
that $f(r,s,t)$ is optimal, and no $\tilde f\in\cF\setminus{f}$ allows us to reach
the same endpoint at equal cost. Rather than providing a list of such cases,
we include pictures of metric spheres $\bbS(c)=\setof{g:\cT(g)=c}$ for some
values of~$c$, cf. Fig.~\ref{fig:metric-spheres}.
To generate these pictures we used the parametrization of~$\SL(2)$, resp.,
its simply connected covering group $\SLt$ described in Appendix~A.
To understand the pictures it suffices to know that $\SL(2)$ is identified
as a set with $\R H+\R T+(-\pi,\pi] U\subseteq\sL(2)\cong\R^3$.
In this parametrization the symmetries $\sigma_H,\sigma_T,\sigma_U$ are 
180-degree rotations around the $H$-, resp., $T,U$-axes. Inversion is simply,
$\imath(X)=-X$, and $\imath\sigma_H$ is reflection in the $TU$-plane.
The horizontal plane $\R H+\R T$ corresponds to symmetric, positive definite
elements of~$\SL(2)$ while the vertical segment $(-\pi,\pi] U$ corresponds to the
circle group $\exp(\R U)=\exp([-\pi,\pi]U)\cong\SO(2)$.

For $c=1$ the relevant maps are $A3,C3,C4a,C4c,S3P$, and $S3Q$.
The four $A3$-patches make up the top and bottom part.
The \textit{thin} sides with the figure eight curves consist of $S3P$- and $S3Q$-patches,
the $S3P$-patch is inside the figure eight. 
The ``flat'' sides consist of eight patches: four $C3$- and two
$C4a$- resp., $C4c$-patches. The $C4a$-patches connect (horizontally) to
an~$S3Q$-patch (because $C4a(0,s,0)=\mex{-sQ,-sP}$) whereas the $C4c$-patches
connect (vertically) to $A3$-patches (because of $C4c(0,s,0)=\mex{sQ,-sP}$).
One can \textit{see} that the $C3$- and $S3$-patches \textit{fit in} perfectly
while the $A3$- and $C4$-patches have nontrivial intersections.

The shape of the sphere~$\bbS(1)$ still reminds us of the fact that the generating
set $\cU=\conv(\pm P,\pm Q)$ is a flat square.

For $c=2$ we also have $S4P$- and $S4Q$-patches emanating from the
$S3$-patches, and connecting them to the $A3$- and $C3$-patches, respectively.
One can also see how larger portions of $C4a$- and $C4c$-patches intersect.

For $c=3$ tiny triangular patches appear, the $S5P$- and $S5Q$-patches.
One also sees how $A4$-extremals arise on the top surface. The $S5a$-patch
is also a tiny triangle adjacent to the $S4P$-patches.

For $c=4$ the $C4$-patches have completely disappeared
 (cf. Proposition~\ref{csp-lagniappe-propo}).
 The only change
compared to $\bbS(3)$ is, of course, that the $S5$- and $A4$-patches are
much larger, hence better visible.

The last two pictures show $\bbS(c)$ for especially critical values of~$c$.
For $c=3\sqrt3$ the sphere touches the plane $\R H+\R T+\pi U$ for the first
time, $A3$-extremals make up only a very thin portion of the top while $A4$- and
$A5$-patches are clearly visible. Nevertheless the $A5$-patches are quite
small which means that the set of matrices for which $A5(r,s,t)$ is
\textit{the} optimal factorization is quite small and close to $-\1\in\SL(2)$.

For $c=\frac{32}7\,\sqrt2$ (ALT)-extremals have disappeared (which means they
are above $\pi U$) and all patches come from (SSP)-extremals. Except for the
$S4Q$- and $S5Q$-patches there are no further intersections. So in particular,
$S3P,S4P,S5P,S5a,S6,S7b$-patches perfectly \textit{fit} together, only the
portions where the vertical coordinate exceeds~$\pi$ are chopped off in
$SL(2)$. They would still be present in~$\SLt$, though.

Finally these pictures also make clear that the set
$\Fix(\imath\sigma_T)
\subseteq\SL(2)$
is a cut locus for the problem. In the pictures $\Fix(\imath\sigma_T)$ is the
plane $\R H+\R U$. An extremal that hits this plane transversally will lose
its global optimality. A glance at Table~\ref{sl2-family-tab} reveals
that this applies to $A3,A5,C4a,S4Q,S5Q$-extremals; it does not apply to
$S5a,S7b,C4c$-extremals because the latter are $\imath\sigma_T$-invariant.

\bigskip
Before we give an outlook on  generalizations and future work
 we  now give a proof for Theorem~\ref{other-2slip-thm}. 

\begin{Proposition}\label{other-2slip-propo}
  Let $S^1,S^2\in\sL(2)$ with $\det(S^1)=\det(S^2)=0$, and assume that
  $[S^1,S^2]\neq0$.
  Let $\tcU=\conv(\pm S^1,\pm S^2)$.
  Then there exist $\mu>0$ and a $g_0\in\SL(2)$ 
  such that $\tcU=\mu\Ad(g_0)\cU$.
  In particular, we have
  $\cT_{\tcU}(g)=\mu^{-1} \cT_{\cU}(g_0^{-1}gg_0)$ for all $g\in\SL(2)$.
\end{Proposition}
\begin{Proof}
  We recall that $\setof{X\in\sL(2):\det(X)=0,X\neq0}$ is the boundary of a
  Lorentzian double cone $C\cup -C$ (minus the vertex), cf. Figure~\ref{square-pic}.
  W.l.o.g. we may assume that $S^1,S^2$ lie on the upper part, say~$C$.
  Finding $\mu$ and $g_0$ is a $3$-step procedure:
\\[2mm]
 \textbf{1. Rotation:} There exist $u\in[0,\pi]$ and $\rho_1>0$
  such that $e^{u\ad U}S^1=\rho_1\,P$. Let $g_1=\exp(uU)$.
\\[1mm]
 \textbf{2. Shearing:} We first observe that (in $HTU$-coordinates:)
$$
e^{\tau\ad(P)}(-Q)=\left( -\tau,\, \frac{\tau^2-1}2,\, \frac{\tau^2+1}2\right)
=\frac{\tau^2+1}2
\left(-\frac{2\tau}{1+\tau^2},\ \frac{\tau^2-1}{\tau^2+1},1\right).
$$
Since $(1+\tau^2)^{-1}\,(2\tau,\,1-\tau^2)$ parametrizes the unit circle
(except for  the point~$(0,1)$), we see that
$(0,\infty)\,e^{\R\ad(P)}(-Q)=C\setminus(0,\infty)P$.
Hence for $X=\Ad(g_1)S^2$ we can find $\tau^*\in\R$ and $\rho_2>0$
such that $e^{\tau^*\ad(P)}X=-\rho_2\,Q$.
We let $g_2=\exp(\tau^*P)$ and
observe that $\Ad(g_2)P=P$.
\\[1mm]
\textbf{3. Hyberbolic rotation:}
  Let
  $h=\frac14\log(\rho_2/\rho_1)\in\R$, so 
$e^{2h}\rho_1=e^{-2h}\rho_2=\sqrt{\rho_1\rho_2}=:\mu>0$.
 Then $e^{h\ad(H)}\rho_1 P=\mu P$, and $e^{h\ad(H)}(-\rho_2 Q)=-\mu Q$.
 Let $g_3=\exp(hH)$.
 Then $\Ad(g_3g_2g_1) S^1=\mu P$ and $\Ad(g_3g_2g_1)S^2=-\mu Q$.
\\[2mm]
Finally let $g_0=(g_3g_2g_1)^{-1}$.
Then $\tcU=\mu\,\Ad(g_0)\,\cU$. Hence
$$
\cT_{\tcU}(g) = \cT_{\mu \Ad(g_0)\, \cU }(g)
=\frac1{\mu}\, \cT_{\Ad(g_0)\, \cU}
=\frac1\mu\, \cT_{\cU}(\Conj{g_0}^{-1} (g))
=\frac1\mu\, \cT_{\cU}(g_0^{-1} g\,g_0).
$$
So the claim of Theorem~\ref{other-2slip-thm} follows
with $\sigma=\Conj{g_0}$ and $\lambda=\mu^{-1}$. 
\end{Proof}

Thus all claims made in the introduction have been proved by now.
For an arbitrary $2$-slip system with symmetric slip rates 
we now know how to find the dissipation distance and geodesics.

\subsection*{Generalizations and future work}

It is clear that two different types of generalizations are of interest,
namely more slip systems and passage to dimension~$d=3$. In a forthcoming
paper we will  treat the $2D$-hexagonal lattice, i.e., we consider slips
along the sides of an equilateral triangle. This is of practical
interest because such systems arise in reality. 

Observing that in the problem analyzed in this paper the optimal controls are
piecewise constant, it is natural to ask whether this is true for general
slip systems, or, for general polytopes, i.e., $\cU=\conv(X_1,\dots,X_m)$
for $X_1,\dots,X_m\in\sL(d)$. The following example shows that the answer is
negative for arbitrary polytopes.

\begin{Example}
  In $\sL(2)$ let $X_1=H+2Q$, $X_2=H-2Q$, $X_3=-H+2P$, $X_4=-H-2P$, and
  $\cU=\conv(X_1,\dots, X_4)$. Then $\cU$ is a simplex and $0\in\inter\cU$. For
  the polar $\cQ=\setof{p:\cH(p)\geq-1}$ we obtain
  $\cQ=\conv(p^1,p^2,p^3,p^4)$ with $p^1=(1,1,-1)$, $p^2=(1,-1,1)$,
  $p^3=(-1,1,1)$, $p^4=(-1,-1,-1)$.  Now consider
  $p(t)\in\conv(p^1,p^2)$, then (PMP:2) implies
  $u(t)\in\conv(X_3,X_4)=-H+[-2,2]P$. 
 Let $\h=\R H+\R P$ denote the subalgebra generated by $X_3,X_4$.
 The associated subgroup is $\left\langle \exp\h\right\rangle=
\setof{\twomatrix{a}{b}{0}{a^{-1}}: a>0,\ b\in\R}$.
 Now $\left\langle (e^{\ad\h})^*\right\rangle$ leaves the line
 $\aff(p^1,p^2)=(1,0,0)+\R(0,1,-1)$ invariant.
 Take $u(t)=-H+v(t) P$ with measurable $v\colon[0,\infty)\to[-2,2]$,
 $g(t)$ the corresponding trajectory with $g(0)=\1$, and
 $p(t)=\Ad(g(t))^*\frac12(p^1+p^2)$.
 Then $g(t)\in\left\langle \exp\h\right\rangle$ and $p(t)\in\aff(p^1,p^2)$ for
 all $t$. By continuity, there exists $t^*>0$ such that
 $p(t)\in\conv(p^1,p^2)$ for $t\in [0,t^*]$.
 Hence $(g,p,u)$ is an extremal in~$[0,t^*]$.
\end{Example}

The geometric reason for the occurence of this degeneracy in the previous
example is that the edge $\conv(p^1,p^2)$ lies in a single coadjoint
orbit---the one-sheeted hyperboloid $\setof{p:\mathcal{C}(p)=-1}$. Considering
the dual picture, we already observed that the opposed edge~$\conv(X_3,X_4)$
lies in a $2$-dimensional subalgebra of~$\sL(2)$. This is a very special
situation that will never occur for slip systems because the $2D$-subalgebras
in~$\sL(2)$ are precisely the hyperplanes tangent to the double cone.
Therefore we state:

\begin{Conjecture}
   Let $\cU=\conv(X_1,\dots,X_m)\subseteq\sL(2)$ with $\det(X_j)=0$.
   Then optimal controls are piecewise constant. 
\end{Conjecture}

A careful look at our analysis of the flow rule actually suggests a much
stronger conjecture. For $X\in\sL(2)\smo$ with $\det(X)=0$, the induced flow
$(p,t)\mapsto p\,e^{t\ad(X)}$ on the opposed face 
$\opp(X)=\setof{p:\left\langle p,X\right\rangle=-1}$ is (up to symmetry)
the same flow as the one we analyzed on p.\pageref{switching-pattern-howto}f.
All flow lines are parabolas, like in Figure~\ref{face-flowpic}.
Consequently, an affine line, say $\ell\subseteq\opp(X)$ will be tangent to at
most one of these. This suggests that we get at most one singular control
for each edge of $\cQ$. At a vertex, say $p$, of~$\cQ$ the stabilizer $\g_p$
is one-dimensional. 
Hence $\g_p\cap\cU\cap\setof{X:\left\langle
    p,X\right\rangle=-1}$ is singleton, or empty.
 Since the polytopes $\cU$ and $\cQ$ are
(combinatorial) duals of each other, we can formulate the following

\begin{Conjecture}
 Consider $2$-dimensional slip systems, so $\cU\subseteq\sL(2)$ is a convex polytope,
 $\cU=\conv(S^1,\dots,S^m)$ with $\det(S^\alpha)=0$.
 Assume that $0\in\inter(\cU)$.
 Let $f_0,f_1,f_2$ denote the number of vertices, edges, and faces
 of the polytope~$\cU$.
 Then optimal controls are piecewise constant, and the number of possible
 values is at most $f_0+f_1+f_2$.
\end{Conjecture}
In other words, for each vertex of~$\cQ$ (resp., face of~$\cU$) and each edge
of~$\cU$, we get at most one singular control.
In order to obtain rigorous proofs for these conjectures,
Proposition~\ref{other-2slip-propo} and our discussion of the flow rule for
the square lattice are the appropriate tools 
(cf. propositions~\ref{switch-propo-u2} and~\ref{switch-propo-t2}).
For if we consider an arbitrary (but fixed) edge of~$\cQ$, we see that
up to automorphism, we may assume that the two active slip-systems are either $-P,Q$, or $-P,-Q$,
as in Proposition~\ref{switch-propo-u2} and~\ref{switch-propo-t2}.
Hence the arguments given there
apply, whence switches must occur except for some isolated points.

\medskip
Formulating reasonable conjectures for  $3D$-slip systems seems to be
more difficult. Of course one can hope that optimal controls must be
piecewise constant, but we do not have such strong evidence as in the
$2D$-case. In particular it is not yet clear, what types of singular controls
one has to expect. 

\subsection*{Computational tools}

Programming and debugging are usually time consuming tasks, therefore 
we will make some  tools available on the web.
In particular, we will provide {\sc Mathematica} notebooks with procedures
for finding the factorization cost and optimal factorizations in~$\SL(2)$.
The pictures of the metric spheres were generated with {\sc Maple}. But since
the numerical computations are too slow, a more efficient C-program
produces a file defining the necessary plot data structures.
An interested reader should follow the links starting at URL\\
 \texttt{http://www.mathematik.uni-stuttgart.de/mathA/lst1/mittenhuber/}


\section*{Appendix A: The group \boldmath{$\SLt$}}

\renewcommand{\thesection}{A}

\label{sec:appendix-sl2t}

In this appendix we provide a summary of facts about the  simply
connected group $\SLt$. It is well-known that $\SLt$ cannot be represented
(faithfully) as a matrix group (subgroup of $\GL(n)$ for some $n\in\N$).
When analyzing $\SL(2)$ and related groups the following functions are ubiquitous:
$$
\begin{array}{cclcl}
  C(z) & = & 
\displaystyle\sum_{n=0}^\infty \frac{z^n}{(2n)!}
&=&
\displaystyle
\begin{cases}
  \cosh(\sqrt z), & \mbox{if $z\geq0$},\\
  \cos(\sqrt{-z}), & \mbox{if $z<0$},
\end{cases}
\\ 
\rule{0pt}{38pt}%
S(z) & = &\displaystyle \sum_{n=0}^\infty \frac{z^n}{(2n+1)!}
&=&
\displaystyle
\begin{cases}
  \frac{\sinh(\sqrt z)}{\sqrt z}, & \mbox{if $z>0$},\\  1, & \mbox{if $z=0$},\\
  \frac{\sin(\sqrt{-z})}{\sqrt{-z}}, & \mbox{if $z<0$}.
\end{cases}
\end{array}
$$
The ubiquity of these functions is due to the fact that $C$ and~$S$ describe
the (matrix) exponential function $\exp_{\SL(2)}:\sL(2)\to\SL(2)$:
$$
\sum_{k=0}^\infty \frac1{k!} X^k=C(-\det(X)\,) \id+S(-\det(X))\, X\qmb{for every $X\in\sL(2)$.}
$$
Due to a lack of letters we let $T(z)=\frac{S(z)}{C(z)}$, assuming that
 the function $T(\cdot)$ will not be confused with the matrix $T=P+Q$.


\subsection*{The covering map and a local inverse}

As a set we identify $\SLt$ with the Lie algebra $\sL(2)$. Using the basis
$\setof{H,T,U}$ and writing $X\in\sL(2)$ as $X=hH+tT+uU$ with $h,t,u\in\R$,
the covering map $f\colon\sL(2)\to\SL(2)$ is defined as
\begin{equation}
  \label{eq:sl2-cover}
  \cover(hH+tT+uU)
=C(h^2+t^2)
\begin{pmatrix}
  \cos(u) & \sin(u)\\ -\sin(u) & \cos(u)
\end{pmatrix}
+S(h^2+t^2)
\begin{pmatrix}
  h & t \\ t & -h
\end{pmatrix}.
\end{equation}

Hilgert and Hofmann introduced this map in~\cite{hofmann-old-new} to
parametrize the group $\SLt$.
Since $f$ is an analytic covering, it allows to define a group
operation~$\circ\colon\sL(2)\times\sL(2)\to\sL(2)$ such that
$$
\cover(X\circ Y)=\cover(X)\,\cover(Y)\quad\mbox{for all $X,Y\in\sL(2)$.}
$$
From now on we will always identify the group $\SLt$ with $(\sL(2),\circ)$.

\bigskip
The identity $\cover(X+2k\pi U)=\cover(X)$ for $X\in\sL(2)$, $k\in\Z$ is obvious.
Conversely, for $g\in\SL(2)$ the preimage  $\cover^{-1}(g)$ has the form $X+2\pi\Z
U$ for some suitable $X$.
Next let $\cE=\R H+\R T$ denote the set of symmetric matrices in $\sL(2)$.
Then the restriction $f\colon \mathcal{E}+(-\pi,\pi]\,U\to\SL(2)$ is
injective. An inverse of this restricition is easily obtained. For
\begin{eqnarray*}
  g & = & 
  \begin{pmatrix}
    a & b \\ c & d
  \end{pmatrix} = \cover(X)
\implies X  =  hH+tT+uU \qmb{with}
\\
u & = & \arg( (a+d)+i(b-c) ) \in (-\pi,\pi],\\
\twovec{h}{t} & = & \frac1{2S(\rho)} \twovec{a-d}{b+c},\ \mbox{where}\ %
\rho  =  \arcosh\left(\sqrt{\frac{(a+d)^2+(b-c)^2}{4}}\right).
\end{eqnarray*}
We denote this local inverse map simply
$\icover:\SL(2)\to\mathcal{E}+(-\pi,\pi]\,U$.
Practically this means that we obtain the group $\SL(2)$ from $\SLt$
simply by taking the $U$-coordinate mod~$2\pi$.

\subsection*{Advantages and disadvantages}

This parametrization of $\SLt$ has \textbf{advantages} and \textbf{disadvantages}.
The \textbf{main disadvantage} is that the group operation is much more complicated
than ordinary matrix multiplication. We will give an explicit expression for
$X\circ Y$, but this expression is practical mainly for numerical computations.
The main \textbf{advantages}, on the other hand, are the rotational symmetry and the fact
that every one-parameter subgroup $\exp(\R X)\subseteq\SLt$ lies in 
$\R X+\R U\subseteq\SLt$.

A particular advantage for the problem (OCP) is that in the above
parametrization all maps from the symmetry group
$\tilde\Gamma$ are simply linear maps on the vector space~$\sL(2)$.
For $X=hH+tT+uU$ one quickly verifies that:
\begin{eqnarray*}
  \cover(-X) & = & \cover(X)^{-1},\\
  \cover(hH-tT-uU)&=&\sigma_H(\cover(X)),\\
  \cover(-hH+tT-uU)&=&\sigma_T(\cover(X)),\\
  \cover(-hH-tT+uU)&=&\sigma_U(\cover(X)).
\end{eqnarray*}
The first equation shows that inversion in $\SLt$ is simply
$\tilde\imath(X)=-X$. 
Similarly, the second equation implies
 that
$\tilde\sigma_H:=(X\mapsto HXH)\colon\SLt\to\SLt$ is an automorphism of~$\SLt$.
So the 180-degree rotations around the $H,T,U$-axes, respectively, are
all group automorphisms.

\subsection*{Group multiplication}

Multiplication and conjugation with elements from $\R U$ is easily
obtained---this is the rotational symmetry we already mentioned.
\begin{eqnarray*}
  (uU)\circ X & = & e^{\frac u2\ad U}X + uU, \\
  X\circ (uU) & = & e^{-\frac u2\ad U}X + uU, \\
  (uU)\circ X\circ(-uU) & = & e^{u\ad U} X,
\end{eqnarray*}
for all $X\in\sL(2)$, $u\in\R$. In particular,
$$
(k\pi U)\circ X =X\circ(k\pi U),\
(2k\pi U)\circ X=X+2k\pi U,\qmb{for all $k\in\Z$}.
$$
In order to find $X_1\circ X_2$ for abitrary $X_1,X_2\in\sL(2)$, we use
the observation that $\mathcal{E}\circ \mathcal{E}\subseteq
\mathcal{E}\times(-\pi/2,\pi/2)$,
cf.~\cite[Lemma~1]{dimi-glob-sl2}.
Hence
$$
X_1,X_2\in\mathcal{E}\implies X_1\circ X_2=\icover( \cover(X_1)\cover(X_2) )
$$
 Writing $X_i=h_iH+t_iT+u_iU$, $i=1,2$,
and observing $(-u_1 U)\circ X_1,X_2\circ(-u_2 U)\in\mathcal{E}$, we
obtain
\begin{eqnarray*}
  X_1\circ X_2 & = & (u_1U)\circ(-u_1 U)\circ X_1\circ X_2\circ(-u_2
  U)\circ(u_2 U)
\\
 & = & (u_1U)\circ \Big(
 \underbrace{((-u_1 U)\circ X_1)}_{\in\mathcal{E}}\circ
 \underbrace{(X_2\circ(-u_2U))}_{\in\mathcal{E}}
 \Big)
\circ(u_2 U)
\\
 & = & (u_1U)\circ \icover\Big(\cover(-u_1
  U)\cover(X_1)\cover(X_2)\cover(-u_2U) \Big)
 \circ(u_2 U).
\end{eqnarray*}
The last expression is as explicit as can be, but of course its usefulness
is mainly restricted to numerical computations.

\subsection*{One-parameter groups}

Let $\expt\colon\sL(2)\to\SLt$ denote the exponential function (of the group
$\SLt$). Our first observation is that the covering map $\cover(X)$ coincides
with the matrix exponential function if $X$ is either symmetric or
skew-symmetric. Hence
$$
\expt(X)=X,\qaq \expt(\R X)=\R X\qmbq{for all $X\in\mathcal{E}\cup\R U$.}
$$
It turns out that $\exp(\R X)\subseteq\R X+\R U$ holds true for all
$X\in\sL(2)$. More precisely, if $X_0=h_0H+u_0U$, $X_0\neq0$,
then (cf.~\cite{neeb-sl2})
$$
\exp(\R X_0)=\setof{hH+uU: u_0\tanh(h)=h_0\sin(u)}.
$$
Qualitatively there are two different cases (for $X_0\neq0$):
 \begin{align*}
   \expt(\R X_0)&= 
   \setof{ uU+\artanh\left(\frac{h_0}{u_0}\sin(u)\right)H: u\in\R},&
   \mbox{if $|u_0|>|h_0|$,}
   \\
   \expt(\R X_0)&=
   \setof{ hH+\arcsin\left(\frac{u_0}{h_0}\tanh(h)\right)\,U : h\in\R},&
   \mbox{if $|u_0|\leq|h_0|$.}
\end{align*}
The other one-parameter groups are obtained via rotation around the $U$-axis.
\begin{figure}[htbp]
  \centering
\nobox{\includegraphics[width=70mm]{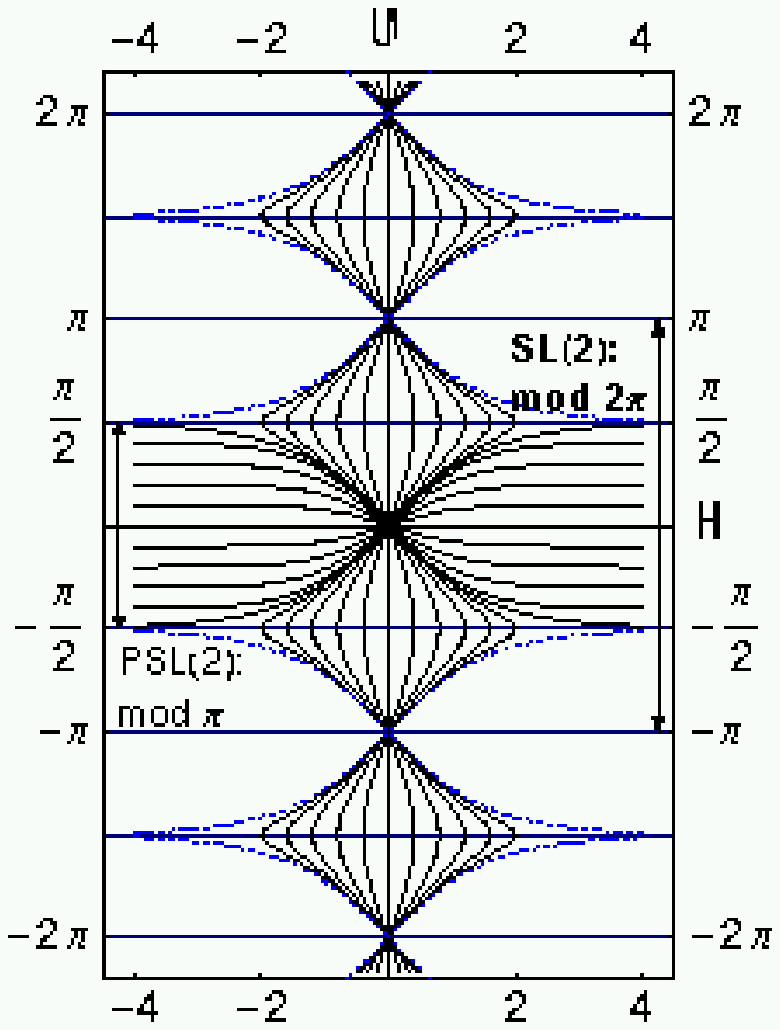}}\quad
\raisebox{10mm}{\includegraphics[width=50mm]{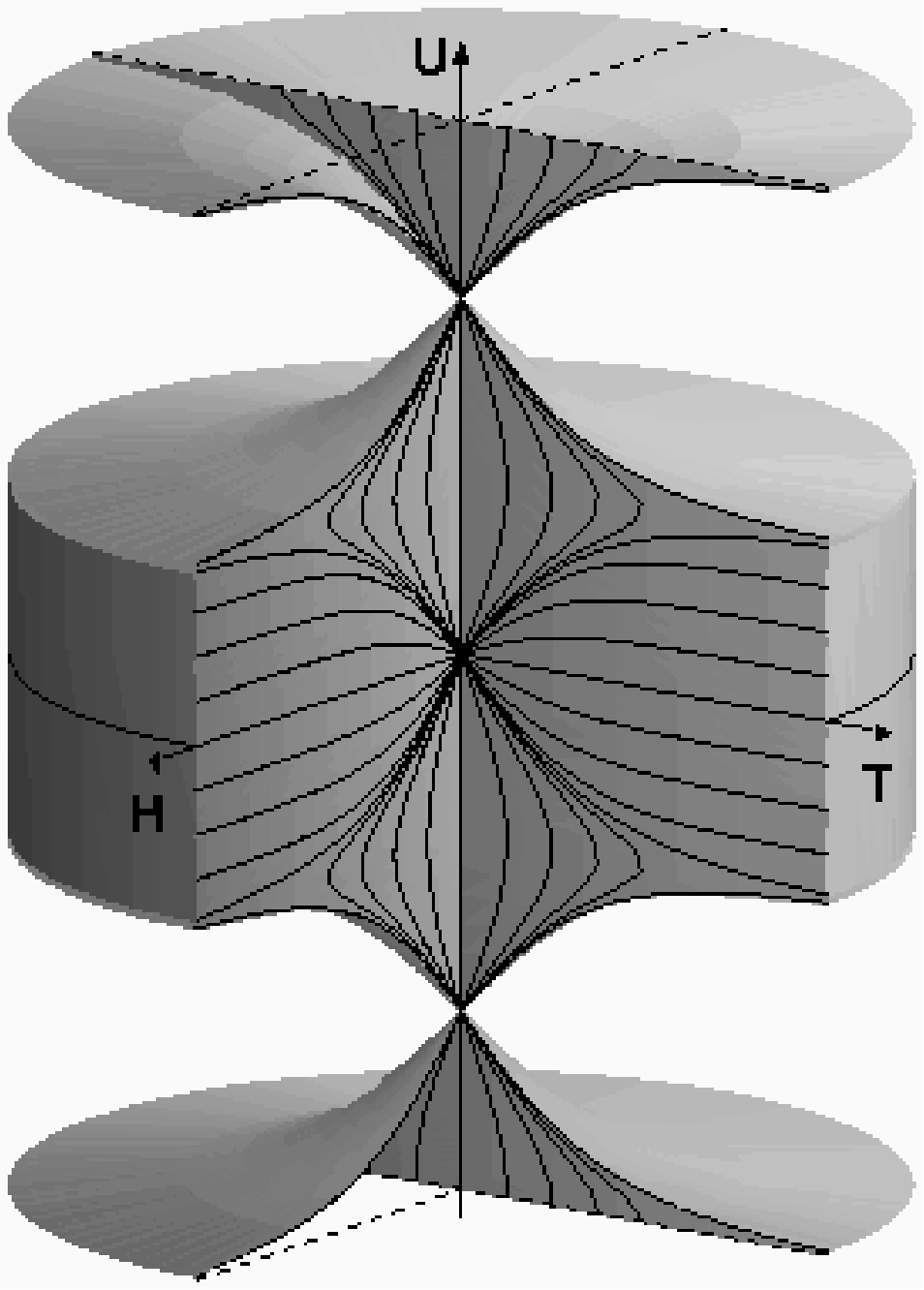}}
  \caption{One parameter groups in the plane $\R H+\R U$}
  \label{fig:sl2-exp-image}
\end{figure}
Figure~\ref{fig:sl2-exp-image} shows one-parameter groups in the $HU$-plane.
The group $\SL(2)$ is obtained by reading the picture modulo $2\pi U$.
The dashed curves indicate the boundary of the complement of the image of the
exponential function.
An explicit expression for $\expt\colon\sL(2)\to\SLt$ is
available, too, cf.~\cite{dimi-sl2,dimi-glob-sl2}.
Recalling that we defined
$$
T(z)=\frac{S(z)}{C(z)}=
\begin{cases}
  \frac{\tanh(\sqrt z)}{\sqrt z}, & z>0,\\
  \frac{\tan(\sqrt{-z})}{\sqrt{-z}}, & z<0,\\
  1, & z=0,
\end{cases}
$$
we observe that
$$
\dom(T)=\R\setminus\setof{-\frac{n^2\pi^2}4 :  \mbox{$n\in\N$, $n$ odd}}.
$$
Now let $X=\rho_0E+u_0U\in\sL(2)$ with
$E=\cos\alpha\,H+\sin\alpha\,T\in\mathcal{E}$, and
let $k(X)=\rho_0^2-u_0^2$.
Then $\exp(X)=\rho\, E+u\,U$ with
\begin{eqnarray*}
\rho & = & \arsinh(\rho_0\,S(k(X)) ),\\
u & = &
\begin{cases}
  \arctan( u_0 T(k(X)) ), & \mbox{$k(X)\geq0$},\\
\sign(u_0)\sqrt{-k(X)}, &
\mbox{$k(X)\not\in\dom(T)$,}\\
  \arctan( u_0 T(k(X)) )+\sign(u_0)\left\lfloor\frac12+\frac{\sqrt{-k(X)}}{\pi}\pi \right\rfloor, & \mbox{otherwise.}
\end{cases}
\end{eqnarray*}
Special cases of particular interest are
\begin{eqnarray*}
  \expt\left(tP\right) & = &
 \expt\left(\frac t2 (T+U) \right) =
 \arsinh\left(\frac t2 \right) T +
 \arctan\left(\frac  t2\right) U.
\\
  \expt(tQ) & = &
 \expt\left( \frac t2 (T-U) \right) =
 \arsinh\left( \frac t2 \right) T 
-\arctan\left(\frac  t2\right) U.
\end{eqnarray*}



\section*{Appendix B: The hyperbolic Reeds-Shepp-Car}

\renewcommand{\thesection}{B}

As we already mentioned (cf.~p.\pageref{hdp-psl2-system}) the Hyperbolic Dubins
Problem (HDP) can be considered as an optimal control problem on~$\PSL(2,\R)$.
Hence solving (HRSCP) is equivalent to finding time-optimal paths of
$$
\dot \gamma(t)=\gamma(t) u(t),\quad u \in\pm\conv(P,Q)
\mbox{ a.e., }\gamma\in\PSL(2).
\eqno{\rm(HRSCP)}
$$
At this point it is not clear, if optimal arcs exist for (HRSCP).
We must convexify the set of admissible control values, i.e.,
we must pass to $\conv(\pm P,\pm Q)$.
The only difference between this Convexified Hyperbolic Reeds-Shepp-Car
Problem (CHRSCP) and the problem (OCP) is that the first one evolves on
$\PSL(2)$ and  the second one on $\SL(2)$.
From the characterization of extremals and the elimination of (U/2)-extremals
we deduce that (HRSCP) is indeed solvable.

By definition, $\PSL(2)=\SL(2)/\{\pm\1\}$. Since the adjoint representation
$\Ad$ has kernel $\ker\Ad=\{\pm\1\}$, we can write $\Ad\colon\SL(2)\to\PSL(2)$
to denote the quotient map.

\begin{Proposition}
  Let $\mathcal{T}_{\PSL(2)}$ denote the factorization cost in $\PSL(2)$. Then
$$
\mathcal{T}_{\PSL(2)}(\Ad(g))= \min\setof{\mathcal{T}(g),\mathcal{T}(-g)}
\qmb{for all $g\in\SL(2)$}. 
$$
\end{Proposition}
This is due to the fact that for a trajectory $\gamma(t)$ in $\SL(2)$
its  projection $\Ad(\gamma(t))$ is a trajectory (for the same control $u(t)$)
in $\PSL(2)$. So ``$\leq$'' follows immediately. Conversely, if
$\eta(t)\in\PSL(2)$ is an optimal path from $\eta(0)=\Ad(\1)$ to
$\eta(t^*)=\Ad(g)$,
 then $\eta(t)$ lifts to a trajectory $\gamma(t)\in\SL(2)$.
Since $\eta(t^*)=\Ad(g)$, $\gamma(t^*)\in\setof{\pm g}$,
hence ``$\geq$'' follows, too.

\smallskip
A practical consequence of this observation is that the sufficient families
listed in Table~\ref{sl2-family-tab} are also sufficient for $\PSL(2)$.
In particular, there always exists an optimal path in~$\HH^2$ with at most $6$
factors, resp. $5$ switches. Of course our goal is to show that we can drop
some maps from Table~\ref{sl2-family-tab}, so we find smaller sufficient
families for $\PSL(2)$. Briefly, we will show that one may
drop~$A5,S5a,S6,S7b$ from $\mathcal{F}$. So optimal paths for (HRSCP) need
at most $5$ pieces, resp. $4$ switches.

Another important observation is that inversion of the
factorization maps is much nicer in $\SL(2)$ because $\exp(tP)=\id+tP$ is
linear in $t$ whereas $\Ad(\exp tP)=e^{t\ad P}$ is quadratic in~$t$.
Thus finding optimal paths for given boundary data is preferably done 
working in $\SL(2)$.

\subsection*{The missing comparison arguments}

We need two new arguments, one for (ALT)- and one for (SSP)-extremals.
The first argument also shows where the, perhaps obscure,
identity used in the proof of Proposition~\ref{s7a-sl2-propo} came from:
%
\begin{figure}
\centering
  \includegraphics[height=50mm,bb=90 490 320 710]{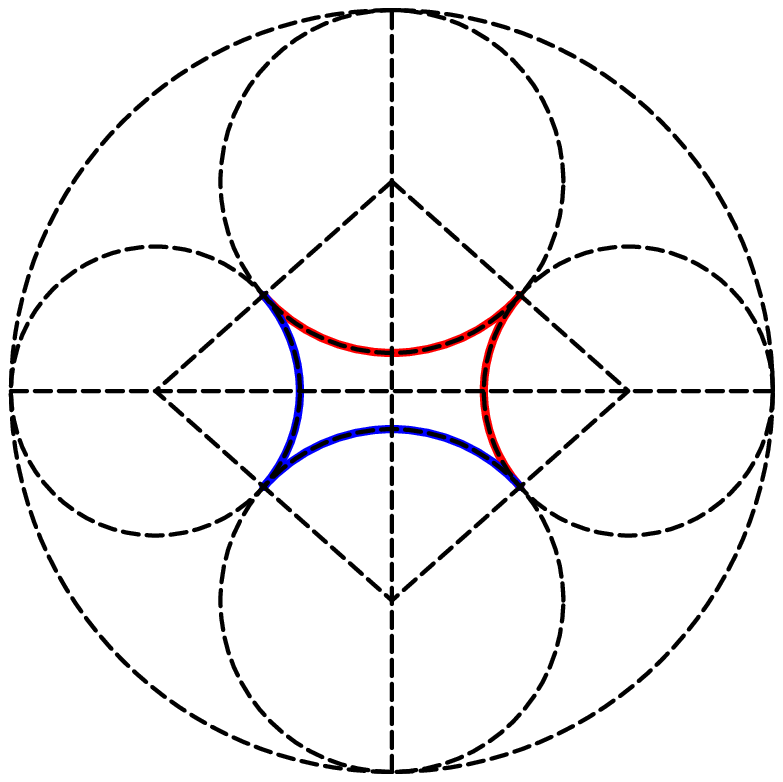}
  \quad 
  \begin{minipage}[b]{7cm}
    \begin{center}
      $\mex{\frac2s P,-sQ}=\mex{s Q,-\frac2s\,P}$
      \\[4pt] in $\PSL(2)$,
      \\[2mm] cf. Proposition~\ref{alt-psl2-propo}
    \end{center}
    \vspace*{1.8cm}
  \end{minipage}
\end{figure}
\begin{Proposition}\label{alt-psl2-propo}
  Let $s>1$ and $r,t>0$ such that $r+t>\frac2s$.
  Then the factorization $\mex{r P,-s Q,t P}$ is not
  optimal in $\PSL(2)$.
 In particular,
$$
\begin{array}{rcll}
s>1      & \implies & \mex{s P,-s Q,s P} & \mbox{not optimal in $\PSL(2)$,} \\
s>\sqrt2 & \implies & \mex{s P,-s Q}     & \mbox{not optimal in $\PSL(2)$,} \\
r>0      & \implies & \mex{r P,-\sqrt2\, Q, \sqrt2 P} & \mbox{not optimal in $\PSL(2)$.}
\end{array}
$$
 Furthermore, optimal \textsc{(ALT)}-extremals in $\PSL(2)$ have at most $4$ factors,
 and  if $A4(r,s,t)$ is optimal in $\PSL(2)$, then $s\in[1,\sqrt2]$.
\end{Proposition}
\begin{Proof}
 Let $\mu(s)=2/s$.  We already observed that 
$$
\mathbb{M}_{\SL(2)}\left({\mu P,-s Q}\right) 
=-\mathbb{M}_{\SL(2)}\left({s Q,-\mu P}\right),
$$
hence $\mex{\mu P,-s Q}=\mex{s Q,-\mu P}$ in $\PSL(2)$.
Thus $\mex{r P,-s Q,t P}=\mex{-(\mu-r) P,s Q,-(\mu-t)P}$. Comparing costs
we find $r+s+t$ for LHS and $s+2\mu-(r+t)$ for RHS. If $r+t>\mu$, then
$s+r+t>s+\mu>s+2\mu-(r+t)$. So RHS is better than LHS.

For  $s>1$ we choose $r=t=s$ and obtain $r+t=2s>2>\frac2s$.
For  $s>\sqrt2$ we choose $r=0$, $t=s$, and obtain $r+t=s>\frac2s$.
Finally, $s=\sqrt{2}=t$ and $r>0$ also implies $r+t>\frac2s$.
\end{Proof}

The second argument is similar, it eliminates $S5a$-extremals.
\begin{Proposition}
  For $r>0$ the factorization 
  $\mex{rP,\sqrt2\,P,-\sqrt2\,Q}$ is not optimal in $\PSL(2)$.
  In particular, $S5a(r,s,t)$ is not optimal for $r>0$,
  and optimal (SSP)-extremals in $\PSL(2)$ have at most~$5$
  factors.
\end{Proposition}
\begin{Proof}
Let $w=\sqrt2$ and $s=r+w>w$. Then $\frac2s<\frac2w=w$.
Let $\varepsilon=w-\frac2{r+w}>0$. In the previous proof we already observed
that $\mex{s P,-\frac2s Q}=\mex{\frac2s Q,-sP}$ holds true in $\PSL(2)$. Hence
\begin{eqnarray*}
\mex{rP,wP,-wQ}&=&\mex{sP,-wQ}= \mex{sP, -\frac2s\,Q,\, -\varepsilon Q} \\
 &=&\mex{\frac2s\,Q, -sP,-\varepsilon Q}.
\end{eqnarray*}
All  these factorizations have equal cost because of $s+\frac2s+\varepsilon=r+2w$.
But the last one cannot be optimal because it does not come from an
extremal.
Indeed, the switching pattern $Q\vdash  -P\vdash -Q$ is \textbf{circular}, but
the switching time is $s=r+\sqrt2>\sqrt2$. And we know that (CSP)-extremals
have switching time less than~$\sqrt2$.
The second claim is obvious because of $S5a(r,s,t)=\mex{rP,wP,-wQ,-sT,-tQ}$.
Finally, if an (SSP)-extremal has six factors, then it contains at most two
singular arcs. Therefore one can always find a subarc like $\mex{rP,wP,-wQ}$,
or $\mex{rP,sT,wQ,-wP}$.
\end{Proof}


\subsection*{Comparison with the euclidean case}

The hyperbolic and the euclidean problem bear similarities, but they also differ in
some respect. The euclidean case has a degeneracy which reflects the geometric
fact that the sum of the angles in a triangle equals $\pi$ in euclidean
geometry. In the hyperbolic case the sum is strictly less than~$\pi$, this
basically
accounts for the fact that we could eliminate the (U/2)-extremals by hand.
Consequently, an optimal path for (CHRSCP) is also admissible for (HRSCP),
and (HRSCP) always has a solution (for all possible boundary data).

Table~\ref{psl2-family-tab} shows a sufficient family for (HRSCP), the symmetry groups
are the same as in Table~\ref{sl2-family-tab}, so we do not list the groups
here. Instead, the last column contains a reference to~\cite{suss-car}, namely
where the corresponding paths appear, resp. where they are eliminated.

It is interesting to observe that the hyperbolic case in some sense swaps
the role of the two types of paths corresponding to $S5P$- and
$S5a$-factorizations. In $\HH^2$ $S5a$ is not optimal and $S5P$ yields some
optimal arcs whereas in the euclidean setting the paths corresponding to~$S5P$
are not optimal (cf.~\cite[Lemma 12,Fig.~9]{suss-car}) but the paths
corresponding to $S5a$ belong to the sufficent family (\cite[Theorem~8, Item~7]{suss-car}).

\begin{table}[htbp]
  \centering
  \begin{tabular}{|c|c|l|l|l|} \hline
    Type & Map & Domain & Remark & $\R^2$ (Sussmann/Tang) \\ \hline
    ALT & $A3$  &                 & $2$ cusps    & \begin{tabular}{l} degenerate case\\
     (LLTV)
   \end{tabular}
    \\
        & $A4$  & $s\in[1\sqrt2]$ & \textbf{new} & not needed \\ \hline
    CSP & $C3$  & &           & [8.3] \\ \hline
        & $C4a$ & & $2$ cusps & [8.4] \\
        & $C4c$ & & $1$ cusp  &       \\ \cline{2-5}
    SSP & $S3P$ & &           & [8.2] \\
        & $S3Q$ & &           &       \\ \cline{2-5}
        & $S4P$ & &           & [8.6] \\
        & $S4Q$ & &           &       \\ \cline{2-5}
        & $S5P$ & & \textbf{new} & \begin{tabular}{l}not optimal in $\R^2$\\
          cf.~Lemma~12, Fig.9
        \end{tabular}
        \\ \cline{2-5}
        & $S5Q$ & &           &  [8.5] \\    \hline
        & $S5a$ & & not optimal in $\HH^2$ & [8.7] \\ \hline
  \end{tabular}
  \caption{A sufficent family for $\PSL(2)$}
  \label{psl2-family-tab}
\end{table}


\subsection*{Visualization}

The problem (HRSCP) is also useful for visualizaton purposes.
Every trajectory $\gamma(t)\in\PSL(2)$ (or $\SL(2)$) yields a path
$\zeta(t)$ in the open unit disc~$\mathbb{D}$. One may assume w.l.o.g.
that the controls $P$ and $Q$ correspond to forward left- and right-turns,
while $\frac12\,(P+Q)$ corresponds to a geodesic arc.

\begin{figure}[htbp]
  \centering
  \includegraphics{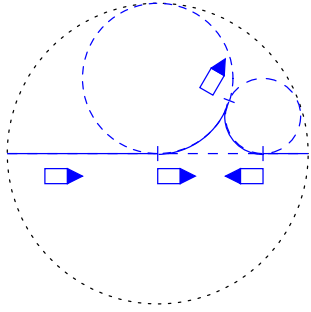}  
  \caption{The geometric meaning of the singular switching time}
  \label{singular-switching-time-pic}
\end{figure}

Some obscure identities and constants suddenly get a clear geometric
interpretation. For example, consider the singular switching time $w=\sqrt2$
and the identities 
$$
\begin{array}{lcl}
\mex{s_1T,wP,-wQ,-s_2T}&=&\mex{(s_1+s_2)T,wP,-wQ}\\
&=&\mex{wP,-wQ,-(s_1+s_2)T}.
\end{array}
\eqno{\rm(*)}
$$
Computational verification is a tedious exercise that gives no insight at all
why these identities hold true.
Now consider the path, say $\zeta(t)\in\bbD$, corresponding to $\mex{wP,-wQ}$.
Let $\zeta_0,v_0,\zeta_1,v_1$ denote the boundary data (positions and
tangents. Assume w.l.o.g. that $\zeta_0=0$ and $v_0=1$. Then $\zeta_1\in(0,1)$
(actually: $\zeta_1=\frac12\,\sqrt2$), and $v_1=-1$. So $\zeta(t)$ is tangent to the geodesic
through~$\zeta_0,\zeta_1$, cf. Figure~\ref{singular-switching-time-pic}.
So geometrically it is clear why these identities hold. We start and finish on
the same geodesic, but with reversed orientation (forward/backward motion).
The $\mex{wP,-wQ}$-maneuver accomplishes the \textit{U-turn} part, and it does
not matter, if we perform this turn at the beginning, in the end, or somewhere
inbetween. If we reflect Figure~\ref{singular-switching-time-pic} along the
horizontal axis, then it becomes obvious that we could have performed the \textit{turning
  maneuver} starting with a right turn (instead of a left turn). Hence we see that
$\mex{wP,-wQ}=\mex{wQ,-wP}$ holds true in $\PSL(2)$.
These observations are actually in full analogy to the euclidean case,
cf.~\cite[Fig.18, p.59]{suss-car}.
It also becomes evident that for (SSP)-extremals the path $\zeta(t)$ will
never intersect itself, so $\zeta(t_1)=\zeta(t_2)$ iff $t_1=t_2$ because
$\zeta(t)$ will consist of subarcs of a fixed geodesic and interspersed
turning maneuvers. W.l.o.g. one may assume that this geodesic is the diameter $(-1,1)$ as
in Fig.~\ref{singular-switching-time-pic}, and then it becomes obvious that
$\zeta(t)$ is doublepoint free.

\begin{Proposition}
 For \textsc{(SSP)}-extremals and \textsc{(CSP)}-extremals with $\leq4$ factors
 the path $\zeta(t)$ has no self-intersections
 {\rm(}i.e. $\zeta(t_1)=\zeta(t_2)$ implies $t_1=t_2${\rm)}.
\end{Proposition}
Instead of a rigorous proof (which would not be  elucidating at all) we simply
provide pictures showing why neither $C4a$-extremals ($2$ cusps) nor
$C4c$-extremals ($1$ cusp) have self-intersections.
\begin{figure}[htbp]
  \centering
  \begin{tabular}{ccc}
  \includegraphics{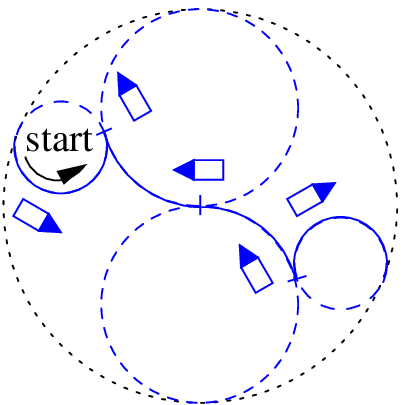} &\quad &
  \includegraphics{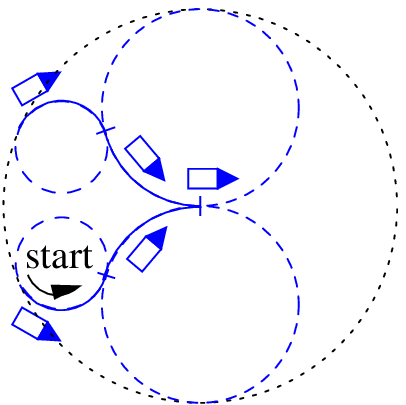} \\
  $C4a$-extremals & & $C4c$-extremals \\
  $\mex{rP,-sQ,-sP,tQ}$ &&   $\mex{rP,sQ,-sP,-tQ}$
  \end{tabular}
  \caption{$C4a$- and $C4c$-extremals have no self-intersections}
  \label{fig:csp4-injectivity-pic}
\end{figure}

The reason why we look at self-intersections is that the corresponding subpath
is nothing but a turning maneuver (same initial and terminal point) by some
angle~$\alpha$. Since the group $\exp(\R U)$ is the stabilizer of the point~$0\in\bbD$,
we see that it is impossible to reach a nontrivial rotation $\exp(\alpha U)$
from the identity~$\1$ along an extremal of
type (SSP), (C4a), or (C4c).

\begin{Proposition}
 For $\alpha\in(0,\pi/2)$ let $s(\alpha)=\sin(\alpha)$ and
 $r(\alpha)=\tan(\alpha/2)$. Then $\mex{rP,-sQ,rP}=\exp(\alpha U)$
 is optimal in $\PSL(2)$ as well as in any other group.
\end{Proposition}
\begin{Proof}
 We know that $g:=\exp(\alpha U)$ has an optimal factorization.
 Since we have a sufficient family, the previous proposition implies
 that the best factorization is of type~$A3$ or $A4$. In view of Table~\ref{sl2-family-tab}
 we must solve  $A3(r,s,t)=\sigma(g)$ and $A4(r,s,t)=\sigma(g)$ for all
 $\sigma\in\Gamma$. Since $\sigma_U(g)=g$, we only need to consider two different
 righthand sides, namely $g$ and $\sigma_T(g)=\exp(-\alpha U)$.
 We already observed $A3(r(\alpha),s(\alpha),r(\alpha))=g$, whereas
 $\sigma_T(g)=A3(r,s,t)$ implies $s=-\sin\alpha<0$, so this solution is not
 feasible (cf.~Table~\ref{sl2-family-tab}).
 In $\SL(2)$ the equation $A4(r,s,t)=\exp(\pm\alpha U)$ leads to:
$$
\begin{pmatrix}
      {*}      & r+s-rs^2 \\
  -s+s^2t-t  & 1-s^2
\end{pmatrix}
=
\begin{pmatrix}
  \cos(\alpha) & \pm\sin(\alpha) \\
  \mp\sin(\alpha) & \cos(\alpha)
\end{pmatrix},
$$
so $s=\sqrt{1-\cos(\alpha)}$.
 Since $\alpha\in(0,\pi/2)$ we deduce $s<1$.
But optimal $A4$-extremals have switching time~${}\geq1$.
Hence none of the $A4$-factorizations is optimal, whence
$g=\mex{rP,-sQ,rP}=\mex{-rQ,sP,-rQ}$ are \textbf{the}
optimal factorizations.
\end{Proof}


\section*{Appendix C: More details for \boldmath{$\SLt$}}

\renewcommand{\thesection}{C}

It is quite natural to ask for a characterization of optimal arcs in other groups
(than $\SL(2)$ and $\PSL(2)$), and in particular in~$\SLt$.
First it is important to notice that there is no
apriori bound for the number of factors of an optimal extremal in $\SLt$.
The reason is very simple: we could never reach all points if we considered only
(ALT)-extremals and (SSP)-extremals with a fixed number, say, $N$, of factors.

\begin{Proposition}
  Let $\cM=\exp(\R\conv(P,Q))\subseteq\SLt$.
  then 
$$
\cM^n=\underbrace{\cM\cdots\cM}_{n}\subseteq
\cE+(2n-1)\left(-\frac\pi2,\frac\pi2\right)U
\qmb{for all $n\in\N$.}
$$
In particular $\cM^n\neq\SLt$ for all $n\in\N$, and there is no apriori bound
on the number of factors of an optimal extremal.
\end{Proposition}
\begin{Proof}
  For $n=1$ the claim holds true, cf. the description of one-parameter groups.
  Now assume the claim holds true for some~$n\in\N$.
  Let $g_0\in\cE+(2n-1)\left(-\frac\pi2,\frac\pi2\right)U$ and
  $g_1\in\cE+(-\frac\pi2,\frac\pi2)U$.
  We recall that
  $\cE\circ\cE\subseteq\cE+(-\frac\pi2,\frac\pi2)U$.
  Writing $g_j=X_j+u_jU$ with $X_j\in\cE$ and $u_j\in\R$ our
  multiplication formula for $\SLt$ yields:
$$
g_0\circ g_1 = (u_0 U)\circ  g' \circ (u_1 U)
$$
with $g'= X'+u'U\in\cE\circ\cE\subseteq\cE+(-\frac\pi2,\frac\pi2)U$.
Hence $g_0\circ g_1 = X''+u'' U$ where 
$u''=u_0+u_1+u'$.
Since $|u_0|\leq(2n-1)\frac\pi2$ and $|u_1|,|\tilde u|\leq\frac\pi2$,
we obtain $|u''|\leq (2n+1)\frac\pi2$, hence
$\cM^{n+1}\subseteq\cE+(2n+1)\left(-\frac\pi2,\frac\pi2\right)U$ follows.
This proves the first claim from which we immediately deduce the second claim.
\end{Proof}

Nevertheless, for finite coverings of $\PSL(2)$ we get an apriori bound on the
numbers of factor of optimal (SSP)-extremals. 
We will also conjecture a bound for (ALT)-extremals and 
present some evidence why this conjecture should be true.
In any case the crucial parameter to consider is the cardinality of the center
$|Z(G)|$ where $G$ denotes an arbitrary covering of~$\PSL(2)$.


\subsection*{Singular extremals}

Reconsidering our results for $\PSL(2)$ and $\SL(2)$ we notice that
we obtained the bounds~$5$, resp., $7$ if we count all virtual switches.
This is actually the appropriate thing to do. Our goal is to show
that for $|Z(G)|=N$ we obtain the bound $2N+3$.
The first step is to see how Proposition~\ref{ssp-eqns-sl2} generalizes to $\SLt$.
Except for Equation~(\ref{sl2:eq-1}) nothing changes.

\begin{Proposition}\label{ssp-eqns-slt}
  Let $w=\sqrt{2}$. Then the following equalitites hold true:
  \begin{eqnarray}
    \mex{s P,-\frac2s\, Q} & = & \mex{\frac2s\, Q,-s P,\pi U},\quad(\forall s>0).\label{eq-1}
\\
    \mex{w P,-w Q} & = & \mex{w Q,-w P,\pi U}.\label{eq-1w}
 \\
    \mex{s T,w P,-w Q} & = & \mex{w P,-w Q,-s T}\quad\mbox{$s\in\R$.}
\label{eq-2}
\\
\label{eq-3}
    \mex{w P,-w Q,-w P,w Q} & = &  \exp(\arsinh(2w) T).
\\
\label{eq-4}
\mex{r P, w Q,-w P} & = &\mex{ \frac{2}{r+w}Q,\, -(r+w)P,\, -\frac{r\,w}{r+w} Q}.
\\
\label{eq-5}
    \mex{w P,-w Q,-w P} & = & \mex{\frac w2 Q,\,-2 wP,\,-\frac w2 Q}.
  \end{eqnarray}
\end{Proposition}
\begin{Proof} 
 From Eqn.~(\ref{sl2:eq-1}) we deduce
 $\mex{s P,-\frac2s Q}=\mex{\frac2s\, Q,-s P,\pi U}\circ(2\pi k)U$ for some $k\in\Z$.
 Since $s>0$, we deduce $\exp(sP), \exp(-2s^{-1}\,Q)\in\cE+(0,\pi/2)U$. Therefore
$$
\mex{s P,-2s^{-1}\, Q}\in\cE+\left(-\frac\pi2,\frac{3\pi}2\right)U,\quad
\mex{2s^{-1} Q,-s P}\in\cE+\left(-\frac{3\pi}2,\frac{\pi}2\right)U.
$$
Hence $\mex{2s^{-1}\, Q,-s P,\pi U}\in\cE+\left(-\frac\pi2,\frac{3\pi}2\right)U$,
whence $k=0$, and Eqn.~(\ref{eq-1}) follows for all  $s>0$.
\\[1mm]
Eqn.~(\ref{eq-2}) is equivalent to $e^{-w\ad(P)}e^{w\ad(Q)} T=-T$. As it
involves only the adjoint action, it holds true regardless of the group~$G$.
To prove Eqn.~(\ref{eq-3}) we observe 
$\mex{wP,-wQ,-wP}=\exp(-w e^{w\ad(P)}Q)\in\cE+(0,\pi/2)U$ and
$\exp(wQ)\in\cE+(-\pi/2,0)U$, hence $\mex{wP,-wQ,-wP,wQ}\in\cE+(-\pi,\pi)U$.
As an elementary computation in $\SL(2)$ yields
$$\cover(\mex{wP,-wQ,-wP,wQ})=\cover(\exp(\arsinh(2w)T)),
$$
we deduce $\mex{wP,-wQ,-wP,wQ}\in\cE\circ 2\pi\Z\,U$, hence~(\ref{eq-3})
follows.
\\[2mm]
Eqn.~(\ref{eq-4}) trivially holds true for $r=0$. As both sides are continuous
in~$r$, it must hold  for all $r\in(-\sqrt2,\infty)$. The special choice $r=w$
yields Eqn.~(\ref{eq-5}).
\end{Proof}

Note that Eqn.~(\ref{eq-1}) is the only equation that has changed (in
comparison to Proposition~\ref{ssp-eqns-sl2}.
In particular, propositions~\ref{ssp-3w-propo} and~\ref{ssp-top-bottom-bound}
hold true in~$\SLt$ because their proofs only required Eqn.~(\ref{eq-5}),
resp., eqns.~(\ref{eq-2},\ref{eq-4}).
The only instance where we
used the result corresponding to Eqn.~(\ref{eq-1})
was in the proof of Proposition~\ref{s7a-sl2-propo}.

\begin{Proposition}\label{ssp-bound-slt}
 Let $|Z(G)|=N\in\N$. Then the factorization
$$
\bbM\big(rP,\underbrace{wP,-wQ,-wQ,wP,\dots}_{\mbox{$2N$ factors}}\big)
$$
 is not optimal in $G$ for $r>0$.
 In particular, optimal (SSP)-extremals in $G$ have at most~$2N+3$ factors.
\end{Proposition}
\begin{Proof}
 Considering $\mex{rP,wP,-wQ,-wQ,wP,\dots}$ we let $\mu=r+w>0$ and observe that
$$
\mex{\mu P,-wQ}=\mex{\frac2\mu\,Q,-\mu P,-\varepsilon Q}\circ(\pi U).
\qmb{with}\ \varepsilon=w-\frac2\mu>0.
$$
Next we observe $\mex{wP,-wQ}=\mex{wQ,-wP}\circ(\pi U)$ and
$\mex{-wQ,wP}=\mex{-wP,wQ}\circ(\pi U)$.
As $|Z(G)|=N$ and $\pi U\in Z(G)$, $(\pi U)^N=\1$. Hence
\begin{eqnarray*}
\lefteqn{\mex{\mu P,-wQ}\circ
\underbrace{\mex{-wQ,wP}\circ\mex{wP,-wQ}\cdots}_{\mbox{$N-1$ factors}} }
\\
& = & \mex{\frac2\mu\,Q,-\mu P,-\varepsilon Q}\circ
\underbrace{\mex{-wP,wQ}\circ\mex{wQ,-wP}\cdots}_{\mbox{$N-1$ factors}}
\end{eqnarray*}
Both factorizations   have equal
cost~$r+2Nw$. But the second factorization  cannot be optimal because it
does not come from an extremal. The switching pattern $Q\vdash
-P\vdash-Q\vdash-P\dots$ is neither alternating nor circular.
For $r>0$ small  the switching time for the second arc
is $\mu=r+w\in(w,2w)$, so it cannot be an (SSP)-extremal either.
Hence the LHS factorization $\mex{\mu P,-w Q,-wQ,wP,\dots}$
cannot be optimal.
An (SSP)-extremal with $2N+4$ factors always contains a subarc of the above
form. 
Note that $2N+3$ arcs are possible, just consider $\mex{-rQ,wP,sT,wP,-wQ,\dots}$.
\end{Proof}




\subsection*{Alternating extremals}

 One may look for a general pattern behind the arguments we used to
 obtain bounds on the number of factors of optimal (ALT)-extremals in $\SL(2)$
 and $\PSL(2)$. And of course, there is one.
 It is very  natural to ask the following questions:
  \begin{quote}
    How can we reach the central elements of $\SLt$, i.e., $k\pi U$ for
    $k\in\N$?  And what is the fastest way to do so?
  \end{quote}
  An answer to the first question (how?) is relatively easy to find. Let
$$
g(s)=\mex{s P,-s Q}=\exP{s}\exQ{-s}=
\begin{pmatrix}
  1-s^2 & s \\ -s & 1
\end{pmatrix}\in\SL(2).
$$
Then $\trace(g(s))=2-s^2$. Thus for $s\in(0,2)$ $\trace(g(s))\in(-2,2)$. Hence
$$
\spec(g(s))=\setof{\cos\alpha \pm i\sin\alpha},\quad\mbox{with}\quad
\cos\alpha=1-\frac{s^2}2.
$$
Let $\alpha(s)=\arccos\left(1-\frac{s^2}2\right)\in(0,\pi)$. Then
$g(s)$ is actually conjugate to $\exp(\alpha(s)U)$, i.e., there exists a matrix
$V\in\SL(2)$ such that
$$
V^{-1}g(s)V =\exU{\alpha(s)}.
$$
For $s=0$ this is true, and for $s\in(0,2)$ it follows by continuity.
In $\SLt$ we therefore have
$$
\gamma(s):=\mex{sP,-sQ}=\mex{X,\alpha(s)U,-X}\qmb{for some $X\in\sL(2)$.}
$$
Next we look for solutions of $\gamma(s)^m=k\pi U$, $m,k\in\N$. As $k\pi U$ is
central, this is equivalent to $m\,\alpha(s)=k\pi$. Since
$\alpha(s)\in(0,\pi)$, solutions exist only for $m>k$:
$$
\alpha_{k,j}(s)=\frac k{k+j}\pi,\quad j\in\N.
$$
Since $\cos\alpha=1-\frac{s^2}2$, we find
$$
s^2=2(1-\cos\alpha)=2\cdot 2\sin^2\left(\frac\alpha2\right),\qmbq{hence}
s(\alpha)=2\sin\left(\frac\alpha2\right).
$$
Thus for $\alpha_{k,j}$ the appropriate switching time is $s_{k,j}=2\sin(\alpha_{k,j}/2)$.
The cost of this factorization is 
$$
2(k+j)\,s_{k,j}=2(k+j)\sin\left(\frac{k\pi}{2(k+j)}\right)
=k\pi\, \frac{\sin\left(\frac{k\pi}{2(k+j)}\right)}{\left(\frac{k\pi}{2(k+j)}\right)}.
$$
Since $x\mapsto\frac{\sin x}{x}$ is decreasing and nonnegative in $[0,\pi]$ and
$\frac{k\pi}{2(k+j)}\in[0,\frac\pi2]$ for $k,j\in\N$, we deduce that $j=1$
provides the best factorization of the form~$\gamma(s)^{k+j}$.
Let us write
\begin{eqnarray*}
\lefteqn{  A(n;r,s,t)  =  \underbrace{\bbM(rP,-sQ,,\dots)}_{\mbox{$n$ factors}}} \\
&=&
\begin{cases}
  \exp(rP)\mex{-sQ,sP}^{\frac n2-1}\exp(-tQ), & \mbox{$n$ even,}\\
  \exp(rP)\mex{-sQ,sP}^{\frac{n-1}2-1}\mex{sP,-tQ}, & \mbox{$n$ odd,}\\
\end{cases}
\end{eqnarray*}
and let $A_n(s)=A(n;s,s,s)$. Then we can prove:

\begin{Proposition}
  For $n\in\N$ let $s_n=2\cos\left(\frac\pi{n}\right)$. Then
$$
A_n(s_n)=(n-2)\frac\pi2\,U\qmb{for all $n\geq3$.}
$$
 Moreover, $A(n+1;t,s_n,s_n)$ is not optimal for $t>0$.
\end{Proposition}
\begin{Proof}
  First we observe that 
$1-\frac{s_n^2}2=-\cos(\frac2n\pi)=\cos\left(\frac{n-2}{n}\pi\right)$, hence
$$
\gamma(s_n)=\mex{X, \frac{n-2}{n}\pi\,U,-X}\qmb{for some $X\in\sL(2)$.}
$$
  Now we distinguish two cases. First, if $n=2k$ is even then we obtain
\begin{eqnarray*}
\gamma(s_n)^k  & = &\mex{X, \frac{n-2}{n}\pi\,U,-X}^k=\mex{X,
  \frac{k(2k-2)}{2k}\pi \,U,-X}
\\
& = &  (k-1)\pi U =\frac{n-2}2 \pi\,U.   
\end{eqnarray*}
 On the other, if $n=2k+1$ is odd, we let $\tilde\gamma=A_n(s_n)$.
 We observe that
 $\imath\sigma_H(\tilde\gamma)=\tilde\gamma$, so $\tilde\gamma=\tau T+u U\in\R T+\R U$.
 Moreover, 
\begin{eqnarray*}
  \tilde\gamma\sigma_U(\tilde\gamma)
&=&\gamma(s_{2k+1})^{2k+1}
   =\mex{X, (2k+1)\frac{2k-1}{2k+1}\pi\,U,-X}
\\
& = & (2k-1)\pi U=(n-2)\pi U.
\end{eqnarray*}
As $\sigma_U(\tau T+uU)=-\tau T+u U$, we compute
\begin{eqnarray*}
\cover(\tau T+uU)\cover(-\tau T+uU)  & = &
 \left(\cosh^2(\tau)\cos(2u)+\sinh^2(\tau)\right)\,\id
\\&&{} +\cosh^2(\tau)\sin(2u)\,U
-2\sinh(\tau) \sin(u)\, H
.
\end{eqnarray*}
Since $\cover((2k-1)\pi U)=-\id$, we deduce $\cos(2u)=-1$, hence
$2u\in\pi+2\pi\Z$, whence $u\in\frac\pi2+\pi\Z$. Thus $\sin(u)\neq0$,
so $\tau=0$ and $\tilde\gamma=uU$ for some $u\in\frac\pi2+\pi\Z$.
Thus $\sigma_U(\tilde\gamma)=\tilde\gamma$, and therefore
$\tilde\gamma\sigma_U(\tilde\gamma)=(n-2)\pi U$ implies
$u=\frac{n-2}2\pi$.

In order to prove that $A(n+1;t,s_{n},s_n)$ is not optimal, we observe
that $A_n(s_n)\in\R U$ and $\Fix(\sigma_U)=\R U$. 
Therefore $A_n(s_n)=\sigma_U(A_n(s_n))=\mex{-s_nQ,s_nP,\dots}$, and we obtain
\begin{eqnarray*}
 A(n+1;t,s_n,s_n)&=&\exp(tP)\mex{-s_nQ,s_nP,\dots}
=\exp(tP)\sigma_U(A_n(s_n))
\\ &=&\exp(t P)A_n(s_n)=\mex{(t+s_n)P,-s_n Q,s_nP,-s_nQ,\dots},
\end{eqnarray*}
and both factorizations have equal cost~$t+n\,s_n$. But since $t>0$, the RHS
does not come from an extremal, hence it cannot be optimal.
\end{Proof}

\begin{Remark}
  It is clear that $s_n$ is algebraic as it is twice the real part of a
  $2n$-th root of~$1\in\C$. For small $n$ we obtain:
  \begin{center}
    \begin{tabular}{|c||c|c|c|c|c|c|c|c|}\hline
      $n$   & $3$ & $4$ & $5$ & $6$ & $7$ & $8$ \\ \hline
      \rule{0pt}{16pt}$s_n$ 
& $1$ & $\sqrt2$ & $\frac{1+\sqrt5}2$ & $\sqrt3$ & $2\cos(\pi/7)$ & $\sqrt{2+\sqrt2}$
\\
 & & $1.41421$ & $1.61803$ & $ 1.73205$ & $1.80194$ & $1.87939$ \\ \hline
    \end{tabular}
  \end{center}
Algebraically, for $n=2k$ even let
$$
p_n(\xi)=
  \sum_{j=0}^{\left\lfloor\frac k2 \right\rfloor}{k\choose 2j} \xi^{k-2j}(\xi^2-1)^j
$$
and for $n=2k+1$ odd let
$$
p_n(\xi)=\sum_{j=0}^{2k}(-1)^j\xi^j
+\sum_{j=1}^{k} {n \choose 2j}\xi^{n-2j}(\xi-1)(\xi^2-1)^{j-1}.
$$
Then $p_n(\frac12\,s_n)=0$ for all $n\in\N$.
In fact for $n=2k$ even, $p_n$ is derived from $\Im(z^k-i)=0$ while
for $n=2k+1$ odd, $p_n$ is derived from
$\Re\frac{1+z^n}{1+z}=\sum_{j=0}^{2k}(-z)^j=0$,
for $z=\xi+i\sqrt{1-\xi^2}$.
For small $n$ we obtain
$$
\begin{array}{rclcrcl}
  p_3(s/2) & = & (s-1)^2, &&
  p_4(s/2) & = & \frac12 (s^2-2), \\[2pt]
  p_5(s/2) & = & (s^2-s-1)^2, &&
  p_6(s/2) & = & \frac12\,s(s^2-3),\\[2pt]
  p_7(s/2) & = & (s^3-s^2-2s+1)^2, &&
  p_8(s/2) & = & \frac12 (s^4-4s^2+2),\\[2pt]
  p_9(s/2) & = & (s-1)^2 (s^3-3s-1)^2, &&
  p_{10}(s/2) & = & \frac12\,s(s^4-5s^2+5).
\end{array}
$$
\end{Remark}

Now we can formulate our conjecture concerning optimal (ALT)-extremals
\begin{Conjecture}
  Let $n\in\N$, $n\geq3$. Then $A_n(s_n)$ is optimal in $\SLt$ and
  in any group with $|Z(G)|\geq n-2$.
 \\ 
  If $A(n;r,s,t)$ is optimal, then $s\geq s_{n-1}$. This is true for
  \textbf{any} group with Lie algebra $\sL(2)$.
\\
  If $|Z(G)|=N$, then optimal \textsc{(ALT)}-extremals have at most~$N+3$
  factors. If $A(N+3;r,s,t)$ is optimal, then $s\in[s_{N+2}, s_{2N+2}]$.
\end{Conjecture}
After our discussion on how to reach the central elements it should be clear
why we expect the first statement to hold true.
The second part of the conjecture is trivial for $n=3$ ($s\geq s_2=0$), and
it has already been proved for $n=4$ and $n=5$, cf.
Propositions~\ref{alt3-s3-propo},\ref{alt4-s4-propo}.
Finally, the last statement generealizes the arguments given in 
Propositions~\ref{alt-sl2-a5-propo}, and~\ref{alt-psl2-propo}.
We conjecture that  in general there exists
a function $\mu_N(s)$ with the following properties:
\begin{itemize}
\item For $g:=A(N+1;\mu_N(s),s,s)$ we have
  \begin{eqnarray*}
    g\circ g&=&N\pi U \qmb{if $N$ is odd,}\\
    g\circ \sigma_U(-g) & = & N\pi U\qmb{if $N$ is even.}
  \end{eqnarray*}
\item $\mu_N(s)$ is well-defined for $s>s_{N+1}$,
  $s<\mu_N(s)<2s$ for $s\in(s_{N+2},s_{2N+2})$, and
      $\mu_N(s)=s$ for $s=s_{2N+2}$, $\mu_N(s)=2s$ for $s=s_{N+2}$.
\end{itemize}
For small $N$ one can verify this explicitly, cf.~Table~\ref{tab:mu-conjecture}
\begin{table}
\begin{center}
  \begin{tabular}{|c|c|c|c|}\hline
 $N$ & $\mu_N(s)$ & $\mu_N(s)=2s$ & $\mu_N(s)=s$ \\ \hline
\rule{0pt}{24pt}$1$ & $\dst\frac2s$ & $1=s_3$ & $\sqrt2=s_4$ \\
\rule{0pt}{24pt}$2$ & $\dst\frac{2s}{s^2-1}$  & $\sqrt2=s_4$ & $\sqrt3=s_6$ \\[10pt] \hline
\rule{0pt}{24pt}$3$ & $\dst\frac{2(s^2-1)}{s(s^2-2)}$ & $\frac{1+\sqrt5}2=s_5$ &
                      $\sqrt{2+\sqrt2}=s_8$ \\
\rule{0pt}{24pt}$4$ & $\dst\frac{2s(s^2-2)}{s^4-3s^2+1}$ & $\sqrt3=s_6$ &
                      $\sqrt{\frac{5+\sqrt5}2}=s_{10}$
\\
\rule{0pt}{24pt}$5$ & $\dst\frac{2(s^4-3s^2+1)}{s(s^4-4s^2+3)} $ & $s_7$ & 
$\frac{1+\sqrt3}{\sqrt2}=s_{12}$ \\ \hline
  \end{tabular}
\end{center}
\caption{The function $\mu_N(s)$ for $N=1,\dots,5$.}\label{tab:mu-conjecture}
\end{table}
Provided we have such a function $\mu_N$ we obtain that $A(N+2;r,s,t)$ is not
optimal if $s>s_{N+2}$ and $r+t>\mu_N(s)$. Hence 
$$
\begin{array}{llcl}
s>s_{2N+2}, & r=s,\,t=0 & \implies & \mbox{$A_{N+1}(s)
$ not optimal.}\\
s=s_{2N+2}, & r=s,\, t>0 & \implies & \mbox{$A(N+2;s,s,t)$ not optimal.}\\
s\in(s_{N+2},s_{2N+2}), & r=t=s & \implies & \mbox{$A_{N+2}(s)$ not optimal.}
\end{array}
$$
Recalling our  conjecture that optimal $A(N+4;\cdot)$-extremals must have
switching time~$s\geq s_{N+3}$, it is clear why
we are convinced that $A(N+4;r,s,t)$ cannot be optimal if $r>0$ or $t>0$.

We conclude with one more 
\begin{Conjecture}
  The factorization  $A_n(s)$ is optimal in $\SLt$ for all $s\in[2,2\sqrt2]$, $n\geq3$.
\end{Conjecture}

\bibliographystyle{plain}


\end{document}